\documentclass[12pt,oneside
]{article}
\usepackage{amsmath}
\usepackage{amssymb}
\numberwithin{equation}{section}

\title{Integrable Lattices: Random Matrices and
Random Permutations\footnote{Expanded version of
lectures at MSRI, Berkeley, Spring 1998. To appear in
"Random Matrices and Their Applications" :
Mathematical Sciences Research Institute Publications
\#40, Cambridge University Press, 2001. } }

\author{
Pierre van Moerbeke\thanks{ Department of Mathematics,
Universit\'e de Louvain, 1348 Louvain-la-Neuve,
Belgium and Brandeis University, Waltham, Mass 02454,
USA. E-mail: vanmoerbeke@geom.ucl.ac.be and
@math.brandeis.edu. The support of a National Science
Foundation grant \# DMS-98-4-50790, a Nato, a FNRS and
a Francqui Foundation grant is gratefully
acknowledged.}}

 \date{March 25, 2000}

\catcode`\@=11
\let\c@equation=\relax
\newcounter{equation}[subsection]

\catcode`\@=12

\newcommand{\MAT}[1]{\left(\begin{array}{*#1c}}
\newcommand{\mat}{\end{array}\right)}
\newcommand{\qed}{\leavevmode\unskip\nobreak\penalty200\hskip2pt\null
\nobreak\hfill\rule{1.1ex}{1.1ex}
\medbreak
}

\newcommand{\rg}{\rightarrow}

\newcommand{\DR}{{\cal D}}

\newcommand{\HR}{{\cal H}}

\newcommand{\LR}{{\cal L}}

\newcommand{\SR}{{\cal S}}
\newcommand{\VR}{{\cal V}}

\newcommand{\BC}{{\mathbb C}}

\newcommand{\BX}{{\mathbb X}}
\newcommand{\BY}{{\mathbb Y}}
\newcommand{\BZ}{{\mathbb Z}}

\newcommand{\iy}{\infty}
\newcommand{\pl}{\partial}
\newcommand{\al}{\alpha}
\newcommand{\vr}{\varepsilon}

\newenvironment{proof}{\medskip\noindent{\it Proof:\/} }{\qed}
\newenvironment
        {remark}{\medskip\noindent\underline{\it Remark:\/} }{\medbreak}
\newenvironment
        {example}{\medskip\noindent\underline{\it Example:\/} }{\medbreak}

\newcommand{\om}{\omega}
\newcommand{\vp}{\varphi}
\newcommand{\la}{\langle}
\newcommand{\ra}{\rangle}
\newcommand{\ga}{\gamma}

\newcommand{\dt}{\delta}
\newcommand{\Dt}{\Delta}
\newcommand{\sg}{\sigma}
\newcommand{\BR}{{\mathbb R}}
\newcommand{\lb}{\lambda}
\newcommand{\Lb}{\Lambda}
\newcommand{\tr}{\mbox{tr}}

\newcommand{\BJ}{{\mathbb J}}

\newcommand{\diag}{\operatorname{diag}}

\def\be#1\ee{\begin{equation}#1\end{equation}}
\def\bea#1\eea{\begin{eqnarray}#1\end{eqnarray}}
\def\bean#1\eean{\begin{eqnarray*}#1\end{eqnarray*}}

\newcommand{\Tr}{\operatorname{\rm Tr}}

\newtheorem{definition}{Definition}[section]
\newtheorem{theorem}[definition]{Theorem}
\newtheorem{lemma}[definition]{Lemma}
\newtheorem{corollary}[definition]{Corollary}
\newtheorem{proposition}[definition]{Proposition}

\makeatletter
\def\ps@X{\let\@mkboth\@gobbletwo
        \def\@oddhead{\tt 
        \hfil\S\thesection, p.\thepage
        }
        \def\@oddfoot{\rm\hfil\thepage\hfil
        }
        \let\@evenhead\@oddhead
        \let\@evenfoot\@oddfoot}
\makeatother

\begin{document}
\maketitle

\begin{abstract}

These lectures present a survey of recent developments
in the area of random matrices (finite and infinite)
and random permutations. These probabilistic problems
suggest matrix integrals (or Fredholm determinants),
which arise very naturally as integrals over the
tangent space to symmetric spaces, as integrals over
groups and finally as integrals over symmetric spaces.
An important part of these lectures is devoted to
showing that these matrix integrals, upon apropriately
adding time-parameters, are natural tau-functions for
integrable lattices, like the Toda, Pfaff and Toeplitz
lattices, but also for integrable PDE's, like the KdV
equation. These matrix integrals or Fredholm
determinants also satisfy Virasoro constraints, which
combined with the integrable equations lead to
(partial) differential equations for the original
probabilities.

\end{abstract}

\newpage

\tableofcontents

\setcounter{section}{-1}

\section{introduction}

The purpose of these lectures is to give a survey of
recent interactions between statistical questions and
integrable theory. Two types of questions will be
tackled here:

(i) Consider a random ensemble of matrices, with
certain symmetry conditions to guarantee the reality
of the spectrum and subjected to a given statistics.
What is the probability that all its eigenvalues
belong to a given subset $E$ ? What happens, when the
size of the matrices gets very large ? The
probabilities here are functions of the boundary
points $c_i$ of $E$.

(ii) What is the statistics of the length of the
largest increasing sequence in a random permutation,
assuming each permutation is equally probable ? Here,
one considers generating functions (over the size of
the permutations) for the probability distributions,
depending on the variable $x$.

 The main emphasis of these lectures is to show that integrable
theory serves as a useful tool for finding equations
satisfied by these functions of $x$, and conversely
the probabilities point the way to new integrable
systems.

These questions are all related to integrals over
spaces of matrices. Such spaces can be classical Lie
groups or algebras, symmetric spaces or their tangent
spaces. In infinite-dimensional situations, the
"$\iy$-fold" integrals get replaced by Fredholm
determinants.

 During the last decade, astonishing
discoveries have been made in a variety of directions.
A first striking feature is that these probabilities
are all related to Painlev\'e equations or interesting
generalizations. In this way, new and unusual
distributions have entered the statistical world.

Another feature is that each of these problems is
related to some integrable hierarchy. Indeed, by
inserting an infinite set of time variables
$t_1,t_2,t_3,...$ in the integrals or Fredholm
determinants - e.g., by introducing appropriate
exponentials $e^{\sum_1^{\iy}t_iy^i}$ in the integral
- this probability, as a function of
$t_1,t_2,t_3,...$, satisfies an integrable hierarchy.
Korteweg-de Vries, KP, Toda lattice equations are only
a few examples of such integrable equations.

Typically integrable systems can be viewed as
isospectral deformations of differential or difference
operators $\LR$. Perhaps, one of the most startling
discoveries of integrable theory is that $\LR$ can be
expressed in terms of a single``$\tau$-function"
$\tau(t_1,t_2,...)$ (or vector of $\tau$-functions),
which satisfy an infinite set of non-linear equations,
encapsulated in a single ``{\sl bilinear identity}".
The $t_i$ account for the commuting flows of this
integrable hierarchy. In this way, many interesting
classical functions live under the same hat:
characters of representations, $\Theta$-functions of
algebraic geometry, hypergeometric functions, certain
integrals over classical Lie algebras or groups,
Fredholm determinants, arising in statistical
mechanics, in scattering and random matrix theory!
They are all special instances of ``{\em
$\tau$-functions}".

The point is that the probabilities or generating
functions above, as functions of $t_1,t_2,...$ (after
some minor renormalization) are precisely such
$\tau$-functions for the corresponding integrable
hierarchy and thus automatically satisfy a large set
of equations.

These probabilities are very special $\tau$-functions:
they happen to be a solution of yet another hierarchy
of (linear) equations in the variables $t_i$ and the
boundary points $c_i$, namely $\BJ^{(2)}_k \tau
(t;c)=0$, where the $\BJ^{(2)}_k$ form -roughly
speaking- a Virasoro-like algebra: $$ \left[~
\BJ_k^{(2)},~ \BJ_{\ell}^{(2)} \right] =(k-\ell)~
 \BJ_{k+\ell}^{(2)} +...
  $$
The point is that each integrable hierarchy has a
natural ``{\em vertex operator}", which automatically
leads to a natural Virasoro algebra.
  Then, eliminating the partial derivatives in $t$
  from the two hierarchy of equations, the integrable
  and the Virasoro hierarchies, and finally setting $t=0$, lead
to PDE's or ODE's satisfied by the probabilities.

In the table below, we give an overview of the
different problems, discussed in this lecture, the
relevant integrals in the second column and the
different hierarchies satisfied by the integrals. To
fix notation, ${\cal H}_{\ell}$, ${\cal S}_{\ell},$
${\cal T}_{\ell}$ refer to the Hermitian, symmetric
and symplectic ensembles, populated respectively by
$\ell \times \ell$ Hermitian matrices , symmetric
matrices and self-dual Hermitian matrices, with
quaternionic entries. ${\cal H}_{\ell}(E)$, ${\cal
S}_{\ell}(E),$ ${\cal T}_{\ell}(E)$ are the
corresponding set of matrices, with all spectral
points belonging to $E$. $U(\ell)$ and $O(\ell)$ are
the unitary and orthogonal groups respectively. In the
table below, $V_t(z):=V_0(z)+\sum t_iz^i$, where
$V_0(z)$ stands for the unperturbed problem; in the
last integral $\tilde V_t(z)$ is a more complicated
function of $t_1,t_2,...$ and $z$, to be specified
later.


\vspace{1cm}

{\small 
\begin{tabular}{l|l|l}
 Probability problem
 &underlying $t$-perturbed   & corresponding \\
 &integral, $\tau$-function of $\longrightarrow$  &integrable\\ &&hierarchies\\

 \hline \\  & & \\

 $P(M\in {\cal H}_n(E)) $ & $\int_{{\cal H}_n(E)
 }e^{Tr(- V( M)+\sum_1^{\iy} t_iM^i)}dM$ &Toda lattice
  \\
&   &  KP hierarchy    \\  & & \\

 $P(M\in {\cal S}_n(E)) $ & $
  \int_{{\cal S}_n(E) }e^{Tr(- V( M)+\sum_1^{\iy} t_iM^i)}dM$ & Pfaff lattice
 \\
 &   & Pfaff-KP hierarchy     \\   & & \\

 $P(M\in {\cal T}_n(E)) $ &
  $\int_{{\cal T}_n(E)}e^{Tr(- V( M)+\sum_1^{\iy} t_iM^i)}dM$& Pfaff lattice
 \\
 &   &Pfaff-KP hierarchy     \\  & & \\

$P((M_1,M_2)\in $ & $\int_{{\cal H}^2_{n}(E)}
dM_1dM_2$
 &
2d-Toda lattice\\ \vspace{-.3cm}
&
 & KP-hierarchy
\\ \vspace{0cm}
 $\hspace{1cm} {\cal H}_n(E_1)\times
{\cal H}_n(E_2))$
  &${\displaystyle e^{-\Tr ( V_t(M_1)-
V_s(M_2)-cM_1M_2)}} $  &  \\  &&\\

$P(M\in {\cal H}_{\iy}(E)) $
  &
 $\det \left(I- K_t(y,z)I_{E^c}(z)\right)  $
  & KdV equation
  \\
 &( Fredholm determinant)
&
\\
&&\\

 longest increasing sequence&
 $ \int_{U(\ell )}e^{Tr\sum_1^{\iy}(t_iM^i-s_i\bar M^i)}dM$
 & Toeplitz lattice   \\

in random permutations && 2d-Toda lattice\\

 &&\\  &&\\

longest increasing sequence&
 $ \int_{O(\ell )}e^{Tr(xM+\tilde V_t(M))}dM$
 & Toda lattice   \\

in random involutions && KP-hierarchy\\

 &&\\  &&\\

\end{tabular}
}

\vspace{1cm}

\noindent {\bf Acknowledgment}: These lectures
represent joint work especially with (but also
inspired by) Mark Adler, Taka Shiota and Emil Horozov.
Thanks also for many informative discussions with
Jinho Baik, Pavel Bleher, Edward Frenkel, Alberto
Gr\"unbaum, Alexander Its, espacially Craig Tracy and
Harold Widom, and with other participants in the
semester at MSRI. I wish to thank Pavel Bleher, David
Eisenbud and Alexander Its for organizing a truly
stimulating and enjoyable semester at MSRI.

\section{Matrix integrals, random matrices and
permutations}

\subsection{Tangent space to symmetric spaces and associated random matrix
ensembles }

Random matrices provided a model for excitation
spectra of heavy nuclei at high excitations (Wigner
\cite{Wigner}, Dyson \cite{Dyson} and Mehta
\cite{Mehta}), based on the nuclear experimental data
by Porter and Rosenzweig \cite{Porter}; they observed
that the occurrence of two levels, close to each
other, is a rare event (level repulsion), showing that
the spacing is not Poissonian, as one might expect
from a naive point of view.

Random matrix ideas play an increasingly prominent
role in mathematics: not only have they come up in the
spacings of the zeroes of the Riemann zeta function,
but their relevance has been observed in the chaotic
Sinai billiard and, more generally, in chaotic
geodesic flows. Chaos seems to lead to the ``spectral
rigidity", typical of  the spectral distributions of
random matrices,  whereas the spectrum of an
integrable system is random (Poisson)! (e.g., see
Odlyzko \cite{Odlyzko} and Sarnak \cite{Sarnak2}).

All these problems have led to three very natural
random matrix ensembles: Hermitian, symmetric and
symplectic ensembles. The purpose of this section is
to show that these three examples appear very
naturally as tangent spaces to symmetric spaces.

A symmetric space $G/K$ is given by a semi-simple Lie
group $G$ and a Lie group involution $\sigma:
G\rightarrow G$ such that $$ K=\{ x\in G,~\sigma
(x)=x\}.$$ Then the following identification holds:
 $$G/K \cong \{ g \sigma(g)^{-1} ~\mbox{with}~ g\in
 G\},
 $$ and the involution $\sigma$ induces a map of the
 Lie algebra
 $$\sigma_*:{\frak g}\longrightarrow {\frak
 g}~~\mbox{such that }~(\sigma_*)^2=1,
 $$
with $$ {\frak g}={\frak k}\oplus {\frak
p}~\mbox{with}
 \left\{\begin{array}{l} {\frak k}=\{ a\in {\frak g}
  ~\mbox{such that} ~ \sigma_*(a)=a\}  \\
 {\frak p} =\{ a\in {\frak g}
  ~\mbox{such that} ~ \sigma_*(a)=-a\}
\end{array}
\right. $$ and
 $$
  [ {\frak k},{\frak k}] \subset {\frak k},~
  [ {\frak k},{\frak p}] \subset {\frak p},~
  [ {\frak p},{\frak p}] \subset {\frak k}.
  $$
  Then $K$ acts on ${\frak p}$ by conjugation: $ k{\frak
  p}k^{-1} \subset {\frak p}$ for all $k \in K$ and
  ${\frak p}$ is the tangent space to $G/K$ at the identity.
 The action of $K$ on ${\frak p}$ induces a root space
 decomposition, with ${\frak a}$ being a maximal abelian
 subalgebra in ${\frak p}$:
 $$
 {\frak p}={\frak a}+\sum_{\alpha \in \Delta}{\frak
 p}_{\alpha},~~\mbox{with}~m_{\alpha}=\dim {\frak
 p}_{\alpha}.
 $$
 Then, according to Helgason \cite{Helgason1}, the volume element on ${\frak p}$
is given by
 $$
 dV=\left(  \prod_{\alpha \in
 \Delta_+}\alpha(z)^{m_{\alpha}}\right)dz_1...dz_n
 ,$$
 where $\Delta_+$ is the set of positive roots; see \cite{Helgason1,
 Helgason2, Terras,Terng}. This will subsequently be worked out for the
 three so-called $A_{n}$-symmetric spaces. See also
 Sarnak's MSRI-lecture \cite{Sarnak} in these proceedings, who deals
 with more general symmetric spaces. I like to thank
 Chuu-Lian Terng for very helpful conversations
 on these matters.

\bigbreak

\noindent\underline{\it Examples}\,:

\subsubsection*{(i) Hermitian ensemble} Consider the
{\em non-compact symmetric space}\footnote{The corresponding
compact symmetric space is given by $(SU(n)\times SU(n))/SU(n)$.}
$SL(n,\BC)/SU(n)$ with $\sigma (g)=
 \bar g^{\top -1}$. Then \bean SL(n,\BC)/SU(n)&=
 &\{g\bar g^{\top} ~|~g \in SL(n,\BC)\}\\
  &=&\{\mbox{positive definite matrices with $\det =1$}\}
  \eean
  with $$K=\{ g \in SL(n,\BC) ~|~ \sigma (g)= g\}=\{ g \in SL(n,\BC) ~|~
g^{-1}=\bar g^{\top}\}=
  SU(n).$$
   Then $\sigma_*(a)=-\bar a ^{\top}$ and the tangent space to $G/K$ is
then given by the space
   ${\frak p}={\cal H}_n$ of Hermitian matrices $$sl(n,\BC)
    = {\frak k}
   \oplus {\frak p}=su(n)\oplus {\cal
  H}_n,~~\mbox{ i.e.},~a=a_1+a_2, ~a_1 \in su(n),~a_2\in
  {\cal H}_n. $$
  If $M\in {\cal H}_n$, then the $M_{ii}$, $\Re M_{ij}$ and $\Im M_{ij}$
  ($1\leq i<j\leq n$) are free variables, so that Haar measure on
$M \in {\cal H}_n$ takes on the following form:
\be
dM:=\displaystyle{\prod_1^n dM_{ii}\prod_{1\leq
i<j\leq n}}(d\Re M_{ij}\,d\Im
M_{ij}). \ee

A maximal abelian subalgebra $a\subset {\frak p}={\cal H}_n$ is given by
real diagonal matrices $z=$ diag$(z_{1},\ldots,z_n)$. Each $M\in{\frak
p}={\cal H}_n$ can be written as
$$
M=e^Az\,e^{-A},\quad e^A\in K=SU(n),
$$
with\footnote{$e_{k \ell}$ is the $n\times n$ matrix with all zeroes,
except for
1 at the $(k,\ell)$th entry.}
\be
A=\sum_{1\leq k\leq\ell\leq
n}\left(a_{k\ell}(e_{k\ell}-e_{\ell
k})+ib_{k\ell}(e_{k\ell}+e_{\ell k}) \right)\in {\frak
k}=su(n),~~a_{\ell\ell}=0. \ee
 Notice that
$e_{k\ell}-e_{\ell k}$ and $i(e_{k\ell}+e_{\ell k})
\in{\frak k} =su(n)$ and that $$ [e_{k \ell}-e_{\ell
k},z]=(z_{\ell}-z_{k})(e_{k\ell}+e_{\ell k}) \in{\frak
p}={\cal H}_{n}$$
\be
[i(e_{k \ell}+e_{\ell k}),z]=(z_{\ell}-z_{k})i(e_{k \ell}-e_{\ell k})
\in{\frak p}={\cal H}_{n}.
\ee
Incidentally, this implies that $e_{k \ell}+e_{\ell k}$ and
$i(e_{k \ell}+e_{\ell k})$ are two-dimensional
eigenspaces\footnote{$ad\,x(y):=[x,y]$.} of $(ad\,z)^2$ with eigenvalue
$(z_{\ell}-z_k)^2$. From (1.1.2) and (1.1.3) it follows
that
\be
[A,z]=(z_{\ell}-z_{k})\sum_{1\leq k<\ell\leq
n}\left(a_{k \ell}(e_{k \ell}+e_{\ell k})+ib_{k
\ell}(e_{k \ell}-e_{\ell k}) \right)\in{\frak p}={\cal
H}_{n} \ee and thus, for small $A$, we
have\footnote{$\Delta_{n}(z)=\displaystyle{\prod_{1\leq
i<j\leq n}(z_{i}-z_{j})}$ is the Vandermonde
determinant.} \bea dM&=&d(e^Az\,e^{-A})\nonumber\\
&=&d(z+[A,z]+\ldots) \nonumber\\
&=&\prod^n_{1}dz_{i}\prod_{1\leq k<\ell\leq
n}d((z_{\ell}-z_{k})a_{k\ell})d((z_{\ell}-z_{k})b_{k\ell}),
~ \mbox{using (1.1.4) and (1.1.1)}\nonumber\\
&=&\prod^n_{1}dz_{i}\Delta^2_{n}(z)\prod_{1\leq
k<\ell\leq n}da_{k\ell}db_{k\ell}. \eea Therefore
$\Delta^2(z)$ is also the Jacobian determinant of the
map $M\rg(z,U)$, such that $M=UzU^{-1}\in {\cal H}_n$,
and thus $dM$ admits the decomposition in polar
coordinates:
\be
dM=\Delta^2_{n}(z)dz_{1}\ldots dz_{n}dU,\quad U\in
SU(n). \ee
 In random matrix theory,  $\HR_n$ is
endowed with the following probability,
 \be
 P(M\in dM)=c_n e^{-tr V(M)}dM, ~~\rho(dz)=e^{-V(z)}dz ,\ee
 where $dM$ is Haar measure (1.1.6) on $\HR_n$ and $c_n$ is
 the normalizing factor. Since $dM$ as in (1.1.6) contains
 $dU$ and since the probability measure (1.1.7) only depends on the trace  of
 $V(M)$, $dU$ completely integrates out.
 Given $E\subset \BR$, define
 \be {\cal
 H}_n(E):=\{M \in {\cal H}_n \mbox{with all spectral
 points} \in E\subset\BR\} \subset{\cal H}_n .\ee
 Then
 \be P(M \in {\cal
 H}_n(E))=\int_{{\cal H}_{n}(E)} c_ne^{-\Tr V(M)}dM=\frac{\int_{E^n}\Delta^2
(z)\prod^n_1\rho(dz_k)}{\int_{\BR^n}\Delta^2
(z)\prod^n_1\rho(dz_k)}.\ee

 As
the reader can find out from the excellent book by
Mehta \cite{Mehta}, it is well known that, if the
probability $P(M\in dM)$ satisfies the following two
requirements: (i) invariance under conjugation by
unitary transformations $M\mapsto UMU^{-1}$, (ii) the
random variables $M_{ii}$, $\Re M_{ij}$, $ \Im M_{ij}$,
$1\leq i<j\leq n$ are independent, then $V(z)$ is
quadratic (Gaussian ensemble).

\subsubsection*{(ii) Symmetric ensemble} Here we
consider the {\em non-compact symmetric
space}\footnote{The compact version is given by $SU(n)/SO(n)$.}
  $SL(n,\BR)/SO(n)$ with $\sigma (g)= g^{\top -1}$
   .
  Then \bean SL(n,\BR)/SO(n)&=&\{gg^{\top} ~|~g \in SL(n,\BR)\}\\
  &=&\{\mbox{positive definite matrices with $\det =1$}\}
  \eean
  with $$K=\{ g \in SL(n,\BR) ~|~ \sigma (g)= g\}=\{ g \in SL(n,\BR) ~|~ g^{\top}=g^{-1}\}= SO(n).$$
   Then $\sigma_*(a)=-a^{\top}$ and the tangent space to $G/K$ is then
given by the space
   ${\frak p}={\cal S}_n$ of symmetric matrices, appearing
   in the decomposition of $sl(n,\BR)$,
   $$sl(n,\BR)= {\frak k}
   \oplus {\frak p}=so(n)\oplus {\cal
  S}_n,~~\mbox{ i.e.},~a=a_1+a_2, ~a_1 \in so(n),~a_2\in
  {\cal S}_n $$
  with Haar measure $dM=\displaystyle{\prod_{1\leq i\leq j\leq
  n}}dM_{ij}$ on ${\cal S}_n$.

  A maximal abelian subalgebra $a\subset{\frak p}={\cal S}_n$ is
  given by real traceless diagonal matrices $z=\diag(z_1,\ldots,z_n)$.
  Each $M\in{\frak p}={\cal S}_n$ conjugates to a diagonal matrix $z$
  $$
  M=e^Az\,e^{-A},\quad e^A\in K=SO(n),\quad A\in so(n).
  $$
  A calculation, analogous to example (i)(1.1.5) leads to
$$ dM=\bigl|\Delta _n(z)\bigr| dz _1\ldots dz
_ndU,\quad U\in SO(n). $$
 Random matrix theory deals with
  the following probability on $\SR_n$:
 \be
 P(M\in dM)=c_n e^{-tr V(M)}dM, ~~\rho(dz)=e^{-V(z)}dz ,\ee
 with normalizing factor $c_n$. Setting as in (1.1.8):
${\cal S}_n(E)\subset {\cal S}_n$ is the subset of matrices with spectrum
$\in E$. Then
\be
P(M\in {\cal S}_n(E))=\int_{{\cal S}_{n}(E)} c_ne^{-\Tr V(M)}dM
=\frac{\int_{E^n}|\Delta
(z)|\prod^n_1\rho(dz_k)}{\int_{\BR^n}|\Delta
(z)|\prod^n_1\rho(dz_k)}.\ee
 As in the Hermitian case, $P(M\in dM)$ is Gaussian, if $P(M\in dM)$ satisfies

\noindent{\em (i)} invariance under conjugation by orthogonal conjugation
$M\rg OMO^{-1}$,
\noindent{\em (ii)} $M_{ii},M_{ij}$ ($i<j$) are independent
random variables.

 \subsubsection*{(iii) Symplectic ensemble} Consider
 the {\em non-compact symmetric space}\footnote{The
 corresponding
 compact symmetric space is $SU(2n)/Sp(n)$.}
  $SU^*(2n)/USp(n)$ with $\sigma
  (g)=J g^{\top -1} J^{-1}$, where $J$ is the $2n \times 2n$
  matrix:
  \be J:=\left(
\begin{array}{cc@{}c@{}cc}
 &\boxed{\begin{array}{cc} 0 & 1 \\ -1 & 0 \end{array}} && & \\
 && \boxed{\begin{array}{cc} 0 & 1 \\ -1 & 0 \end{array}} &&\\
 &&& \boxed{\begin{array}{cc} 0 & 1 \\ -1 & 0 \end{array}} & \\
 &&&& \ddots
 \end{array}
 \right)\quad\mbox{with $J^2=-I$},
\ee and
 \bean
  G&=& SU^*(2n)=\{ g \in SL(2n,\BC)~|~
g=J\bar g J^{-1}\},\\ &&\\
  K&=&\{g \in SU^*(2n)~|~ \sigma (g)=g\}:= Sp(n,\BC)\cap U(2n)\\
&=& \{g\in SL(2n,\BC)~\bigl|~g^{\top}Jg=J\}\cap\{g\in
 SL(2n,\BC)~\bigl|~g^{-1}=\bar g^{\top}\}\\
 &=& \{g\in SL(2n,\BC)~\bigl|~g^{-1}=\bar g^{\top}\mbox{and}~
 g=J\bar g J^{-1}\}\\&=:&USp(n).\\
 \eean

Then, $\sigma_*(a)=-Ja^{\top}J^{-1}$ and
  \bean
 {\frak k}&=& \{ a \in su^*(2n)~\big |~
 \sigma_*(a)=a\}=sp(n,\BC) \cap u(2n)\\
 &=&\{ a \in \BC^{2n\times 2n}~\big |~ a^{\top}=-\bar
 a,~a=J\bar a J^{-1} \}\\ && \\
 {\frak p}&=& \{ a \in su^*(2n)~\big |~
 \sigma_*(a)=-a\}= su^*(2n) \cap iu(2n)\\
 &=&\{ a \in \BC^{2n\times 2n}~\big |~ a^{\top}=\bar
 a,~a=J\bar a J^{-1} \}\\
 &=& \left\{ M=(M_{k \ell})_{1\leq k ,\ell \leq n} ,
M_{k \ell}=\left(\begin{array}{ll} M^{(0)}_{k
\ell}&M^{(1)}_{k \ell}\\
 &  \\
 -\bar M^{(1)}_{k \ell}&\bar M^{(0)}_{k \ell}
\end{array}
\right)
        \mbox{ with }~M_{ \ell k}=\bar M_{k \ell}^{\top}
\in \BC^{2\times 2}\right\}\\
 &\cong&\{ \mbox{self-dual
 $n\times n$ Hermitean
matrices, with quaternionic entries}\}\\&=:& {\cal
T}_{2n}.
 \eean
The condition on the $2\times 2 $ matrices $M_{k
\ell}$ implies that $M_{kk}=M_k I$, with $M_k\in \BR$
and the $2\times 2$ identity $I$. Notice $USp(n)$ acts
naturally by conjugation on the tangent space ${\frak
p} $ to $G/K$. Haar measure on ${\cal T}_{2n}$ is
given by
\be
dM = \prod^n_1 dM_k
\prod_{1\leq k < \ell \leq n}
 dM_{k \ell}^{(0)} {d\bar M_{k \ell}^{(0)}}
   dM_{k \ell}^{(1)}  {d\bar M_{k \ell}^{(1)}},
   \ee
   since these $M_{ij}$ are the only free variables in the matrix $M\in{\cal
   T}_{2n}$. A maximal abelian subalgebra
in $\frak p$ is given by real diagonal matrices of the
form $z= \diag (z_1,z_1, z_2,z_2,...,z_n,z_n)$. Each
$M \in {\frak p}={\cal T}_{2n}$ can be written as
\be
M=e^A z e^{-A},~~ e^A\in K=USp(n), \ee with
$(a_{k\ell},b_{k\ell},c_{k\ell},d_{k\ell}\in\BR)$
 $$
A=\sum_{1 \leq k\leq \ell \leq n}
a_{k\ell}(e^{(0)}_{k\ell}-e^{(0)}_{\ell k}) +b_{k\ell}
(e^{(1)}_{k\ell} +e^{(1)}_{\ell k})
+c_{k\ell}(e^{(2)}_{k\ell} -e^{(2)}_{\ell k})
+d_{k\ell}( e^{(3)}_{k\ell} +e^{(3)}_{\ell k}) \in
{\frak k} $$
 \vspace{-1.2cm}
  \be \ee
  in terms of the four $2 \times 2$
matrices \footnote{$e_{k\ell}^{(i)}$ in (1.1.15)
refers to putting the $2\times 2$ matrix $e^{(i)} $ at
place $(k, \ell)$.}
 $$
e^{(0)}= \left(\begin{array}{ll} 1&0\\ 0&1
 \end{array}\right),~
 e^{(1)}=\left(\begin{array}{ll} i&0\\ 0&-i
 \end{array}\right),~
e^{(2)}= \left(\begin{array}{ll} 0&1\\ -1&0
 \end{array}\right),~
e^{(3)}= \left(\begin{array}{ll} 0&i\\ i&0
\end{array}\right).
 $$
 Since
  \bea
 \bigl[e^{(0)}_{k\ell}-e^{(0)}_{\ell k},z\bigr]& =&
 (z_{\ell}-z_k)(e^{(0)}_{k\ell}+e^{(0)}_{\ell
 k})\in {\frak p} \nonumber \\
 \bigl[ e^{(1)}_{k\ell} +e^{(1)}_{\ell k},z \bigr] &=&
(z_{\ell}-z_k)(e^{(1)}_{k\ell}-e^{(1)}_{\ell
 k})\in {\frak p}   \nonumber \\
\bigl[ e^{(2)}_{k\ell} -e^{(2)}_{\ell k},z \bigr] &=&
(z_{\ell}-z_k)(e^{(2)}_{k\ell}+e^{(2)}_{\ell
 k})\in {\frak p}  \nonumber\\
\bigl[ e^{(3)}_{k\ell} +e^{(3)}_{\ell k},z \bigr] &=&
(z_{\ell}-z_k)(e^{(3)}_{k\ell}-e^{(3)}_{\ell
 k})\in {\frak p},
 \eea
 $[A,z]\in {\frak p}$ has the following form: it has $2\times 2$
 zero blocks along the diagonal and from (1.1.16) and (1.1.15),
 \be
 ((k,\ell)\mbox{th block in }[A,z])=(z_{\ell}-z_k)
  \left(\begin{array}{ll} a_{k\ell}+ib_{k\ell}
   &c_{k\ell}+id_{k\ell}\\
   -c_{k\ell}+id_{k\ell}&a_{k\ell}-ib_{k\ell}
 \end{array}\right) ,(k<\ell).
\ee
Therefore, using (1.1.17), Haar measure $dM$ on ${\cal T}_{2n}$ equals
 \bean dM
&=&d(e^A z e^{-A})\\ &=&d(I + A + \ldots)z (I - A +\ldots)\\
&=&d(z + [A,z]+ \ldots)\\
  &=&\prod_{1\leq k \leq n} d z_k  \prod_{1
\leq k < \ell \leq n}
d((z_{\ell}-z_k)(a_{k\ell}+ib_{k\ell}))
d((z_{\ell}-z_k)(a_{k\ell}-ib_{k\ell}))\\
&&~~~~~~~~~~~ d((z_{\ell}-z_k)(c_{k\ell}+id_{k\ell}))
d((z_{\ell}-z_k)(-c_{k\ell}+id_{k\ell})) \\
 &=&  \Delta^4 (z)dz_1\cdots dz_n~ \prod_{1
\leq k < \ell \leq n}4 da_{k\ell}db_{k\ell}
 dc_{k\ell}dd_{k\ell}.
 \eean

 As before, define ${\cal
T}_{2n}(E)\subset {\cal T}_{2n} $ as the subset of matrices
with spectrum $\in E$ and define the probability:
 \be P(M\in{\cal
T}_{2n}(E))=\int_{{\cal T}_{2n}(E)} c_ne^{-\Tr V(M)}dM=
 \frac{\int_{E^n} \Delta^4 (z)
\prod^n_1\rho(dz_k)}{\int_{\BR^n} \Delta^4 (z)
\prod^n_1\rho(dz_k)}. \ee

\remark Notice ${\cal T}_{2n}$ is called the
symplectic ensemble, although the matrices in ${\frak
p}={\cal T}_{2n}$ are not at all symplectic; but
rather the matrices in ${\frak k}$ are.

\subsection{Infinite Hermitian matrix ensembles}

Consider now he limit of the probability
 \be P(  M\in {\cal H}_{n}(E))=\frac{\int_{E^n} \Delta^2 (z)
 \prod^n_1\rho(dz_k)}{\int_{\BR^n} \Delta^2 (z)
 \prod^n_1\rho(dz_k)},~ \mbox{when} ~~ n\nearrow
\infty. \ee
 Dyson \cite{Dyson} (see also Mehta
\cite{Mehta}) used the following trick, to circumvent
the problem of dealing with $\iy$-fold integrals.
Using the orthogonality of the {\em monic orthogonal
polynomials} $p_{k}=p_{k}(z)$ for the weight
$\rho(dz)$ on $\BR$, and the $L^2$-norms
$h_k=\int_{\BR} p_k^2(z)\rho(dz)$ of the $p_k$'s, one
finds, using $(\det A)^2=\det(AA^{\top})$,

\bigbreak

\noindent
$\displaystyle{\int_{\BR^n}\Delta^2(z)\prod_1^n\rho(dz_i)}$
\begin{eqnarray}
&=&\int_{\BR^n}\det( p_{i-1}(z_j))_{1\leq i,j\leq
n}\det( p_{k-1}(z_{\ell}))_{1\leq k,\ell\leq
n}\prod^n_{k=1}\rho(dz_k)\nonumber\\
&=&\sum_{\pi,\pi'\in\sg_n}(-1)^{\pi+\pi'}\prod^n_{k=1}
 \int_{\BR}
p_{\pi(k)-1}(z_k)
p_{\pi'(k)-1}(z_k)\rho(dz_k)\nonumber\\
&=&n!\prod^{n-1}_{0}\int_{\BR}
p^2_k(z)\rho(dz)=n!\prod^{n-1}_{0}h_k.
\end{eqnarray}
For the integral over an arbitrary subset $E\subset
\BR$, one stops at the second equality, since the
$p_n$'s are not necessarily orthogonal over $E$. This
leads to the probability (1.2.1),

\bigbreak

\noindent $P( M\in {\cal
H}_{n}(E)))$
\begin{eqnarray}
&=&\frac{1}{n!\prod_1^nh_{i-1}}\int_{E^n}\det\left(\sum_{1
\leq j\leq n}  p_{j-1}(z_k)
p_{j-1}(z_{\ell})\right)_{1\leq
 k,\ell\leq n}\prod^n_1\rho(dz_i)\nonumber\\
&=&\frac{1}{n!}\int_{E^n}\det(K_n(z_k,z_{\ell}))_{1\leq
k,\ell\leq n} \prod_1^n\rho(dz_i),
\end{eqnarray}
 in terms of the kernel
\be
K_n(y,z):=\sum^n_{j=1}\frac{p_{j-1}(y)}{\sqrt{h_{j-1}}}\frac{p_{j-1}(z)}{\sqrt{h
_{j-1}}}.
 \ee
 The orthonormality relations of the
$p_k(y)/\sqrt{h_k}$ lead to the reproducing property
for the kernel $K_n(y,z)$:
\be
\int_{\BR} K_n(y,z)K_n(z,u)\rho(dz)=K_n(y,u),
\quad\quad\int_{\BR} K_n(z,z)\rho(dz)=n. \ee Upon
replacing $E^n$ by $
\prod^k_1 dz_i\times \BR^{n-k}$ in (1.2.3), upon
integrating out all the remaining variables
$z_{k+1},...,z_n$ and using the reproducing property
(1.2.5), one finds the $n$-point correlation function
\bea &&\hspace{-2cm}P(\mbox{one eigenvalue in
each\,\,}[z_i,z_i+dz_i],~i=1,...,k)\nonumber \\
 &=&c_n\det\left(K_n(z_i,z_j)\right)_{1\leq i,j\leq
k}\prod_1^k\rho (dz_i). \eea Finally, by Poincar\'e's
formula for the probability $\displaystyle{P\bigl(\cup
E_i\bigr)}$, the probability that no spectral point of
$M$ belongs to $E$ is given by a Fredholm determinant
\bean P( M\in {\cal H}_{n}(E^c))&=& \det(I-\lb
K_n^E)\\ &=& 1+\sum^{\iy}_{k=1}(-\lb)^k\int_{z_1\leq
...\leq z_k}\det\Bigl(K_n^E(z_i,z_j)\Bigr)_{1\leq
i,j\leq k}\prod_1^k\rho(dz_i), \eean for the kernel $
K_n^E(y,z)=K_n(y,z)I_E(z)$.

\bigbreak

\noindent $\bullet$ {\em Wigner's semi-circle law}:
    For this ensemble (defined by a large class of
     $\rho$'s, in particular for the Gaussian ensemble)
      and for very large $n$, the
density of eigenvalues tends to Wigner's semi-circle
distribution on the interval $[-\sqrt{2n},\sqrt{2n}]$:
$$ \mbox{density of
eigenvalues}\left\{\begin{array}{ll}
=\displaystyle{\frac{1}{\pi}}\sqrt{2n-z^2}dz,& |z |
\leq\sqrt{2n}\\
 & \\
=0,&|z|>\sqrt{2n}.
\end{array}
\right. $$

\noindent $\bullet$ {\em Bulk scaling limit}: From the
formula above, it follows that the average number of
eigenvalues per unit length near $z=0$ (``the bulk")
is given by $\sqrt{2n}/\pi$ and thus the average
distance between two consecutive eigenvalues is given
by $\pi/\sqrt{2n}$. Upon using this rescaling, one
shows (\cite{Kamien,Mahoux,Nagao,Pastur,JMMS})
 $$ \lim_{n\nearrow\infty}\frac{\pi
}{\sqrt{2n}}K_{n}\left(\frac{\pi
x}{\sqrt{2n}},\frac{\pi y}{\sqrt{2n}}
\right)=\frac{\sin\pi(x-y)}{\pi (
x-y)}\quad\quad\mbox{(Sine kernel)} $$ and $$
P(\mbox{exactly $k$ eigenvalues\,\,}\in
[0,a])=\frac{(-1)^k}{k!}\left. \left(\frac{\pl}{\pl
\lambda}\right)^k \det(I-\lb
KI_{[0,a]})\right|_{\lambda =1}$$ with
 \be \det
(I-\lambda KI_{[0,a]})= \exp\int_{0}^{\pi
a}\frac{f(x;\lb)}{x}dx,
 \ee
  where $f(x,\lb)$ is a
solution to the following differential equation, due
to the pioneering work of Jimbo, Miwa, Mori, Sato
\cite{JMMS}, (${}^{\prime}=\pl / \pl x$)
  \bea (xf'')^2&=&4(xf'-f)(-f^{\prime 2}-xf'+f)
 ~,~\mbox{with\,\,}f(x;\lb)\cong
-\frac{\lb}{\pi}x \mbox{\,\,for\,\,} x\simeq 0.
\nonumber \\ && \hspace{7cm}\mbox{({\bf Painlev\'e
V})}  \eea

\noindent $\bullet$ {\em Edge scaling limit}: Near the
edge $\sqrt{2n}$ of the Wigner semi-circle, the
scaling is $\sqrt{2}n^{1/6}$ and thus the scaling is
more subtle: (see \cite{BB,Forrester,Moore,Mehta,TW1})
\be y=\sqrt{2n}+\frac{u}{\sqrt{2}n^{1/6}} ,\ee
 and so for the kernel $K_n$ as in (1.2.4), with the
 $p_n$'s being Hermite polynomials,
  $$
\lim_{n\nearrow\infty}\frac{1}{\sqrt{2}n^{1/6}}K_n\left(\sqrt{2n}+\frac{u}{\sqrt{2}n^{1/
6}}, \sqrt{2n}+\frac{v}{\sqrt{2}n^{1/6}}
\right)=K(u,v), $$ where $$
K(u,v)=\int^{\iy}_{0}A(x+u)A(x+v)dx,~~~A(u)=\int^{\iy}_{
-\iy}e^{iux-x^3/3}dx. $$ Relating $y$ and $u$ by
(1.2.9), the statistics of the largest eigenvalue for
very large $n$ is governed by the function,
   \bean
    P(\lb_{\max}\leq
y)&=&P\left(2n^{\frac{2}{3}}\Bigl(\frac{\lb_{\max}}{\sqrt{2n}}-1
\Bigr)\leq u\right),\mbox{  for $n\nearrow \iy$,}\\ 
 &=& \det
(I-KI_{(-\iy,u]})=\exp\left(-\int^{\iy}_{u}(\al-u)g^2(\al)
d\al \right), \eean
 with  $g(x)$ a solution of
 \be
\left\{\begin{array}{l} g''=xg+2g^3\\ g(x)\cong
 -\frac{
  e^{-\frac{2}{3}  x^{\frac{3}{2}}}}{2\sqrt \pi x^{1/4}}
\mbox{\,\,for\,\,}x\nearrow \infty.
\end{array}\right. (\mbox{{\bf Painlev\'e II}})
\ee
  The latter is essentially the asymptotics of the
Airy function. In section 5, I shall derive, via
Virasoro constraints, not only this result, due to
Tracy-Widom \cite{TW1}, but also a PDE for the
probability that the eigenvalues belong to several
intervals, due to Adler-Shiota-van Moerbeke
\cite{ASV1,ASV2}.


 \noindent $\bullet$ {\em Hard edge scaling limit}:
 Consider the ensemble of $n\times n$ random matrices
for the Laguerre probability distribution, thus
corresponding to
 (1.1.9) with $\rho(dz)=z^{\nu/2}e^{-z/2}dz$. One shows the density of
 eigenvalues near $z=0$ is given by $4n$ for
 very large $n$. At this
 edge, one computes for the kernel (1.2.4) with Laguerre polynomials
 $p_n$ \cite{Nagao, Forrester}:
 \be
  \lim_{n\nearrow\infty}\frac{1}{4n}K^{(\nu)}_n
   \left(\frac{u}{4n}, \frac{v}{4n}
\right)=K^{(\nu)}(u,v), \ee
  where $K^{(\nu)}(u,v)$ is the
{\em Bessel kernel}, with Bessel functions $J_{\nu}$,
\begin{eqnarray}
K^{(\nu)}(u,v)&=& {1\over2}\int_0^1
xJ_\nu(xu)J_\nu(xv)dx \nonumber \\
 &=& \frac{J_{\nu} (u)
  \sqrt{u}
  J'_{\nu} (v) -
J_{\nu} (\sqrt{v}) \sqrt{v} J'_{\nu} (\sqrt{u})} {2(u
- v)}.
\end{eqnarray}
 Then
$$ P(~\mbox{no eigenvalues}~ \in [0,x])=
 \exp \left(-\int_0^x \frac{f(u)}{u} du\right),
 $$
with $f$ satisfying
 \be
  (xf^{\prime\prime})^2-4\bigl(xf^{\prime}-f \bigr
)f^{\prime 2}+\bigl((x-\nu^2)f^{\prime}-f \bigr
)f^{\prime}=0 .\, \mbox{ \bf (Painlev\'{e} V)}
 \ee
This result due to Tracy-Widom \cite{TW11} and a more
general statement, due to \cite{ASV1,ASV2} will be
shown using Virasoro constraints in section 5.

\subsection{Integrals over classical groups}

 The integration on a compact
semi-simple simply connected Lie group $G$ is given by the
formula (see Helgason \cite{Helgason2})
\be
 \int_Gf(M)dM=\frac{1}{|W|}\int_T\left| \prod_{\alpha
\in \Delta}2 \sin \frac{\alpha (iH)}{2}\right| dt
\int_Uf(utu^{-1})du,~~~t=e^H,
 \ee
  where $A\subset G$ is a maximal subgroup,
  with ${\frak g}$ and ${\frak a}$ being the Lie algebras
  of $G$ and $A$. Let $du $ and $dt$ be Haar
measures on $G$, $A$ respectively such that
$$\displaystyle{\int_{A}dt=\int_{U}du=1};$$
 $\Delta$
denotes the set of roots of ${\frak g}$ with respect to
${\frak a}$; $|W|$ is the order of the Weyl group of $G$.

{\em Integration formula (1.3.1) will be applied to
integrals of $f=e^{ \sum_1^{\iy} t_i \Tr M^i}$ over
the groups SO$(2n)$, SO$(2n+1)$ and Sp$(n)$.} Their
Lie algebras (over $\BC$) are given respectively by
$\frak d_n$, $\frak b_n$, $\frak c_n$, with sets of
roots: (e.g., see \cite{Bourbaki})
 $$
\Dt_n=\{\pm\vr e_i,1\leq i\leq
k,\pm(e_i+e_j),\pm(e_i-e_j),1\leq i<j\leq n\} ,$$ with

\hspace{3cm}$\vr=0$ for $\frak d_n=so(2n)$

\hspace{3cm}$\vr=1$ for $\frak b_n=so(2n+1)$

\hspace{3cm}$\vr=2$ for $\frak c_n=sp(n)$.

\bigbreak

 \noindent Setting $H=i\theta$, we have, in
view of formula (1.3.1),

\medbreak

$\displaystyle{ \left|\prod_{\al\in\Dt}2
 \sin\frac{\al(iH)}{2}\right|dt}$

 \bean&=&
\left\{\begin{array}{ll}c_n \left(\prod_{1\leq j<k\leq
n}\sin
\frac{\theta_j-\theta_k}{2}\sin\frac{\theta_j+\theta_k}{2}
\right)^2 \prod_{1}^n 
d\theta_j
   & \mbox{for ${\frak d}_n$}\\
 c_n \left(\prod_{1\leq j<k\leq n}\sin
\frac{\theta_j-\theta_k}{2}\sin\frac{\theta_j+\theta_k}{2}
\right)^2\prod_{1}^n\sin^2\frac{\vr\theta_j}{2}d\theta_j
  & \mbox{for  ${\frak b}_n$,  ${\frak c}_n$}
\end{array}
\right.\\
\\ &=&c'_n\prod_{1\leq j<k'\leq
n}(\cos\theta_j-\cos\theta_k)^2\left\{\begin{array}{ll}
\displaystyle{\prod_{1\leq j\leq n}}d\theta_j
&\mbox{for $\frak d_n$}\\ \displaystyle{\prod_{1\leq
j\leq
n}\left(\frac{1-\cos\theta_j}{2}\right)}d\theta_j
&\mbox{for $\frak b_n$}\\ \displaystyle{\prod_{1\leq
j\leq n}}(1-\cos^2\theta_j)d\theta_j &\mbox{for $\frak
c_n$}
\end{array}
\right.\\ 
&=&\left\{\begin{array}{ll}
 c'_n\Dt^2(z)\displaystyle{\prod_{1\leq j\leq
n}\frac{dz_j}{\sqrt{1-z^2_j}}}& \mbox{for ${\frak
d}_n$}\\ c'_n\Dt^2(z)\displaystyle{\prod_{1\leq j\leq
n}(1-z_j)\frac{dz_j}{\sqrt{1-z^2_j}}}& \mbox{for
${\frak b}_n$}\\
c'_n\Dt^2(z)\displaystyle{\prod_{1\leq j\leq
n}(1-z^2_j)\frac{dz_j}{\sqrt{1-z^2_j}}} & \mbox{for
${\frak c}_n$}\\
\end{array}
\right.
\\ &=&c''_n\Dt^2(z)\prod_{1\leq j\leq
n}(1-z_j)^{\al}(1+z_j)^{\beta}dz_j\mbox{\,\,with\,\,}
\left\{\begin{array}{ll}
 \al=\beta
=-1/2&\mbox{\,\,for $\frak d_n$}\\
 \al=1/2,\beta
=-1/2&\mbox{\,\,for $\frak b_n$}\\ \al=\beta
=1/2&\mbox{\,\,for $\frak c_n$}
\end{array}
\right. \eean
\medbreak

 For $M \in SO(2n),~Sp(n)$, the
eigenvalues are given by $e^{i\theta_j}$ and
$e^{-i\theta_j},~1 \leq j \leq n $; therefore, setting
$f=\exp (\sum t_k  tr M^k)$ in formula (1.3.1), leads
to
 \be
 e^{\sum_{1}^{\iy} t_k\Tr M^k}
=e^{\sum^{\iy}_{1} t_k
\sum^n_{j=1}(e^{ik\theta_j}+e^{-ik\theta_j}) }=
\prod_{j=1}^ne^{2\sum^{\iy}_{k=1} t_k \cos k
\theta_j}=\prod_{j=1}^ne^{2\sum  t_k T_k(z_j)},
 \ee
 where $T_n(z)$ are the Tchebychev polynomials,
 defined by $T_n(\cos \theta):=\cos n\theta$;
 in particular $T_1(z)=z$.

 For $M \in SO(2n+1)$, the eigenvalues are given by
$1,~e^{i\theta_j}$ and $e^{-i\theta_j},~1 \leq j \leq
n $, which is responsible for the extra-exponential
$e^{\sum t_i}$ appearing in (1.3.2).

 Before listing various integrals, define
the Jacobi weight
\be
\rho_{\alpha\beta}(z)dz:=(1-z)^{\alpha}(1+z)^{\beta}
dz, \ee 
and the formal sum
 $$g(z):= 2\sum_1^{\iy}t_iT_i(z).$$
The arguments above lead to the following integrals,
originally due to H. Weyl \cite{Weyl}, and in its
present form, due to
 Johansson \cite{Johansson2}; besides the
 integrals over $SO(k)=O_+(k)$, the integrals over $O_-(k)$ and
 $U(n)$
 will also be of interest in the theory of random
 permutations:
\bea
 \int_{O(2n)_+}e^{\sum_{1}^{\iy} t_i tr M^i}dM
  &=&
 \int_{[-1,1]^n}\Delta_n(z)^2
 \prod_{k=1}^n e^{g(z_k)}
  \rho_{(-\frac{1}{2},-\frac{1}{2})}(z_k) dz_k\nonumber\\
   \int_{O(2n+1)_+}e^{\sum_{1}^{\iy} t_i tr
M^i}dM&=& e^{\sum_{1}^{\iy} t_i}
 \int_{[-1,1]^n}\Delta_n(z)^2
 \prod_{k=1}^n e^{g(z_k)}
  \rho_{(\frac{1}{2},-\frac{1}{2})}(z_k) dz_k\nonumber\\
  \int_{Sp(n)}e^{\sum_{1}^{\iy} t_i tr M^i}dM
  &=&
 \int_{[-1,1]^n}\Delta_n(z)^2
 \prod_{k=1}^n e^{g(z_k)}
  \rho_{(\frac{1}{2},\frac{1}{2})}(z_k) dz_k.\nonumber\\
   \int_{O(2n)_-}e^{\sum_{1}^{\iy} t_i tr M^i}dM
  &=& e^{\sum_{1}^{\iy} 2 t_{2i}}
 \int_{[-1,1]^{n-1}}\Delta_{n-1}(z)^2
 \prod_{k=1}^{n-1} e^{g(z_k)}
  \rho_{(\frac{1}{2},\frac{1}{2})}(z_k) dz_k\nonumber\\
  \int_{O(2n+1)_-}e^{\sum_{1}^{\iy} t_i tr M^i}dM
  &=& e^{\sum_{1}^{\iy} (-1)^i t_i}
 \int_{[-1,1]^n}\Delta_n(z)^2
 \prod_{k=1}^n e^{g(z_k)}
  \rho_{(-\frac{1}{2},\frac{1}{2})}(z_k) dz_k\nonumber
  \eea
  \bea
  \int_{U(n)}e^{\sum_1^{\iy}tr (t_iM^i-s_i\bar
M^i)} dM &=&
 \frac{1}{n!}
 \int_{(S^1)^{n}}|\Dt_n(z)|^{2}
 \prod_{k=1}^n
e^{\sum_1^{\iy}(t_i z_k^i-s_iz_k^{-i})}
 \frac{dz_k}{2\pi i z_k}
\nonumber\\
 \eea


\subsection{Permutations and integrals over groups}

Let $S_n$ be the group of permutations $\pi$ and
$S_{2n}^0$ the subset of fixed-point free involutions
$\pi^0$(i.e., $(\pi^0) ^2=I$ and $\pi^0 (k) \neq k$
for $1\leq k \leq 2n$ ). Put the uniform distribution
on $S_n$ and $S_{2n}^0$ ; i.e., all permutations or
involutions have equal probability: \be P(\pi_n)=1/n!
~~\mbox{and}~~ P(\pi^0_{2n})=
 \frac{2^nn!}{(2n)!} ; \ee
$\pi_n$ refers to
 a permutation in $S_n$ and $\pi^0_{2n}$ to
 an involution in $S^0_{2n}$.

An {\em increasing subsequence} of $\pi \in S_n$ or
$S_{n}^0$ is a sequence $1\leq j_1<...< j_k\leq n$,
such that $\pi (j_1)<...<\pi(j_k)$. Define
\be
L(\pi_n) =  \mbox{ length of the longest increasing
subsequence of $\pi_n$ }.
 \ee
 Example: for $\pi=( \underline{3},1,\underline{4},2,\underline{6}
  ,\underline{7},5)$, we have
  $L(\pi_7)=4$.

Around 1960 and based on Monte-Carlo methods, Ulam
\cite{Ulam} conjectured that
 $$
\lim_{n\rightarrow \iy} \frac{E(L_n)}{\sqrt n}=c
~~\mbox{exists}.
 $$
An argument of Erd\"os \& Szekeres \cite{Erdos and
Szekeres}, dating back from 1935 showed that
 $E(L_n)\geq \frac{1}{2}\sqrt{n-1}$, and thus $c\geq
1/2$. In '72, Hammersley \cite{Hammersley} showed
rigorously that the limit exists. Logan and
Shepp \cite{Logan and Shepp} showed the limit $c\geq 2$, and finally
Vershik and Kerov
\cite{Vershik and Kerov} that $c=2$. In
1990, I. Gessel \cite{Gessel} showed that the
following generating function is the determinant of a
Toeplitz matrix:
 \be
 \sum_{n=0}^{\iy}\frac{t^n}{n!} P(L_n\leq \ell)=
 \det \left(\int _0^{2\pi} e^{2\sqrt
 t \cos \theta }e^{i(k-m )\theta} d\theta\right)
 _{0 \leq k,m\leq \ell-1}
 .\ee
The next major contribution was due to Johansson
\cite{Johansson1} and Baik-Deift-Johansson
\cite{Baik-Deift-Johansson}, who prove that for
arbitrary $x \in \BR$, we have a "{\em law of large
numbers}" and a "{\em central limit theorem}", where
$F(x)$ is the statistics (1.2.10),
 $$
 \lim_{n\rightarrow \iy} \frac{L_n}{2\sqrt
 n}=1,~~\mbox{and}~~
 P\left( \frac{L_n-2\sqrt n}{n^{1/6}} \leq x   \right)\longrightarrow F(x)
,~~~ \mbox{ for } n\longrightarrow \iy.
 $$

\noindent A next set of ideas is due to Diaconis \&
Shashahani \cite{DS}, Rains \cite{RThesis,R}, Baik \&
Rains \cite{BR}. For a nice state-of-the-art account,
see Aldous \& Diaconis \cite{AD}.
An illustration is contained in the following
proposition; the first statement is essentially
Gessel's and the next statement is due to
\cite{DS,R,BR}.

\begin{proposition} The following holds
 \bea
\sum^{\iy}_{n=0}\frac{t^n}{n!}P(L(\pi_n)\leq\ell)
 &=&\int_{U(\ell)}e^{\sqrt{t}\Tr (M+\bar M)}dM\\
&=&\int_{[0,2\pi]^{\ell}}\prod_{1\leq j<k\leq\ell}
|e^{i\theta_j}-e^{i\theta_k}|^2\prod_{1\leq
k\leq\ell}e^{2\sqrt{t}\cos\theta_k}\frac{d\theta_k}{2\pi}.
\nonumber \eea
 \be
  \sum^{\iy}_{n=0}\frac{(t^2/2)^n}{n!}P(
L(\pi^0_{2n})\leq \ell)=\int_{O(\ell)}e^{t\Tr M}dM.\ee

\end{proposition}

The \underline{\sl proof}  of this statement will be
sketched later. The connection with integrable systems
goes via the following chain of ideas:

 $$ \mbox{Combinatorics}$$ $$\downarrow
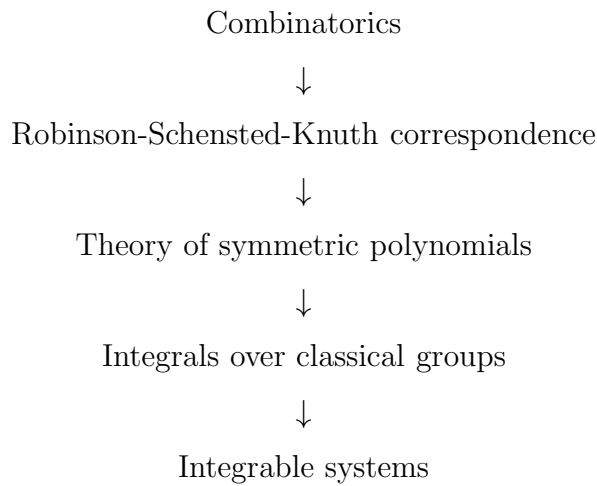
$$
 $$\mbox{Robinson-Schensted-Knuth correspondence}$$
 $$\downarrow$$
 $$\mbox{Theory of symmetric polynomials}$$
 $$\downarrow$$
 $$ \mbox{Integrals over classical groups}$$
 $$\downarrow$$
   $$\mbox{Integrable systems}$$

All the arrows, but the last one, will be explained in
this section; the last arrow will be discussed in
sections 7 and 8. We briefly sketch a few of the basic
well known facts going into these arguments. They can
be found in MacDonald \cite{MacDonald}, Knuth
\cite{Knuth}, Aldous-Diaconis \cite{AD}. Useful facts
on symmetric functions, applicable to integrable
theory, can be found in the appendix to \cite{AvM0}.
Let me mention a few of these facts:

\begin{itemize}
  \item A {\em Young diagram} $\lambda$ is a finite sequence of
  non-increasing, non-negative integers $\lambda_1\geq
  \lambda_2 \geq ...\geq \lambda_{\ell}\geq 0$; also
  called a {\em partition } of
  $n=|\lambda|:=\lambda_1+...+\lambda_{\ell}$, with
  $|\lambda|$ being the weight.  It can be
  represented by a diagram, having $\lambda_1$ boxes
  in the first row, $\lambda_2$ boxes in the second
  row, etc..., all aligned to the left.
 A {\em dual Young diagram} $\hat \lambda=(\hat
  \lambda_1 \geq \hat \lambda_2 \geq...)$ is the
  diagram obtained by flipping the diagram $\lambda$
  about its diagonal.

  \item A {\em Young tableau} of shape $\lambda$
  is an array of positive
   integers $a_{ij}$ (at place $(i,j)$ in the
   Young diagram) placed in the Young diagram $\lambda$, which are
   non-decreasing from left to right {\em and} strictly
   increasing from top to bottom.

  \item A {\em standard Young tableau} of shape $\lambda$ is an array of
   integers $1,...,n$ placed in the Young diagram, which are
   strictly
   increasing from left to right {\em and} from top to bottom.
   The number of Young
   tableaux of a given shape $\lambda=
   (\lambda_1\geq...\geq\lambda_m)$ is given
   by a number of formulae (for the Schur polynomial
    $s_{\lb}$, see below)\footnote{$h^{\lb}_{ij}:=\lb_i +\hat \lb_j-i-j+1$
 is the {\em hook length} of the $i,j$th box in the
 Young diagram; i.e., the length of the hook formed
 by drawing a horizontal line emanating from the center
 of the box to the right and a vertical line emanating
 from the center of the box to the bottom of the diagram.
  }
  \bea f^{\lb}&=&\# \{\mbox{standard tableaux of shape
$\lb$}\}\nonumber\\ &=&\mbox{coefficient of
$x_1x_2\ldots x_n$ in $s_{\lb}(x)$} \nonumber
\\&=&\frac{|\lb
|!}{\prod_{\mbox{\footnotesize{all}}~i,
j}h^{\lb}_{ij}}~=~ |\lb |
!\det\left(\frac{1}{(\lb_{i}-i+j)!}\right)\quad
\\
  &=& |\lb | ! \prod_{1\leq i< j \leq
m}(h_i-h_j)\prod_1^m
\frac{1}{h_i!},~~\mbox{with}~h_i:=\lambda_i-i+m,~m:=\hat
\lambda_1 \nonumber
 . \eea

  \item The {\em Schur polynomial $s_{\lambda}$} associated with a Young
  diagram $\lambda$ is a symmetric function in the
  variables $x_1,x_2,...$ (finite or infinite),
  defined by
  \be
  s_{\lambda}(x_1,x_2,...):=\sum_{\{a_{ij}\}~ \mbox{\tiny{tableaux of
  }}\lambda}~\prod_{ij}x_{a_{ij}}.
  \ee

 \item The linear {\em space $\Lambda_n$ of symmetric
  polynomials}
  in $x_1,...,x_n$ with rational coefficients comes
  equipped with the inner product
  \bea
  \la f,g \ra&=&\frac{1}{n!}\int _{(S_1)^n} f(z_1,...,z_n)
   g(\bar z_1,...,\bar z_n)\prod_{1 \leq k<\ell \leq
  n}|z_k-z_{\ell}|^2 \prod _1^n \frac{dz_k}{2\pi iz_k}
  \nonumber\\
  &=& \int_{U(n)}f(M){g(\bar M)}dM.
  \eea

  \item An {\em orthonormal basis of the space
  $\Lambda_n$}
  is given by the Schur polynomials
  $s_{\lambda}(x_1,...,x_n)$, in which the numbers
  $a_{ij}$ are restricted to $1,...,n$. Therefore,
  each symmetric function admits a ``{\em Fourier
  series}"
  \be \hspace{-1cm}
  f(x_1,...,x_n)=
 \sum_{\stackrel{\lambda~\mbox{\footnotesize{with}}}{\hat \lambda_1 \leq
  n}}  \la f, s_{\lambda} \ra ~ s_{\lambda}(x_1,...,x_n),
 ~ \mbox{with}~ \la s_{\lambda},
  s_{\lambda'}\ra=\delta_{\lambda\lambda'}.
  \ee
   In particular, one
proves (see (1.4.6) for the definition of
$f^{\lb}$)\be (x_1+\ldots+x_n)^k=\sum_{{|\lb
|=k}\atop{\hat\lb_1\leq n}}f^{\lb}s_{\lb}, \ee
  If
$\lambda=(\lambda_1\geq ...\geq \lambda_{\ell}>0)$,
with\footnote{Remember, from the definition of the
dual Young diagram, that $\hat\lb_1=\mbox{the length
of the first column of } \lb$} $\hat \lambda_1=\ell >
 n$, then obviously $s_{\lambda}=0$.

  \item

\noindent{\em Robinson-Schensted-Knuth
correspondence}: There is a 1-1 correspondence $$
S_n\longrightarrow\left\{\begin{array}{l} (P,Q),
\mbox{two standard Young}\\ \mbox{tableaux from
$1,\ldots,n$, where}\\ \mbox{$P$ and $Q$ have the same
shape}
\end{array}
\right\} $$

Given a permutation $i_1,...,i_n$, the correspondence
constructs two standard Young tableaux $P,Q$ having
the same shape $\lambda$. This construction is
inductive. Namely, having obtained two equally shaped
Young diagrams $P_k,Q_k$ from $i_1,...,i_k$, with the
numbers $(i_1,...,i_k)$ in the boxes of $P_k$ and the
numbers $(1,...,k)$ in the boxes of $Q_k$, one creates
a new diagram $Q_{k+1}$, by putting the {\em next
number $i_{k+1}$ in the first row of $P$}, according
to the following rule:
\begin{description}
  \item[(i)] if $i_{k+1} \geq $ all numbers appearing
  in the first row of $P_k$, then one creates a new box
  with $i_{k+1}$ in that box to the right of the first
  column,
  \item[(ii)] if not, place $i_{k+1}$ in the box
  (of the first row) with
  the smallest number higher than $i_{k+1}$. That number then gets
  pushed down to the second row of $P_k$ according to
  the rule (i) or (ii), as if the first row had been removed.
\end{description}
The diagram Q is a bookkeeping device; namely, add a
box (with the number $k+1$ in it) to $Q_k$ exactly at
the place, where the new box has been added to $P_k$.
This produces a new diagram $Q_{k+1}$ of same shape as
$P_{k+1}$.

The inverse of this map is constructed essentially by
reversing the steps above.

\end{itemize}

\begin{example}  $\pi=(5,1,4,3,2)\in S_5$,

{\footnotesize $$
\begin{array}{cccccccccccccccccccc}
5  &&&&1&&&&1&4&&&  &1&3&&&&1&2\\
   &&&&5&&&&5& &&&  &4& &&&&3& \\
   &&&& &&&& & &&&  &5& &&&&4& \\
   &&&& &&&& & &&&  & & &&&& 5& \\
   &&&& &&&& & &&&  & & &&&& &  \\
 1 &&&&1&&&&1&3&&&  &1&3&&&&1&3\\
   &&&&2&&&&2& &&&  &2& &&&&2 &\\
   &&&& &&&& & &&&  &4& &&&&4& \\
   &&&& &&&& & &&&  & & &&&&5 &
\end{array}
$$}

Hence $\displaystyle{\pi\longrightarrow
(P(\pi),Q(\pi))=\left(\left(\begin{array}{cc} 1&2\\ 3&
\\ 4& \\ 5&
\end{array}\right),
\left(\begin{array}{cc} 1&3\\ 2& \\ 4& \\ 5&
\end{array}\right)\right)}$

\noindent and so $ L_5(\pi)=2=\# \mbox{columns of $P$
or $Q$.}
$

\end{example}

\vspace{.3cm}

The Robinson-Schensted-Knuth correspondence has the
following properties
\begin{itemize}
  \item $\pi\mapsto (P,Q)$, then $\pi^{-1}\mapsto (Q,P)$
  \item length (longest increasing subsequence of $\pi$) $=\#$ (columns
  in $P$)
  \item length (longest decreasing subsequence of $\pi$) $=\#$ (rows
  in $P$)
  \item $\pi^2=I$, then $\pi\mapsto (P,P)$
  \item $\pi^2=I$, with $k$ fixed points, then $P$ has exactly $k$
  columns of odd length.
 \vspace{-.7cm}
  \be \ee
\end{itemize}

\noindent From representation theory (see Weyl
\cite{Weyl} and especially Rains {\cite{RThesis}), one
proves:

\begin{lemma} The following perpendicularity relations
hold:

 \noindent(i)  $\displaystyle{ \int_{U(n)}s_{\lb}(M)  s_{\mu}(\bar M) dM=\la
s_{\lb}, s_{\mu}\ra= \dt_{\lb\mu} }$

  \bean \hspace{-.5cm} (ii)
\int_{O(n)}s_{\lb}(M)dM&=&1\quad\mbox{for
$\lb=(\lb_1\geq  \ldots\geq\lb_k\geq 0)$,
$k\leq n$,}
~~\mbox{$\lb_i$ even}\nonumber\\
&=&0\quad\mbox{otherwise.} \eean
 \bea \hspace{-5cm}
(iii) \int_{Sp(n)}s_{\lb}(M)dM&=&1\quad\mbox{for $\hat
\lb_i$ even, $\hat\lb_1 \leq 2n$
,}
\nonumber\\
&=&0\quad\mbox{otherwise.} \eea

\end{lemma}

\noindent{\it Proof of Proposition 1.1 :\/} On the one
hand,

\medbreak

 $\la (x_1+\ldots+x_n)^k,(x_1+\ldots+x_n)^k\ra$
\bea &=& \sum_{{|\lb |=|\mu
|=k}\atop{\hat\lb_1,\hat\mu_1\leq n}}f^{\lb}f^{\mu}
 \la s_{\lb},s_{\mu}\ra \nonumber\\ &=&\sum_{{|\lb
|=k}\atop{\hat\lb_1\leq n}}(f^{\lb})^2\nonumber\\ & &
\nonumber\\ &=&\sum_{{|\lb |=k}\atop{\lb_1\leq
n}}(f^{\lb})^2\nonumber\\ & & \nonumber\\
&=&\#\{(P,Q),\mbox{\,\,standard Young tableaux, each
of arbitrary shape $\lb$}\nonumber\\ &
&\hspace{7cm}\mbox{ with $|\lb |=k,\,\,\lb _1 \leq
n$}\}\nonumber\\ &=&\#\{\pi_k\in S_k\mbox{\,\,such
that $L(\pi_k)\leq n$}\}. \eea

On the other hand, notice that, upon setting
$\theta_{j}=\theta'_j+\theta_1$ for $2\leq j\leq n$,
the expression $\displaystyle{\prod_{1\leq j<k\leq
n}|e^{i\theta_j}-e^{i\theta_k}|^2}$ is independent of
$\theta_1$. Then, setting $z_k=e^{i\theta_k}$, one
computes:


\medbreak

 $\la (x_1+\ldots+x_n)^k,(x_1+\ldots+x_n)^{\ell}\ra$
\bea &=&\frac{1}{n!}\int_{[0,2\pi]^n}(z_1+\ldots+z_n)^k(\bar
z_1+\ldots+\bar z_n)^{\ell}\prod_{1\leq j<k\leq
n}|e^{i\theta_j}-e^{i\theta_k}|^2d\theta_1\ldots
d\theta_n\quad \nonumber\\
&=&\frac{1}{n!}\int_{[0,2\pi]^n}e^{ik\theta_1}(1+z'_2+\ldots+z'_n)^k
 e^{-i\ell\theta_1}(1+\bar z'_2+\ldots+\bar z'_n)^{\ell}
\nonumber \\&&\hspace{7cm} \prod_{1\leq j<k\leq
n}|e^{i\theta_j}-e^{i\theta_k}|^2d\theta_1\ldots
d\theta_n \nonumber \\ & &\hspace{3cm}\mbox{upon
setting $\theta_j=\theta'_j+\theta_1,\quad\mbox{for }
j\geq 2$ and $z_k^{\prime}=e^{i\theta_k^{\prime}},$}
\nonumber \\
&=&\frac{1}{n!}\int_0^{2\pi}e^{i(k-\ell)\theta_1}d\theta_1
~~~{\times }~~~\mbox{\,\,(an~ $n-1$-fold
integral)}\nonumber\\ &=&\dt_{k\ell}\la
(x_1+\ldots+x_n)^k,(x_1+\ldots+x_n)^k\ra 
 =\dt_{k\ell}\int_{U(n)} |\Tr M|^{2k}dM. \eea
 It follows that
 \bea \lefteqn{\int_{U(n)}(\tr(M+\bar
M))^kdM}\nonumber\\&=&\sum_{0\leq j\leq
k}\left(\begin{array}{c}k\\j\end{array}\right)\int_{U(n)}((\tr
M)^j(\overline{\tr M})^{k-j})dM \nonumber\\
&=&\left\{\begin{array}{l} 0, \mbox{\,\,if $k$ is odd
(because then $j\neq k-j$ for all $0\leq j\leq k$}) \\
\left(\begin{array}{c} k\\ k/2
\end{array}\right)\int_{U(n)} |\tr M|^kdM,\mbox{\,\,if $k$ is even}.
\end{array}
\right. \eea Combining the three identities (1.4.13),
(1.4.14) and (1.4.15) leads to \be \#\{\pi_k\in
S_k\mbox{\,\,such that\,\,}L(\pi_k)\leq
n\}=\left(\begin{array}{c} 2k\\ k
\end{array}\right)^{-1}\int_{U(n)}(\Tr(M+\bar
M))^{2k}dM. \ee  Finally

\medbreak

\bean
 \lefteqn{\sum^{\iy}_{n=0}\frac{t^n}{n!}P(L(\pi_n)\leq\ell)}\\
 &=&\sum^{\iy}_{n=0}
\frac{t^n}{n!}\frac{\#\{\pi_n\in S_n ~\bigl |~L(\pi_n)\leq
\ell \} }{n!}\\
&=&\sum^{\iy}_{n=0}\frac{t^n}{(n!)^2}\left(\begin{array}{c}
2n\\ n
\end{array}\right)^{-1}\int_{U(\ell)}(\tr(M+\bar M))^{2n}dM\\
&=&\sum^{\iy}_{n=0}\frac{(\sqrt{t})^{2n}}{(2n)!}\int_{U(\ell)}(\tr(M+\bar
M))^{2n}dM\\ &=&\int_{U(\ell)}e^{\sqrt{t}\Tr(M+\bar
M)}dM\\ &=&\frac{1}{\ell
!}\int_{[0,2\pi]^{\ell}}e^{\sqrt{t}(z_1+z^{-1}_1+\ldots
+z_{\ell}+z^{-1}_{\ell})}\prod_{1\leq
j<k\leq\ell}|e^{i\theta_j}-e^{i\theta_k}
|^2\prod_{1\leq k\leq\ell}\frac{d\theta_k}{2\pi},\\ &
&\hspace{8cm}\mbox{where}~ z_k=e^{i\theta_k},\\
&=&\frac{1}{\ell !}\int_{[0,2\pi]^{\ell}}\prod_{1\leq
j<k\leq\ell}|e^{i\theta_j}-e^{i\theta_k}
|^2\prod^{\ell}_{k=1}
e^{2\sqrt{t}\cos\theta_k}\frac{d\theta_k}{2\pi},\\
\eean
 showing (1.4.4) of Proposition 1.1. The latter
also equals:
 \bean &=&
 \frac{1}{{\ell}!}\int_{(S^1)^{\ell}}\Dt_{\ell}(z)\Dt_{\ell}(\bar z)
 \prod_{k=1}^{\ell}
\left(e^{ \sqrt{t}(z_k+\bar z_k)}
 \frac{ dz_k }{2\pi i z_k}\right)\nonumber\\
&=&
 \frac{1}{{\ell}!}\int_{(S^1)^{{\ell}}}\sum_{\sigma\in S_{\ell}}\det
 \left(z_{\sigma(m)}^{k-1}
 \bar z_{\sigma(m)}^{m-1}  \right)_{1\leq k,m\leq {\ell}}
 \prod_{k=1}^{\ell}
\left(e^{ \sqrt{t}(z_k+\bar z_k)}
 \frac{dz_k}{2\pi i z_k}\right)
 \nonumber\\
&=&
 \frac{1}{{\ell}!}\sum_{\sigma\in S_{\ell}}
  \det\left(\int_{S^1}z_{k}^{k-1}
 \bar z_{k}^{m-1}
e^{ \sqrt{t}(z_k+\bar z_k)}
 \frac{dz_k}{2\pi i z_k}\right)_{1\leq k,m\leq {\ell}}
  \nonumber\\
  &=&
 \det \left(\int _0^{2\pi} e^{2\sqrt
 t \cos \theta }e^{i(k-m )\theta} d\theta\right)
 _{1 \leq k,m\leq \ell},\nonumber\\
 \eean
confirming Gessel's result (1.4.3).
\medbreak

The proof of the second relation (1.4.5) of
Proposition 1.1 is based on the following computation:

\bean \int_{O(n)}(\Tr M)^kdM&=&\sum_{{ |\lb  |
=k}\atop{\hat\lb_1\leq n}}f^{\lb}
\int_{O(n)}s_{\lb}(M)dM,~~\mbox{using (1.4.10),}\\ & &
\\ &=&\renewcommand{\arraystretch}{0.5}
\begin{array}[t]{c}
\sum\\ {\scriptstyle |\lb  | =k}\\ {\scriptstyle
\hat\lb_1\leq n}\\ {\scriptstyle
\lb_i\mbox{\,\,\small{even}}}
\end{array}
\renewcommand{\arraystretch}{1}f^{\lb}
~,~~\mbox{using Lemma 1.2,}\\
& & \\ &=&\renewcommand{\arraystretch}{0.5}
\begin{array}[t]{c}
\sum\\ {\scriptstyle |\lb  | =k}\\ {\scriptstyle
\lb_1\leq n}\\ {\scriptstyle
\hat\lb_i\mbox{\,\,\small{even}}}
\end{array}
\renewcommand{\arraystretch}{1}f^{\lb}~,~~\mbox{using duality,}\\
& & \\ &=&\#\left\{\begin{array}{l} (P,P), P
\mbox{\,\,standard Young tableau of shape $\lb$}\\
\mbox{with $|\lb |=k$, $\lb_1\leq n$, $\hat\lb_i$
even}
\end{array}
\right\}\\ & & \\ &=&\#\left\{ \pi_k^{0}\in S^{0}_k,
\mbox{ no fixed points and $L(\pi_k^{0})\leq n$}
\right\}. \eean
\vspace{-0.8cm}
\be \ee

In the last equality, we have used property (1.4.11):
an involution has no fixed points iff all columns of
$P$ have even length. Since all columns $\hat\lb_i$
have even length, it follows that $|\lb|=k$ is even
and then only is $\int_{O(n)}(\Tr M)^kdM>0$; otherwise
$=0$ . Finally, one computes,
 \bean
  \lefteqn{\sum^{\iy}_{k=0}\frac{(t^2/2)^k}{k!}
 P\left( L(\pi_{2k}^{0})\leq n,~
\pi_{2k}^{0}\in S^{0}_{2k} \right)}\\
&=&\sum^{\iy}_{k=0}\frac{t^{2k}}{2^kk!}\frac{2^kk!}{(2k)!}
 \#\{\pi_{2k}^0\in
S^0_{2k},~L(\pi_{2k}^0)\leq n\}~,~~\mbox{using
(1.4.1)},\\ & & \\
&=&\sum^{\iy}_{k=0}\frac{t^k}{k!}\#\{\pi^0_{k}\in
S^0_{k},\,\,L(\pi^0_{k})\leq n\}\\ & & \\
&=&\sum^{\iy}_{k=0}\frac{t^k}{k!}\int_{O(n)}(\Tr
M)^kdM ~,~~\mbox{using (1.4.17)},\\ & & \\
&=&\int_{O(n)}e^{t \Tr M} dM, \eean
 ending the proof of Proposition 1.1.\qed

\newpage

\section{Integrals, vertex operators and Virasoro
relations}

In section 1, we discussed random matrix problems over
different finite and infinite matrix ensembles,
generating functions for the statistics of the length
of longest increasing sequences in random permutations
and involutions. One can also consider two Hermitian
random matrix ensembles, coupled together. All those
problems lead to matrix integrals or Fredholm
determinants, which we list here: ($\beta=2,1,4$)
\begin{itemize}
  \item ~~~~$
  \displaystyle{\int_{{\cal H}_n(E),~{\cal S}_n(E)~
  \mbox{\tiny{or}}
  ~{\cal T}_n(E)}e^{-Tr V( M)}dM}=
  c_n \int_{E^n} |\Delta_n(z)|^{\beta}
\prod^n_1\rho(z_k)dz_k
  $

  \item ~~~~$\displaystyle{\int\int_{{\cal H}_n^2(E_1\times
  E_2)}}dM_1dM_2\,e^{-\frac{1}{2}{\rm
Tr}(M^2_1+M_2^2-2cM_1M_2)}$

\item ~~~~$\displaystyle{\int_{O(n)}e^{x\Tr M} dM}$

  \item ~~~~$\displaystyle{\int_{U(n)}e^{\sqrt{x}\Tr(M+\bar
M)}dM  }  $

  \item ~~~~$\displaystyle{\det \bigl(I-\lambda K(y,z)
   I_E(z)\bigr)}$
  with $K(y,z)$ as in (1.2.4).
  \vspace{-0.78cm}  \be  \ee
\end{itemize}
The point is that each of these quantities admit a
natural deformation, by inserting time variables
$t_1,t_2,...$ and possibly a second set $s_1,s_2,...$,
in a seemingly {\sl ad hoc} way. Each of these
integrals or Fredholm determinant is then a fixed
point for a natural {\em vertex operator}, which
generates a Virasoro-like algebra. These new integrals
in $t_1,t_2,...$ are all annihilated by the precise
subalgebra of the Virasoro generators, which
annihilates $\tau_0$. This will be the topic of this
section.

 \subsection{Beta-integrals}

\subsubsection{Virasoro constraints for Beta-integrals}

Consider weights of the form $\rho(z)dz:=e^{-V(z)}dz$
on an interval $F=[A,B]\subseteq\BR$, with rational
logarithmic derivative and subjected to the following
boundary conditions:
\be
-\frac{\rho'}{\rho}=V^{\prime}=\frac{g}{f}=\frac{\sum_0^{\iy}b_iz^i}{\sum_0^{\iy}a_iz^i},
\quad \lim_{z\rightarrow
A,B}f(z)\rho(z)z^k=0\mbox{\,\,for all\,\,}k\geq 0,
  \ee
and a disjoint union of intervals,
  \be
E=\bigcup_1^{2r}~[c_{2i-1},c_{2i}]\subset F\subset
\BR. \ee These data define an algebra of differential operators
 \be{\cal
B}_k=\sum_1^{2r} c_i^{k+1}f(c_i)\frac{\pl}{\pl c_i}
 .\ee
   Take the first type
of integrals in the list (2.0.1) for general $\beta >0$,
thus generalizing the integrals, appearing in the probabilities
 (1.1.9),(1.1.11) and (1.1.18).
 Consider $t$-deformations of such
integrals, for general (fixed) $\beta >0$:
($t:=(t_1,t_2,...)$, $c=(c_1,c_2,...,c_{2r})$ and
$z=(z_1,...,z_n)$)
\be
I_n(t,c;\beta):=\int_{E^n}|\Dt_n(z)|^{\beta}\prod_{k=1}^n
\left(e^{\sum_1^{\iy}t_i z_k^i}\rho(z_k)dz_k\right)
~\mbox{for} ~~n>0.\ee
 The main statement of this section is
 Theorem 2.1, whose proof will be outlined in the
 next subsection. The central charge (2.1.9) has
 already appeared in the work of Awata et al.
 \cite{Awata}.

 \begin{theorem} {\em (Adler-van Moerbeke \cite{AvM2,AvM3})}
 The multiple integrals
 \be
I_n(t,c;\beta):=\int_{E^n}|\Dt_n(z)|^{\beta}\prod_{k=1}^n
\left(e^{\sum_1^{\iy}t_i z_k^i}\rho(z_k)dz_k\right)
,~\mbox{for} ~~n>0
 \ee
 and
 \be
I_n(t,c;\frac{4}{\beta}):=\int_{E^n}|\Dt_n(z)|^{4/\beta}\prod_{k=1}^n
\left(e^{\sum_1^{\iy}t_i z_k^i} \rho(z_k)dz_k\right)
,~\mbox{for} ~~n>0,
 \ee
 with $I_0=1$, satisfy respectively the following
 Virasoro constraints\footnote{
 When $E$ equals the whole range $F$, then the
  ${\cal B}_k$'s are absent in the formulae (2.1.7).} for all $k\geq -1$:
\bea && \hspace{-1cm}\left(-{\cal B}_k+ \sum_{i\geq
0}\left( a_i~ {}^{\beta}\BJ_{k+i,n}^{(2)}(t,n)-b_i ~
{}^{\beta} \BJ_{k+i+1,n}^{(1)}(t,n)\right)
\right)I_n(t,c;\beta) = 0, \nonumber \\ &&
\hspace{-1cm} \left(-{\cal B}_k + \sum_{i\geq 0}\left(
a_i~ {}^{\beta}\BJ_{k+i,n}^{(2)}\bigl(-\frac{\beta t
}{2},-\frac{2n}{\beta}\bigr)+\frac{\beta b_i}{2} ~
{}^{\beta} \BJ_{k+i+1,n}^{(1)}\bigl(-\frac{\beta t
}{2},-\frac{2n}{\beta}\bigr)\right)
\right)I_n(t,c;\frac{4}{\beta}) = 0, \nonumber \\
 \eea
  in terms of the coefficients
$a_i,~b_i$ of the rational function $(-\log \rho)'$
and the end points $c_i$ of the subset $E$, as in
(2.1.1) to (2.1.2). For all $n\in \BZ$, the
 $
 ~ {}^{\beta}\BJ_{k,n}^{(2)}(t,n) $
and $
 ~{}^{\beta}\BJ_{k,n}^{(1)}(t,n)$ form a Virasoro and a
Heisenberg algebra respectively, interacting as
follows \bea \left[~ {}^{\beta}\BJ_{k,n}^{(2)},~
{}^{\beta}\BJ_{\ell,n}^{(2)} \right] &=&(k-\ell)~
{}^{\beta}\BJ_{k+\ell,n}^{(2)} +c\left(
\frac{k^3-k}{12} \right)\dt_{k,-\ell}\nonumber\\
\left[~ {}^{\beta}\BJ_{k,n}^{(2)},~
{}^{\beta}\BJ_{\ell,n}^{(1)} \right] &=&-\ell ~ ~
{}^{\beta}\BJ_{k+\ell,n}^{(1)}+c'k(k+1)\delta_{k,-\ell}.
\nonumber\\ \left[~ {}^{\beta}\BJ_{k,n}^{(1)},~
{}^{\beta}\BJ_{\ell,n}^{(1)}  \right]
&=&\frac{k}{\beta}\delta_{k,-\ell}, \eea with central
charge
 \be
c=1-6\left(\left(\frac{\beta}{2}\right)^{1/2}
-\left(\frac{\beta}{2}\right)^{-1/2}  \right)^2 ~~~
\mbox{and}~~~c'=\left(\frac{1}{\beta}- \frac{1}{2}
\right). \ee

\end{theorem}

\medskip\noindent\underline{\it Remark 1:\/}
The $ {}^{\beta}\BJ_{k,n}^{(2)}$'s are defined as
follows:
 \be
 {}^{\beta}\BJ_{k,n}^{(2)}=\frac{\beta}{2}\sum_{i+j=k}
:~ {}^{\beta}\BJ_{i,n}^{(1)} ~
{}^{\beta}\BJ_{j,n}^{(1)}:
+\left(1-\frac{\beta}{2}\right)\left((k+1)~
{}^{\beta}\BJ_{k,n}^{(1)} -k\BJ_{k,n}^{(0)}\right)
.\ee
 Componentwise, we have $$ ~
{}^{\beta}\BJ_{k,n}^{(1)}(t,n)= ~
{}^{\beta}J_k^{(1)}+nJ_k^{(0)}
 ~\mbox{and}~ ~ {}^{\beta}\BJ_{k,n}^{(0)}=nJ_k^{(0)}= n\dt_{0k}
$$ and \bean \lefteqn{ ~{}^{\beta}\BJ_{k,n}^{(2)}(t,n)
}\\&=&\Bigl(\frac{\beta}{2}\Bigr) ~~
{}^{\beta}J_k^{(2)} + \Bigl(n\beta
+(k+1)(1-\frac{\beta}{2})\Bigr) ~~ {}^{\beta}J_k^{(1)}
+ n\Bigl((n-1)\frac{\beta}{2}+1\Bigr) J_k^{(0)},
 \eean
  where
  \bea
 ~ {}^{\beta}J_k^{(1)}&=&\frac{\pl}{\pl
t_k}+\frac{1}{\beta}(-k)t_{-k}
 \\
  ^{\beta}J^{(2)}_{k}&=&\sum_{i+j=k}\frac{\pl^2}{\pl
 t_{i}\pl t_{j}}+\frac{2}{\beta}\sum_{-i+j=k}it_{i}\frac{\pl}{\pl
 t_{j}}+\frac{1}{\beta^2}\sum_{-i-j=k}it_{i}jt_{j}.
\nonumber \eea We put $n$ explicitly in
 $ {}^{\beta}\BJ_{\ell,n}^{(2)}(t,n)$ to indicate
  the
 $n$th component contains $n$ explicitly, besides $t$.
 Note that for $\beta=2$, (2.1.10) becomes
 particularly simple:
$$
 \left.{}^{\beta}\BJ_{k,n}^{(2)}\right|_{\beta=2}=
  \sum_{i+j=k}
:~ {}^{2}\BJ_{i,n}^{(1)} ~ {}^{2}\BJ_{j,n}^{(1)}: .$$

\medskip\noindent\underline{\it Remark 2:\/}
The Heisenberg and Virasoro generators satisfy the
following {\em duality} properties:
  \bea
 {}^{\frac{4}
{\beta}  }\BJ_{\ell,n}^{(2)}\Bigl(t,n \Bigr)
&=&{}^{\beta}\BJ_{\ell,n}^{(2)}\left(-\frac{\beta
t}{2},-\frac{2n}{\beta}\right ) ,~n\in \BZ\nonumber\\
 {}^{\frac{4} {\beta} }\BJ_{\ell,n}^{(1)}\Bigl(t,n
\Bigr) &=&-\frac{\beta}{2}~~
{}^{\beta}\BJ_{\ell,n}^{(1)}\left(-\frac{\beta
t}{2},-\frac{2n}{\beta}\right ),~n>0.
  \eea
 In (2.1.7),
${}^{\beta}\BJ_{\ell,n}^{(2)}\left(- \beta
t/2,-2n/\beta \right )$ means that the variable $n$,
which appears in the $n$th component, gets replaced by
$2n/\beta$ and $t$ by $-\beta t/2$.

\medskip\noindent\underline{\it Remark 3:\/}
 Theorem 2.1 states that the
integrals (2.1.5) and (2.1.6) satisfy two sets of
differential equations (2.1.7) respectively. Of
course, the second integral also satisfies the first
set of equations, with $\beta$ replaced by $4/\beta$.

\subsubsection{Proof: $\beta$-integrals as fixed points of vertex operators}

Theorem 2.1 can be established by using the invariance
of the integral under the transformation $z_i\mapsto
z_i+\vr f(z_i)z_i^{k+1}$ of the integration variables.
However, the most transparent way to prove Theorem 2.1
is via vector vertex operators, for which the
$\beta$-integrals are fixed points. This is a
technique which has been used by us, already in
\cite{AvM1}. Indeed, define the (vector) vertex
operator, for $t=(t_1,t_2,...)\in \BC^{\iy},~u \in
\BC$, and setting $\chi(z):=(1,z,z^2,...)$.
\be
\BX_{\beta}(t,u)=\Lb^{-1}e^{\sum_1^{\iy}t_i u^i}
e^{-\beta\sum_1^{\iy}\frac{ u^{-i}}{i} \frac{\pl}{\pl
t_i}}\chi(|u|^{\beta}), \ee which acts on vectors
$f(t)=(f_0(t),f_1(t),...)$ of functions,
 as follows\footnote{For $\alpha \in \BC$, define $[\alpha]:=(\alpha, \frac{\alpha^2}{2},
 \frac{\alpha^3}{3},...) \in \BC^{\iy}$. The operator $\Lb$ is the
 shift matrix, with zeroes everywhere, except for $1$'s just above
 the diagonal, i.e., $(\Lambda v)_n=v_{n+1}$. }
$$
\bigl(\BX_{\beta}(t,u)f(t)\bigr)_n=e^{\sum_1^{\iy}t_i
u^i}\left(|u|^{\beta }\right)^{n-1} f_{n-1}(t-\beta [u^{-1}]).
$$
 For the sake of these arguments, it is
convenient to introduce the following vector Virasoro
generators: $~ {}^{\beta}\BJ_k^{(i)}(t):=~ ( ~
{}^{\beta}\BJ_{k,n}^{(i)}(t,n))_{n\in \BZ} $.

\begin{proposition} The multiplication operator $z^k$ and the
differential operators $\frac{\pl}{\pl z}z^{k+1}$ with
$z \in \BC^*$, acting on the vertex operator
$\BX_{\beta}(t,z)$, have realizations as commutators,
in terms of the Heisenberg and Virasoro generators $~
{}^{\beta}\BJ_k^{(1)}(t,n)$ and $~
{}^{\beta}\BJ_k^{(2)}(t,n)$ : \bea
z^{k}\BX_{\beta}(t,z)&=& \left[ ~
{}^{\beta}\BJ_k^{(1)}(t), \BX_{\beta}(t,z)\right]
\nonumber\\ \frac{\pl}{\pl
z}z^{k+1}\BX_{\beta}(t,z)&=& \left[~ {}^{\beta}
\BJ_k^{(2)}(t), \BX_{\beta}(t,z)\right]. \eea

\end{proposition}

\begin{corollary}
 Given a weight $\rho(z)dz$ on $\BR$ satisfying (2.1.1), we
 have
 \be
 \frac{\pl}{\pl z}z^{k+1}f(z)\BX_{\beta}(t,z)\rho(z)
 =
\left[ \sum_{i\geq 0}\left( a_i ~ {}^{\beta}\BJ_{k+i}^{(2)}(t)-b_i
 ~ {}^{\beta}\BJ_{k+i+1}^{(1)}(t)\right)
 , \BX_{\beta}(t,z)\rho(z)\right].
\ee

\end{corollary}

\proof Using (2.1.13) in the last line, compute

\noindent $\displaystyle{\frac{\pl}{\pl z}
z^{k+1}f(z)\BX_{\beta}(t,z)\rho(z)
 }$
\bea
&=& \left(\frac{\rho'(z)}{\rho(z)}f(z)
\right)z^{k+1}\BX_{\beta}(t,z)\rho(z)+
 \rho(z)\frac{\pl}{\pl z}\left(
z^{k+1}f(z)\BX_{\beta}(t,z)\right)
\nonumber\\
&=&-\left(\sum_0^{\iy} b_i z^{k+i+1} \BX_{\beta}(t,z)
\right)\rho(z)+\rho(z)\frac{\pl}{\pl z}
 \left(\sum_0^{\iy} a_i z^{k+i+1} \BX_{\beta}(t,z)
\right) \nonumber\\ &=&-\left[\sum_0^{\iy} b_i
~{}^{\beta}\BJ_{k+i+1}^{(1)} ,
  \BX_{\beta}(t,z)\rho(z) \right]
 +\left[\sum_0^{\iy} a_i ~{}^{\beta}\BJ_{k+i}^{(2)}
,\BX_{\beta}(t,z)\rho(z)\right] \nonumber\\ \eea
establishing (2.1.14).\qed


Giving the weight $\rho_E(u)du=\rho(u)I_E(u)du$, with
$\rho$ and $E$ as before, define the integrated vector
vertex operator
\be
\BY_{\beta}(t, \rho_E):= \int_{E} du
\rho(u)\BX_{\beta}(t,u). \ee and the vector operator
\bea \DR_k&:=&{\cal B}_k~~~-~~~\VR_k\nonumber\\
&:=&\sum_1^{2r} c_i^{k+1}f(c_i)\frac{\pl}{\pl c_i}-
\sum_{i\geq 0}\left( a_i ~{}^{\beta}
\BJ_{k+i}^{(2)}(t)-b_i
 ~{}^{\beta}\BJ_{k+i+1}^{(1)}(t)\right),
\eea
consisting of
 a $c$-dependent boundary part ${\cal B}_k$ and a $(t,n)$-dependent Virasoro
 part $\VR_k$.

\begin{proposition} The following commutation relation holds:

\be
\left[\DR_k,\BY_{\beta}(t, \rho_E)\right]=0
\ee
\end{proposition}

\proof
Integrating both sides of (2.1.14) over $E$, one computes:
\begin{eqnarray}
\int_Edz\frac{\pl}{\pl z}\left( z^{k+1}f(z)\BX_{\beta}
(t,z) \rho(z)\right)
&=& \sum_1^{2r} (-1)^ic_i^{k+1}f(c_i)
 \BX_{\beta}(t,c_i) \rho(c_i)\nonumber\\
 &=& \sum_1^{2r}  c_i^{k+1}f(c_i)
 \frac{\pl}{\pl c_i}\int_{E}
 \BX_{\beta}(t,z) \rho(z)dz\nonumber\\
 &=& \left[ {\cal B}_k, \BY_{\beta}(t,\rho_E) \right],
\end{eqnarray}
while on the other hand

\bigbreak

\noindent $\displaystyle{
\int_E dz \left[\sum_{i\geq 0}\left( a_i ~{}^{\beta}\BJ_{k+i}^{(2)}-b_i
 ~{}^{\beta}\BJ_{k+i+1}^{(1)}\right),\BX_{\beta}
(t,z) \rho(z)\right]}
$
 \bea
 &=&\left[\sum_{i\geq 0}\left( a_i ~{}^{\beta}\BJ_{k+i}^{(2)}-b_i
 ~{}^{\beta}\BJ_{k+i+1}^{(1)}\right),\int_{\BR} dz\rho_E(z) \BX_{\beta}
(t,z) \right]\nonumber\\
 &=&\left[ \VR_k, \BY_{\beta}(t, \rho_E) \right].
\eea
Subtracting both expressions (2.1.19) and (2.1.20) yields
$$
0=\left[{\cal B}_k- \VR_k, \BY_{\beta}(t, \rho_E) \right]
 = \left[\DR_k, \BY_{\beta}(t, \rho_E) \right],
 $$
concluding the proof of proposition 2.4.\qed

\begin{proposition}
The column vector,
$$I(t):=\left(\int_{E^n}|\Dt_n(z)|^{\beta}\prod_{k=1}^n
e^{\sum_1^{\iy}t_i z_k^i}\rho(z_k)dz_k\right)_{n\geq
0}$$ is a fixed point for the vertex operator
$\BY_{\beta}(t, \rho_E)$: (see definition (2.1.17))
\be
\left(\BY_{\beta}(t, \rho_E)I\right)_n=I_n,~
n\geq 1. \ee

\end{proposition}

\proof We have
 \bea
I_n(t)&=&\int_{\BR^n}|\Dt_n(z)|^{\beta}
 \prod_{k=1}^n
\left(e^{\sum_1^{\iy}t_i z_k^i}I_E(z_k)\rho(z_k)dz_k\right)
\nonumber\\
&&\nonumber\\
 &=&\int_{\BR}du\rho_E(u)e^{\sum_1^{\iy}t_i u^i}
 |u|^{\beta (n-1)}\nonumber\\
 &&~~~~~~~\int_{\BR^{n-1}}
 \prod_{k=1}^{n-1}\left|1-\frac{z_k}{u}\right|^{\beta}
  |\Dt_{n-1}(z)|^{\beta}\prod_{k=1}^{n-1}
\left(e^{\sum_1^{\iy}t_i z_k^i}
\rho_E(z_k)dz_k\right)\nonumber\\
&=&\int_{\BR}du\rho_E(u)e^{\sum_1^{\iy}t_i u^i}
 |u|^{\beta (n-1)}\nonumber\\
&&~~~~e^{-\beta
\sum_1^{\iy}\frac{u^{-i}}{i}\frac{\pl}{\pl t_i}}
\int_{\BR^{n-1}}|\Dt_{n-1}(z)|^{\beta}
\prod_{k=1}^{n-1} \left(e^{\sum_1^{\iy}t_i z_k^i}
\rho_E(z_k)dz_k\right)\nonumber\\
&=&\int_{\BR}du\rho_E(u)|u|^{\beta
(n-1)}e^{\sum_1^{\iy}t_i u^i}
  e^{-\beta \sum_1^{\iy}\frac{u^{-i}}{i}
\frac{\pl}{\pl t_i}}\tau_{n-1}(t)\nonumber\\
&&\nonumber\\
 &=&\Big(\BY_{\beta}(t, \rho_E)I(t)\Big)_n.
\eea
\qed

{\underline{\sl Proof of Theorem 2.1}: }From
proposition 2.4 it follows that for $n\geq 1$,
  \bea
0&=&\left[\DR_k,\left( \BY_{\beta}(t,\rho_E) \right)^n
\right]I \nonumber\\ &=&\DR_k \BY_{\beta}(t,\rho_E)^n
I- \BY_{\beta}(t,\rho_E)^n\DR_k I.
\eea Taking the $n$th component for $n\geq 1$ and
$k\geq -1$, and setting
 $X_{\beta}(t,u)=e^{\sum t_i u^i}
 e^{-\beta \sum \frac{u^{-i}}{i}\frac{\pl}{\pl t_i}}$,
 and using (2.1.21)
\bea
0&=& \left(\DR_k I-\BY_{\beta}(t,\rho_E)^n\DR_k
 I \right)_n\nonumber\\
 &=&(\DR_k I )_n- \int du\rho_E(u)X_{\beta}(t;u)
 (|u|^{\beta})^{n-1}...\int du\rho_E(u)X_{\beta}(t;u)
 (\DR_k I)_0 \nonumber\\
 &=& (\DR_k I)_n. \nonumber
\eea Indeed $(\DR_k I)_0=0$ for $k\geq -1$, since $\tau_0=1$ and
$\DR_k$ involves $~ {}^{\beta}J_{k}^{(2)},
 ~ {}^{\beta}J_{k}^{(1)}$ and $J_{k}^{(0)}$ for $k \geq -1$:

$\left\{\begin{array}{l}
 {}^{\beta}J_{k}^{(2)}\mbox{ is pure
differentiation for $k\geq -1$;}\hspace{7.9cm} \\
 {}^{\beta} J_{k}^{(1)}\mbox{is pure
differentiation,
 except for $k=-1$; but
}
 \\ \mbox{$~ {}^{\beta}J_{-1}^{(1)}$ appears with
coefficient $n\beta$, which vanishes for
 $n=0$;}
  \\
 \mbox{$J_{k}^{(0)}$ appears with coefficient}~
 n((n-1)\frac{\beta}{2} +1),~\mbox{vanishing for
 $n=0$.}
\end{array}    \right.
   $

\qed

\subsubsection{Examples}

\subsubsection*{Example 1 : Gaussian $\beta$-integrals}
 The weight and the $a_i$ and $b_i$, as in
(2.1.1), are given by
 $$\rho(z)=e^{-V(z)}= e^{-z^2},
~~V'=g/f=2z,$$ $$ a_0=1,b_0=0,b_1=2,~\mbox{and all
other}~a_i,b_i=0 .$$
  From (2.1.5), the integrals
\be
I_n= \int_{E^n}|\Dt_n(z)|^{\beta}\prod_{k=1}^n
e^{-z_k^2+\sum_{i=1}^{\iy}t_iz_k^i}dz_k
 \ee
satisfy the Virasoro constraints, for $k\geq -1$,
 \be
- {\cal B}_k I_n=- \sum_1^{2r}c_i^{k+2}\frac{\pl}{\pl
c_i}I_n=\left(-~{}^{\beta}\BJ_{k+1,n}^{(2)} -a
~{}^{\beta}\BJ_{k+1,n}^{(1)}+~{}^{\beta}\BJ_{k+2,n}^{(1)}
 \right) I_n.\ee
 Introducing the following notation
 \be \sigma_i=
  (n-\frac{i+1}{2})\beta+i+1-b_0=(n-\frac{i+1}{2})\beta+i+1,
\ee
 the first three constraints have the following
form, upon setting $F_n=\log I_n$,
 {\footnotesize \bea
 -{\cal B}_{-1} F&=&
 \left(
2\frac{\pl}{\pl t_1}-\sum_{i\geq 2} it_i
\frac{\pl}{\pl t_{i-1}}  \right)F - nt_1 ~,~ - {\cal
B}_0 F  =
 \left(
2\frac{\pl}{\pl t_2}- \sum_{i\geq 1} it_i
\frac{\pl}{\pl t_i} \right)F - \frac{ n}{2}\sigma_1
\nonumber
\\
 -{\cal B}_1 F &=&
 \left(
2\frac{\pl}{\pl t_3} -\sigma_1
\frac{\pl}{\pl t_1}- \sum_{i\geq 1} it_i
\frac{\pl}{\pl t_{i+1}} \right)F. \eea}
 For later use, take the linear combinations
 \be
{\cal D}_1=-\frac{1}{2}{\cal B}_{-1},~~{\cal
D}_2=-\frac{1}{2}{\cal B}_0,~~{\cal
D}_3=-\frac{1}{2}\left({\cal B}_1+\frac{\sigma_1
}{2}{\cal B}_{-1}\right), \ee
 such that each ${\cal D}_i$ contains the pure term  $\pl F/\pl t_i$, i.e.,
$ {\cal
D}_i F =\frac{\pl F}{\pl t_i}+...$

\subsubsection*{Example 2  (Laguerre $\beta$-integrals)}
 Here, the
weight and the $a_i$ and $b_i$, as in (2.1.1), are
given by
 $$e^{-V}=z^ae^{-z},~~V'=\frac{g}{f}=\frac{z-a}{z}, $$
 $$ a_0=0,~a_1=1,~ b_0=-a,~b_1=1,~\mbox{and all
other}~a_i,b_i=0 .$$ Thus from (2.1.4), the integrals
 \be
I_n= \int_{E^n}|\Dt_n(z)|^{\beta}\prod_{k=1}^n
z_k^{a}e^{-z_k+\sum_{i=1}^{\iy}t_iz_k^i}dz_k
 \ee
satisfy the Virasoro constraints, for $k\geq -1$,
 \be
- {\cal B}_k I_n=- \sum_1^{2r}c_i^{k+2}\frac{\pl}{\pl
c_i}I_n=\left(-~{}^{\beta}\BJ_{k+1,n}^{(2)} -a
~{}^{\beta}\BJ_{k+1,n}^{(1)}+~{}^{\beta}\BJ_{k+2,n}^{(1)}
 \right) I_n.\ee
Introducing the following notation
 \bean \sigma_i&=&=
  (n-\frac{i+1}{2})\beta+i+1-b_0=(n-\frac{i+1}{2})\beta+i+1+a,\\
\eean
 the first three have the form, upon setting $F=F_n=\log
I_n$,
 {\footnotesize \bean
- {\cal B}_{-1} F&=& \left(\frac{\pl}{\pl
t_1}-\sum_{i\geq 1} it_i \frac{\pl}{\pl t_{i}}
 \right) F - \frac{n}{2}(\sigma_1+a)\\
- {\cal B}_0 F   &=&   \left(\frac{\pl}{\pl t_2} -
\sigma_1\frac{\pl}{\pl t_1}-\sum_{i\geq 1} it_i
 \frac{\pl}{\pl t_{i+1}} \right) F \\
  -{\cal B}_1 F  &=& \left(\frac{\pl}{\pl
t_3}-\sigma_2\frac{\pl}{\pl t_2}-\sum_{i\geq 1} it_i
\frac{\pl}{\pl t_{i+2}}- \frac{\beta
}{2}\frac{\pl^2}{\pl t_1^2} \right) F 
 -
\frac{\beta }{2}\left(\frac{\pl F}{\pl t_1} \right)^2.
 \eean}
 Again, replace the operators ${\cal B}_i$ by linear
combinations $ {\cal D}_i$, such that $ {\cal D}_i F
=\frac{\pl F}{\pl t_i}+...$,
 \be
  {\cal D}_1=-{\cal
B}_{-1},~~{\cal D}_2=-{\cal B}_0-\sigma_1{\cal
B}_{-1},~~{\cal D}_3=-{\cal B}_1-
 \sigma_2{\cal B}_0 -\sigma_1\sigma_2{\cal B}_{-1}.
 \ee

\subsubsection*{Example 3  (Jacobi $\beta$-integral)}
 This case is particularly important, because it covers
 not only the first integral, but also
 the second integral in the list (2.0.1),
 used in the problem of random permutations.
The weight and the $a_i$ and $b_i$, as in (2.1.1), are
given by
 $$
\rho (z):=e^{-V}=(1-z)^{a}(1+z)^{b} ,
V'=\frac{g}{f}=\frac{a-b+(a+b)z}{1-z^2}
 $$
 $$
 a_0=1,a_1=0,a_2=-1,b_0=a-b,b_1=a+b
 ,~\mbox{and all
 other}~a_i,b_i=0 .$$
 The integrals
 \be
I_n= \int_{E^n}|\Dt_n(z)|^{\beta}\prod_{k=1}^n
(1-z_k)^{a}(1+z_k)^{b}e^{\sum_{i=1}^{\iy}t_iz_k^i}dz_k
 \ee
satisfy the Virasoro constraints $(k\geq-1)$:
 \bea  - {\cal B}_k I_n&=& -
\sum_1^{2r}c_i^{k+1}(1-c_i^2)\frac{\pl}{\pl
c_i}I_n
 \nonumber\\&=&
\left(~{}^{\beta}\BJ_{k+2,n}^{(2)}-~{}^{\beta}\BJ_{k,n}^{(2)}+
b_0~{}^{\beta}\BJ_{k+1,n}^{(1)}
+b_1~{}^{\beta}\BJ_{k+2,n}^{(1)}\right)I_n. \eea
 Introducing
  $ \sigma_i=(n-\frac{i+1}{2})\beta+i+1+b_1,
$
 the first four have the form ($k=-1,0,1,2$)
{\footnotesize \bea
 - {\cal B}_{-1} F&=& \left(\sigma_1
 \frac{\pl}{\pl
t_1}+\sum_{i\geq 1} it_i \frac{\pl}{\pl t_{i+1}}-
 \sum_{i\geq 2} it_i \frac{\pl}{\pl t_{i-1}}
 \right) F 
 +n(b_0-t_1) \nonumber \\ &&  \nonumber\\
 -{\cal B}_0 F   &=&   \left(
 \sigma_2 \frac{\pl}{\pl t_2}
+  b_0\frac{\pl}{\pl t_1}
 +\sum_{i\geq 1}it_i
 (\frac{\pl}{\pl t_{i+2}}-\frac{\pl}{\pl t_{i}})
   +\frac{\beta}{2}\frac{\pl^2 }{\pl t_1^2}
   \right) F  +\frac{\beta}{2}\left(\frac{\pl F}{\pl
 t_1}\right)^2-\frac{n}{2}(\sigma_1-b_1)\nonumber\\
 &&
\nonumber\\
 - {\cal B}_1 F  &=& \left(\sigma_3
 \frac{\pl}{\pl t_3}+b_0 \frac{\pl}{\pl t_2}
 -(\sigma_1-b_1)\frac{\pl}{\pl t_1}  +\sum_{i\geq 1} it_i (\frac{\pl}{\pl
t_{i+3}}
 -\frac{\pl}{\pl t_{i+1}})
  +\beta \frac{\pl^2}{\pl t_{1}\pl t_2}\right)F
   \nonumber \\&&
 \hspace{8cm} +\beta \frac{\pl F}{\pl t_{1}} \frac{\pl F}{\pl t_{2}}
 \nonumber\\
  - {\cal B}_2 F  &=& \left(\sigma_4
 \frac{\pl}{\pl t_4}+b_0 \frac{\pl}{\pl t_3}
 -(\sigma_2-b_1)\frac{\pl}{\pl t_2} +\sum_{i\geq 1} it_i (\frac{\pl}{\pl
t_{i+4}}
 -\frac{\pl}{\pl t_{i+2}})
 \right.   \nonumber \\
&&\hspace{-0.5cm} \left. +\frac{\beta}{2}(
\frac{\pl^2}{\pl t_{2}^2}-
  \frac{\pl^2}{\pl t_{1}^2}+
  2\frac{\pl^2}{\pl t_{1}\pl t_3})\right)F
   +\frac{\beta}{2}\left(( \frac{\pl F}{\pl t_{2}} )^2-
  ( \frac{\pl F}{\pl t_{1}} )^2
  +2\frac{\pl F}{\pl t_{1}}\frac{\pl F}{\pl t_{3}}\right)
.\nonumber\\ \eea }

\subsection{Double matrix integrals}

Consider now weights of the form \be\rho(x,y)=e^{
\sum_{i,j \geq 1} r_{ij}x^i y^j}\rho(x)
\tilde\rho(y),\ee defined on a product of intervals
 $F_1\times F_2\subset\BR^2$,  with
rational logarithmic derivative
 $$
-\frac{\rho'}{\rho}=\frac{g}{f}=
 \frac{\sum_{i\geq 0} b_i x^i}{\sum_{i\geq 0} a_i x^i}
 ~~\mbox{and}~~
-\frac{\tilde\rho'}{\tilde\rho}=\frac{\tilde g}{\tilde
f}=
 \frac{\sum_{i\geq 0} \tilde b_i y^i}{\sum_{i\geq 0}
  \tilde a_i y^i}
 ,$$
  satisfying
\be
\lim_{x\rightarrow\pl  F}  f(x) \rho(x)x^k=
\lim_{y\rightarrow\pl \tilde F} \tilde f(y)\tilde
\rho(y)y^k=0 ~~\mbox{for all}~k\geq 0.
 \ee
 Consider subsets of the form \be E=E_1\times
E_2:=\cup^r_{i=1}[c_{2i-1},c_{2i}]\times
\cup^s_{i=1}[\tilde c_{2i-1},\tilde c_{2i}]\subset F_1
\times  F_2\subset\BR^2. \ee
 A natural deformation of the second integral
 in the list (2.0.1) is given by the following
 integrals:
\be
I_n(t,s,r;E)=\int\!\!\int_{E^n} \Delta_n(x)\Delta_n(y)
\prod^n_{k=1} e^{\sum^{\iy}_{i=1}(t_ix_k^i-s_iy^i_k)}
\rho(x_k,y_k)dx_kdy_k
 \ee
In the theorem below, $\BJ_{k,n}^{(i)}$ and
$\tilde\BJ_{k,n}^{(i)}$ are vectors of operators,
whose components are given by the operators (2.1.10)
for $\beta =1$; i.e., $$ \BJ_{k,n}^{(i)}(t)=\left.
{}^{\beta}\BJ_{k,n}^{(i)}(t)\right|_{\beta=1},
 \tilde\BJ_{k,n}^{(i)}(s):
  = \left.{}^{\beta}
\BJ_{k,n}^{(i)}(t)\right|_{\beta=1, ~t\mapsto -s} ;$$
thus, from (2.1.10) and (2.1.11), one finds:
 \be
\BJ_{k,n}^{(2)}(t)
 =\frac{1}{2}\left(J_k^{(2)}(t)
+ (2n+k+1)J_k^{(1)}(t)
+n(n+1)J_k^{(0)}\right)
,
 \ee
satisfying the Heisenberg and Virasoro relations
(2.1.8), with {\em central charge} $c=-2$ and
$c^{\prime}=1/2$.

The $a_i,\tilde a_i,b_i,\tilde b_i,c_i,\tilde c_i,
r_{ij}$ given by (2.2.1), (2.2.2) and (2.2.3) define
differential operators:
  \bea \DR_{k,n}:&=&
 \sum_1^{2r} c_i^{k+1}f(c_i)\frac{\pl}{\pl c_i}-
\sum_{i\geq 0}\Bigl( a_i (\BJ_{k+i,n}^{(2)}+\sum
 _{m,\ell \geq 1}mr_{m\ell}\frac{\pl}{\pl r_{m+k+i,\ell}})-b_i
 \BJ_{k+i+1,n}^{(1)}
\Bigr) \nonumber\\
 \tilde\DR_{k,n}:&=&
  \sum_1^{2r} \tilde c_i^{k+1}\tilde
 f(\tilde c_i)\frac{\pl}{\pl \tilde c_i}-
\sum_{i\geq 0}\Bigl( \tilde a_i
(\tilde\BJ_{k+i,n}^{(2)}
 +\sum
 _{m,\ell \geq 1}\ell r_{m\ell}\frac{\pl}{\pl r_{m,\ell +k+i}})
 -\tilde b_i
\tilde \BJ_{k+i+1,n}^{(1)}
\Bigr).\nonumber\\
 \eea

\begin{theorem}{\em (Adler-van Moerbeke \cite{AvM2,AvM21})}
 Given $\rho(x,y)$ as in (2.2.1), the integrals
 \be
 I_n(t,s,r;E):=\int\!\!\!\int_{E^n}
\Delta_n(x)\Delta_n(y) \prod^n_{k=1}
e^{\sum^{\iy}_{i=1}(t_ix_k^i-s_iy^i_k)}
\rho(x_k,y_k)dx_kdy_k
 \ee
 satisfy two families of Virasoro equations
  for $k\geq -1$:
  \be
\DR_{k,n} I_n(t,s,r;E)=0
~~~\mbox{and}~~~\tilde\DR_{k,n} I_n (t,s,r;E)=0.
 \ee

\end{theorem}

\vspace{0.5cm}

 The proof of this statement is very
similar to the one for $\beta$-integrals. Namely,
define the vector vertex operator,
\be
\BX_{12}(t,s;u,v)=\Lb^{-1}e^{\sum_1^{\iy}(t_i u^i-s_i
v^{i})} e^{-\sum_1^{\iy}(\frac{ u^{-i}}{i}
\frac{\pl}{\pl t_i}-
 \frac{ v^{-i}}{i}
\frac{\pl}{\pl s_i})}\chi(uv), \ee
 which, as a consequence of Proposition 2.2, interacts with
 the operators
 $\BJ_k^{(i)}(t)=
 \left(\BJ_{k,n}^{(i)}(t,n)\right)_{n\in \BZ}$,
 as follows:
  \bea
  u^{k}\BX_{12}(t,s;u,v)&=& \left[
\BJ_k^{(1)}(t), \BX_{12}(t,s;u,v)\right] \nonumber\\
\frac{\pl}{\pl u}u^{k+1}\BX_{12}(t,s;u,v)&=& \left[
\BJ_k^{(2)}(t), \BX_{12}(t,s;u,v)\right]. \eea
 A similar statement can
be made, upon replacing the operators $u^k$ and
$\frac{\pl}{\pl u}u^{k+1}$ by $v^k$ and
$\frac{\pl}{\pl v}v^{k+1}$, and upon using the
$\tilde\BJ_k^{(i)}(s)$'s.

Finally, one checks that the integral vertex operator
\be
\BY(t,s; \rho_E):=\int\!\!\int_{E} dx dy \rho(x,y)
\BX_{12}(t,s;x,y) \ee
 commutes with the two vectors of differential operators
 ${\cal D}_k=({\cal D}_{k,n})_{n\in \BZ}$, as in  (2.2.6):
  $$
  \Bigl[  \DR_k, \BY(t,s;\rho_E) \Bigr]= \left[  \tilde\DR_k, \BY(t,s;\rho_E) \right]=0,
   $$
    and that the vector $I=(I_0=1,I_1,...)$ of
    integrals (2.2.7) is a fixed point for
$\BY(t,s; \rho_E)$, $$
 \BY(t,s;\rho_E)
 I(t,s,r;E)= I(t,s,r;E).
 $$
Then, as before, the proof of Theorem 2.6 hinges
ultimately on the fact that $\DR_{k,0}$ annihilates
$I_0=1$.

\subsection{Integrals over the unit circle}

We now deal with the fourth type of integral in the
list (2.0.1), which we deform, this time, by inserting
two sequences of times $t_1,t_2,...$ and
$s_1,s_2,...$. The following theorem holds:

\begin{theorem} {\em (Adler-van Moerbeke \cite{AvM4})}
The multiple integrals 
 over the unit circle $S^1$,
\be
I_n(t,s)= \int_{(S^1)^{n}}|\Dt_n(z)|^{2}
 \prod_{k=1}^n
e^{\sum_1^{\iy}(t_i z_k^i-s_iz_k^{-i})}
 \frac{dz_k}{2\pi i z_k},~~~n>0,\ee with $I_0=1$,
 satisfy an SL(2,$\BZ$)-algebra of Virasoro constraints:
\be
\DR^{\theta}_{k,n}I_n(t,s)=0,
 ~~~\mbox{for}~ \left\{\begin{array}{l} k
=-1 ,~\theta= 0\\ k=0 , ~~\theta~~ {arbitrary}\\
  k=1 ,~
\theta= 1\end{array}\right\} \mbox{\,\,only,
 }
 \ee
where the operators $\DR^{\theta}_{k,n}:=
 \DR^{\theta}_{k,n}(t,s,n),~~k \in \BZ, ~n\geq 0$ are given
by
 \be
 {\cal D}^{\theta}_{k,n}:=\BJ_{k,n}^{(2)}(t,n)-
  \BJ_{-k,n}^{(2)}(-s,n) -k\left( \theta
\BJ_{k,n}^{(1)}(t,n)+(1-\theta)\BJ_{-k,n}^{(1)}(-s,n)
   \right),
 \ee
 with $\BJ_{k,n}^{(i)}(t,n):=
 \left.{}^{\beta}\BJ_{k,n}^{(i)}(t,n)\right|_{\beta=1} $, as in (2.1.11).

\end{theorem}

The explicit expressions are {\footnotesize
\begin{eqnarray} {\cal
D}_{-1}I_n&=&\left(\sum_{i\geq
1}(i+1)t_{i+1}\frac{\pl}{\pl t_{i}}-\sum_{i\geq
2}(i-1)s_{i-1}\frac{\pl}{\pl
s_{i}}+n\left(t_1+\frac{\pl}{\pl
s_{1}}\right)\right)I_n=0\nonumber\\ {\cal
D}_{0}I_n&=&\sum_{i\geq 1}\left(it_{i}\frac{\pl}{\pl
t_{i}}-is_{i}\frac{\pl}{\pl s_{i}}\right) I_n=0 \\
{\cal D}_{1}I_n&=&\left(-\sum_{i\geq
1}(i+1)s_{i+1}\frac{\pl}{\pl s_{i}}+\sum_{i\geq
2}(i-1)t_{i-1}\frac{\pl}{\pl
t_{i}}+n\left(s_1+\frac{\pl}{\pl t_1}
\right)\right)I_n=0.\nonumber
\end{eqnarray}  }

 Here the key vertex
operator is a reduction of $\BX_{12}(t,s;u,v)$,
defined in the previous section (formula (2.2.9)). For
all $k \in \BZ$, the vector of operators
 ${\cal D}^{\theta}_{k}(t,s)=
 \left({\cal D}^{\theta}_{k,n}(t,s,n)\right)_{n\in \BZ}$ form a realization of
 the first order differential operators $\frac{d}{du}u^{k+1}$,
  using the vertex operator
 $\BX_{12}(t,s;u,u^{-1})$, namely
\be
\frac{d}{du}u^{k+1} \frac{\BX_{12}(t,s;u,u^{-1})}{u}
 =\left[ {\cal D}^{\theta}_{k}(t,s),
  \frac{\BX_{12}(t,s;u,u^{-1})}{u}\right]
  . \ee
Indeed,

  \noindent$\displaystyle{
 u\frac{d}{du}u^k \BX_{12}(t,s;u,u^{-1})}$
\begin{eqnarray*}
 &=&\left. \left(\frac{\pl}{\pl u}u^{k+1}-\frac{\pl}{\pl v}v^{1-k}
 -k\theta u^k-k(1-\theta)v^{-k}\right)\BX_{12}(t,s;u,v) \right|_{v=-u}\\
  &=&\left[ \BJ_{k}^{(2)}(t)-\BJ_{-k}^{(2)}(-s)
-k\left( \theta
\BJ_{k}^{(1)}(t)+(1-\theta)\BJ_{k}^{(1)}(-s)
   \right),\BX_{12}(t,s;u,-u)\right]\\
    &=&
    \left[\DR_k^{\theta}(t,s),\BX_{12}(t,s;u,u^{-1}\right].
   \end{eqnarray*}
The ${\cal D}^{\theta}_{k}:={\cal
D}^{\theta}_{k}(t,s)$ satisfy Virasoro relations with
central charge $=0$,
\be
 \left[{\cal D}^{\theta}_{k},{\cal D}^{\theta}_{\ell}\right]=
 (k-\ell){\cal D}^{\theta}_{k+\ell},
 \ee
and, from (2.3.5) the following commutation relation
holds:
\be
 \left[{\cal D}^{\theta}_{k}(t,s),
 \BY(t,s)\right]=0,~~\mbox{with}~~ {\BY}(t,s):=\int_{S^1}\frac{du}{2\pi
iu}\BX_{12}(t,s;u,u^{-1}).
 \ee

The point is that the column vector
$I(t,s)=(I_0,I_1,...)$ of integrals (2.3.1), is a
fixed point for ${\BY}_{}(t,s)$:
\be
({\BY}_{}(t, s)I)_n=I_n, ~n\geq 1,\ee
 which is shown in a way, similar to Proposition 2.5.

{\noindent{\sl Proof of Theorem 2.7}: } Here again the
proof is similar to the one of Theorem 2.1.  Taking
the $n$th component and the $n$th power of $\BY(t,s)$,
with $n\geq 1$, and noticing from the explicit
formulae (2.3.4) that $\left({\cal
D}^{\theta}_k(t,s)I\right)_{0}=0$, we have, by means
of a calculation similar to the proof of Theorem 2.1,
  \bean
0&=&\left([{\cal D}^{\theta}_{k},\BY(t,s)^n]I\right)_n
\\
 &=&\left({\cal D}^{\theta}_{k} \BY(t,s)^n
I-\BY(t,s)^n{\cal D}^{\theta}_{k}I\right)_n
 \\
 &=&\left({\cal D}^{\theta}_{k}I-\BY(t,s)^n{\cal
D}^{\theta}_{k}I\right)_n
 =\left({\cal D}^{\theta}_{k} I\right) _{n} .
 \eean

\vspace{-1cm} \qed


\section{Integrable systems and associated matrix integrals}
\subsection {Toda lattice and Hermitian matrix
integrals}

\subsubsection{Toda lattice, factorization of
 symmetric matrices and orthogonal polynomials}

Given a weight $\rho(z)=e^{-V(z)}$ defined as in
(2.1.1), the inner-product over $E\subseteq \BR$,
\be
\la
f,g\ra_t=\int_{E}f(z)g(z)\rho_t(z)dz,~~~\mbox{with}
~\rho_t(z):=e^{\sum t_iz^i} \rho(z), \ee leads to a
moment matrix
\be
m_{n}(t)=(\mu_{ij}(t))_{0\leq i,j<n}=(\la
z^i,z^j\ra_t)_{0\leq i,j<n}, \ee which is a {\em
H\"ankel matrix}\footnote{ H\"ankel means $\mu_{ij}$
depends on $i+j$ only}, thus symmetric. This is
tantamount to $\Lambda m_{\iy}=m_{\iy}\Lambda^{\top}$,
where $\Lb$ denotes the shift matrix; see footnote 11
. As easily seen, the semi-infinite moment matrix
$m_{\iy}$ evolves in $t$ according to the equations
\be
\frac{\pl\mu_{ij}}{\pl
t_{k}}=\mu_{i+k,j},\mbox{\,\,and thus\,\,} \frac{\pl
m_{\iy}}{\pl
t_{k}}=\Lambda^km_{\iy}.~~~~~\left(\begin{array}{l}
\mbox{commuting}\\ \mbox{vector fields}
\end{array}
\right) \ee Another important ingredient is the
factorization of $m_{\iy}$ into a lower- times an
upper-triangular matrix\footnote{This factorization is
possible for those $t$'s for which
$
\tau_n(t):=\det m_n(t)\neq 0 $ for all $n>0$.} \bean
&& m_{\iy}(t)=S(t)^{-1}S(t)^{\top -1},\quad~~~~~~
S(t)=\mbox{\,lower triangular with}\\ &&
\hspace{7cm}\mbox{non-zero diagonal elements.} \eean
The main ideas of the following theorem can be found
in \cite{AvM1, AvM22}. Remember $c=(c_1,...,c_{2r})$
denotes the boundary points of the set $E$. $dM$
refers to properly normalized Haar measure on ${\cal
H}_n$.

\begin{theorem} The determinants of the moment matrices
 \bea
  \tau_n(t,c): =\det m_n(t,c)
&=&\frac{1}{n!}\int_{E^n}\Delta^{2}_n(z)
\prod^n_{k=1}\rho_t(z_k)dz_k \nonumber \\ & =&
\int_{{\cal H}_n(E)}
e^{tr(-V(M)+\sum_1^{\iy}t_iM^i)}dM, \eea
  satisfy
\begin{description}

  \item[(i)] \underline{Virasoro constraints} (2.1.7)
for $\beta =2$,
 \be
\left(-\sum_1^{2r} c_i^{k+1}f(c_i)\frac{\pl}{\pl c_i}+
\sum_{i\geq 0}\left( a_i~ \BJ_{k+i,n}^{(2)}-b_i ~
 \BJ_{k+i+1,n}^{(1)}\right) \right)\tau_n(t,c) =
0. \ee

 \item[(ii)] The \underline{KP-hierarchy}\footnote{Given a polynomial
 $p(t_1,t_2,...)$, define the
customary Hirota symbol $p(\pl_t)f\circ g:=
p(\frac{\pl}{\pl y_1},\frac{\pl}{\pl
y_2},...)f(t+y)g(t-y) \Bigl|_{y=0}$. The ${\bf
s}_{\ell}$'s are the elementary Schur polynomials
$e^{\sum^{\iy}_{1}t_iz^i}:=\sum_{i\geq 0} {\bf
s}_i(t)z^i$ and ${\bf s}_{\ell}(\tilde \pl):={\bf
s}_{\ell}(\frac{\pl}{\pl
t_1},\frac{1}{2}\frac{\pl}{\pl
t_2},\ldots).$}($k=0,1,2,\ldots$)
  $$
\left({\bf s} _{k+4}\bigl(\frac{\pl}{\pl
t_1},\frac{1}{2}\frac{\pl}{\pl
t_2},\frac{1}{3}\frac{\pl}{\pl
t_3},\ldots\bigr)-\frac{1}{2}\frac{\pl^2}{\pl t_1\pl
t_{k+3}}\right)\tau_n \circ\tau_n=0, $$
 of
which the first equation reads: \be \hspace{-1cm}
\left(\left(\frac{\pl}{\pl t_1}
\right)^4+3\left(\frac{\pl}{\pl
t_2}\right)^2-4\frac{\pl^2}{\pl t_1 \pl
t_3}\right)\log\tau_n+6\left(\frac{\pl^2}{\pl
t^2_1}\log\tau_n \right)^2=0. \ee

\item [(iii)] The \underline{standard Toda lattice}, i.e., the
symmetric tridiagonal matrix
 \be
  L(t):=S(t)\Lb S(t)^{-1}=
\left(\begin{tabular}{lllll} $\frac{\pl}{\pl t_1}\log
\frac{\tau_1}{\tau_0}$ &
 $ \left(\frac{\tau_{0}\tau_{2}}{\tau_{1}^2}\right)^{1/2}$
 & ~~ $0$      &    \\
$\left(\frac{\tau_{0}\tau_{2}}{\tau_{1}^2}\right)^{1/2}$&
$\frac{\pl}{\pl t_1}\log \frac{\tau_2}{\tau_1}$ &
$\left(\frac{\tau_{1}\tau_{3}}{\tau_{2}^2}\right)^{1/2}$&
\\~~~$0$&$\left(\frac{\tau_{1}\tau_{3}}{\tau_{2}^2}\right)^{1/2}$
&
 $\frac{\pl}{\pl t_1}\log \frac{\tau_3}{\tau_2}$&
   \\
 & &  &      $\ddots$\\
\end{tabular}
\right)
 \ee
satisfies the commuting equations\footnote{$()_{\frak
s}$ means: take the skew-symmetric part of $()$; for
more details, see subsection 3.1.2.}
 \be\frac{\pl L}{\pl
t_k}=\left[\frac{1}{2}(L^k)_{\frak s},L\right].
 \ee


\item [(iv)]  \underline{Eigenvectors of $L$}:
The tridiagonal matrix $L$ admits two independent
eigenvectors:

\begin{itemize}
  \item 
 $p(t;z)=(p_n(t;z))_{n\geq 0}$ satisfying
$(L(t)p(t;z))_n=zp_n(t;z),~n\geq 0$, where $p_n(t;z)$
are $n$th degree polynomials  in $z$, depending on $t
\in \BC^{\iy}$, orthonormal with respect to the
$t$-dependent inner product\footnote{The explicit
dependence on the boundary points $c$ will be omitted
in this point (iv). } (3.1.1) $$ \la
p_k(t;z),p_{\ell}(t;z)\ra_t =\delta_{k\ell}; $$
 they are eigenvectors of $L$, i.e., $L(t)p(t;z)=zp(t;z)$, and
 enjoy the following representations:
 $(\chi(z)=(1,z,z^2,...)^{\top})$
 \bea 
p_n(t;z)&:=&(S(t)\chi(z))_n\nonumber\\&=&
\frac{1}{\sqrt{\tau_n(t)\tau_{n+1}( t)}} \det\left(
\begin{array}{lll|l}
  & & &1\\
  &m_n(t)& &z\\
  & & &\vdots\\
 \hline
 \mu_{n,0}&\ldots&\mu_{n,n-1}&z^n
\end{array}
\right) \nonumber
 \\
 &=&z^nh_{n}^{-1/2}\frac{\tau_n(t-[z^{-1}])}{\tau_{n}
(t)}, ~~~~h_{n}:=\frac{\tau_{n+1}(t)}{\tau_{n}(t)}
\eea 
 \item
$
 q(t,z)=(q_n(t;z))_{n\geq 0}$, with $q_n(t;z):=
 z\int_{\BR^n}\frac{p_n(t;u)}{z-u}\rho_t(u)du
  $, satisfying $(L(t)q(t;z))_n=zq_n(t;z)
 ,~~n\geq 1$; $q_n(t;z)$ enjoys
the following
 representations:
 \bea
q_n(t;z):=z\int_{\BR^n}\frac{p_n(t;u)}{z-u}\rho_t(u)du
&=&\left(S^{\top -1}(t)\chi(z^{-1})\right)_n
\nonumber\\
&=&\left(S(t)m_{\iy}(t)\chi(z^{-1})\right)_n\nonumber\\
&=&z^{-n}h_n^{-1/2}
\frac{\tau_{n+1}(t+[z^{-1}])}{\tau_n(t)}.\nonumber\\
\eea
 \end{itemize}
\end{description}

\end{theorem}

  In the case $\beta=2$, the Virasoro generators
  (2.1.11) take on a particularly
  elegant form, namely for ${n \geq 0}$,
  \bean \BJ_{k,n}^{(2)}(t)&=&\sum_{i+j=k}
   :~ \BJ_{i,n}^{(1)}(t) ~
\BJ_{j,n}^{(1)}(t):~=~  J_k^{(2)}(t) + 2n J_k^{(1)}(t)
+ n^2\dt_{0k}\\
 \BJ_{k,n}^{(1)}(t) &=& ~  J_k^{(1)}(t)+
 n\dt_{0k} ,
  \eean
 with\footnote{The expression $J_k^{(1)}=0$ for $k=0$. }
\be
J_k^{(1)}=\frac{\pl}{\pl t_k}+\frac{1}{2}(-k)t_{-k}~,~
J^{(2)}_{k}=\sum_{i+j=k}\frac{\pl^2}{\pl
 t_{i}\pl t_{j}}+ \sum_{-i+j=k}it_{i}\frac{\pl}{\pl
 t_{j}}+\frac{1}{4}\sum_{-i-j=k}it_{i}jt_{j}.
 \ee
Statement (i) is already contained in section 2,
whereas statement (ii) will be established in
subsection 3.1.2, using elementary methods.

\subsubsection{Sketch of Proof}

\noindent {\bf Orthogonal polynomials and
$\tau$-function representation:} The representation
(3.1.4) of the determinants of moment matrices as
integrals follows immediately from the fact that the
square of a Vandermonde determinant can be represented
as a sum of determinants $$
\Delta^2(u_1,\ldots,u_n)=\sum_{\sg\in
S_n}\det\left(u_{\sg(k)}^{\ell +k-2}\right)_{1\leq
k,\ell\leq n}. $$
  Indeed,
 \bean n! \tau_n(t)&=& n! \det m_n(t)\\
 &=& \sum_{\sigma \in S_n}\det
 \left( \int_E z_{\sigma(k)}^{\ell +k-2} \rho_t(z_{\sigma(k)})
 dz_{\sigma(k)}  \right)_{1\leq
k,\ell\leq n}\\
 &=&   \sum_{\sigma \in S_n} \int_{E^n}  \det
 \left(z_{\sigma(k)}^{\ell +k-2}   \right)_{1\leq
k,\ell\leq n}\rho_t(z_{\sigma(k)})
 dz_{\sigma(k)}\\
 &=& \int_{E^n}\Delta^{2}_n(z)
\prod^n_{k=1}\rho_t(z_k)dz_k,
  \eean
whereas the representation (3.1.4) in terms of
integrals over Hermitian matrices follows from section
1.1.

The Borel factorization of $m_{\iy}$ is responsible
for the orthonormality of the polynomials
$p_n(t;z)=(S(t)\chi(z))_n$; indeed,
 $$ \la
p_k(t;z),p_{\ell}(t;z)\ra_{0\leq k,\ell
<\iy}=\int_{E}S\chi(z)(S\chi(z))^{\top}\rho_t(z)dz=S
m_{\iy}S^{\top}=I. $$ Note that
$S\chi(z)(S\chi(z))^{\top}$ should be viewed as a
semi-infinite matrix obtained by multiplying a
semi-infinite column and row. The determinantal
representation (3.1.9) follows at once from noticing
that $\la p_n(t;z), z^k \ra=0 $ for $0\leq k \leq
n-1,$ because taking that inner-product produces two
identical columns in the matrix thus obtained. From
the same representation (3.1.9), one has
$p_n(t;z)=h_n^{-1}z^n+...$ , where
$h_n:=(\tau_{n+1}/\tau_n(t))^{1/2}$.

The ``Sato" representation (3.1.9) of $p_n(t;z)$ in
terms of the determinant $\tau_n(t)$ of the moment
matrix can be shown by first proving the Heine
representation of the orthogonal polynomials, which
goes as follows:

 \noindent $h_n p_n(t;z)$
\bean
 &=& \frac{1}{\tau_n}\det\left(
\begin{array}{lll|l}
  & & &1\\
  &m_n(t)& &z\\
  & & &\vdots\\
 \hline
 \mu_{n,0}&\ldots&\mu_{n,n-1}&z^n
\end{array}
\right)\\
  &=&\frac{1}{\tau_n} \int_{E^n}\det\left(
\begin{array}{llll|l}
 u_1^0 &u_2^1 &\cdots &u_n^{n-1} &1\\
  u_1^1&  u_2^2&\cdots  & u_n^n& z\\
  \vdots& & &\vdots&\vdots\\
  u_1^{n-1}& u_2^{n} &\cdots  &u_n^{2n-2}& z^{n-1}\\
 \hline
 u_1^n& u_2^{n+1}& \cdots & u_n^{2n-1} &z^n
\end{array}
\right) \prod_1^n \rho_t(u_i)du_i \\
&=&\frac{1}{\tau_n} \int_{E^n}\det\left(
\begin{array}{llll|l}
 u_1^0 &u_2^0 &\cdots &u_n^{0} &1\\
  u_1^1&  u_2^1&\cdots  & u_n^1& z\\
  \vdots& & &\vdots&\vdots\\
  u_1^{n-1}& u_2^{n-1} &\cdots  &u_n^{n-1}& z^{n-1}\\
 \hline
 u_1^n& u_2^{n}& \cdots & u_n^{n} &z^n
\end{array}
\right) u_1^0 u_2^1 \cdots u_n^{n-1} \prod_1^n
\rho_t(u_i)du_i\\
 &=&\frac{1}{\tau_n} \int_{E^n}\det\left(
\begin{array}{llll|l}
 u_{\sigma(1)}^0 &u_{\sigma(2)}^0 &\cdots &u_{\sigma(n)}^{0} &1\\
  u_{\sigma(1)}^1&  u_{\sigma(2)}^1&\cdots  & u_{\sigma(n)}^1& z\\
  \vdots& & &\vdots&\vdots\\
  u_{\sigma(1)}^{n-1}& u_{\sigma(2)}^{n-1} &\cdots
   &u_{\sigma(n)}^{n-1}& z^{n-1}\\
 \hline
 u_{\sigma(1)}^n& u_{\sigma(2)}^{n}& \cdots & u_{\sigma(n)}^{n} &z^n
\end{array}
\right)\\ && \hspace{3cm} u_{\sigma(1)}^0
u_{\sigma(2)}^1 \cdots u_{\sigma(n)}^{n-1}
 \prod_1^n \rho_t(u_{\sigma (i)}) du_{\sigma (i)}
~,\\&& \hspace{5cm} ~~~~\mbox{for any permutation
$\sigma\in S_n$}\\
  &=&\frac{1}{\tau_n} \int_{E^n}\det\left(
\begin{array}{llll|l}
 u_1^0 &u_2^0 &\cdots &u_n^{0} &1\\
  u_1^1&  u_2^1&\cdots  & u_n^1& z\\
  \vdots& & &\vdots&\vdots\\
  u_1^{n-1}& u_2^{n-1} &\cdots  &u_n^{n-1}& z^{n-1}\\
 \hline
 u_1^n& u_2^{n}& \cdots & u_n^{n} &z^n
\end{array}
\right)  \\ && \hspace{3cm}(-1)^{\sigma}
u_{\sigma(1)}^0 u_{\sigma(2)}^1 \cdots
u_{\sigma(n)}^{n-1}  \prod_1^n \rho_t(u_{i}) du_{i}\\
 &=&\frac{1}{n!\tau_n}\int_{E^n}\Delta_n^2(u)
  \prod^n_{k=1}(z-u_k)\rho_t(u_k)du_k,~~\mbox{
  upon summing over all $\sigma$.}
\eean Therefore, using again the representation of
$\Delta^2(z)$ as a sum of determinants, Heine's
formula leads to
 \bean \lefteqn{h_np_n(t,z)}\\
&=&\frac{z^n}{n!\tau_n}\int_{E^n}\sum_{\sg\in
S_n}\det\left(u_{\sg(k)}^{\ell +k-2}\right)_{1\leq
k,\ell\leq
n}\prod^n_{k=1}\left(1-\frac{u_{\sg(k)}}{z}\right)
\rho_t(u_{\sg(k)})du_{\sg(k)},\\
&=&\frac{z^n}{n!\tau_n}\int_{E^n}\sum_{\sg\in
S_n}\det\left(u_{\sg(k)}^{\ell
+k-2}-\frac{1}{z}u_{\sg(k)}^{\ell +k-1}\right)_{1\leq
k,\ell\leq n }\rho_t(u_{\sg(k)})du_{\sg(k)}\\
 &=&\frac{z^n}{\tau_n}\det\left(\mu_{ij}-\frac{1}{z}\mu_{i,j+1}\right)_{0\leq
i,j\leq n-1}\\
 &=& \frac{z^n}{\tau_n}
\det\left(\mu_{ij}(t-[z^{-1}])\right)_{0\leq i,j\leq
n-1} \\ &=& z^n\frac{ \tau_n(t-[z^{-1}])}
 {\tau_n(t)}, \eean \vspace{-1.7cm}\be  \ee
 invoking
the fact that \bean \mu_{ij}(t-[z^{-1}])=\int
u^{i+j}e^{\sum_1^{\iy}\left(t_i-\frac{z^{-i}}{i}\right)u^i}\rho(u)du&=&\int
u^{i+j}\left(1-\frac{u}{z}\right)\rho(u)du\\
&=&\mu_{i+j}-\frac{1}{z}\mu_{i+j+1}. \eean


Formula (3.1.10) follows from computing on the one
hand $S(t)m_{\iy}\chi(z)$ using the explicit moments
$\mu_{ij}$, together with (3.1.12), and on the other
hand the equivalent expression $S^{\top
-1}(t)\chi(z^{-1})$. Indeed, using
$(S(t)\chi(z))_n=p_n(t;z)=\sum_{0}^n p_{nk}(t)z^k$,
\begin{eqnarray*}
\sum_{j\geq 0} \left(  S m_{\iy} \right)_{nj} z^{-j}
 &=&\sum_{j\geq 0} z^{-j} \sum_{\ell \geq 0}
 p_{n\ell}(t)\mu_{\ell j} \\
 &=&  \sum_{j\geq 0} z^{-j}
\sum_{\ell \geq 0} p_{n\ell}(t)\int_{E}
u^{\ell+j}\rho_t(u)du
\\ &=&  \int_{E} \sum_{\ell \geq 0} p_{n\ell}(t)u^{\ell}
\sum_{j\geq 0}\left(\frac{u}{z} \right)^j\rho_t(u)du\\
&=&z  \int_{E} \frac{ p_n(t,u)\rho_t(u)}{z-u}du.\\
\end{eqnarray*}
 Mimicking computation (3.1.12), one shows
\begin{eqnarray*}
h_n\sum_{j\geq 0} \left(S^{\top -1}(t) \right)_{nj}
z^{-j} &=& \frac{\tau_{n+1}(t+[z^{-1}])}{\tau_n(t) }
z^{-n} ,
\end{eqnarray*}
from which (3.1.10) follows, upon using $S
m_{\iy}=S^{\top -1}$. Details of this and subsequent
derivation can be found in \cite{AvM22,AvM3}.

 {\em The vectors $p$ and $q$ are eigenvectors
of $L$}. Indeed, remembering
$\chi(z)=(1,z,z^2,...)^{\top}$, and the shift
$(\Lambda v)_n=v_{n+1}$, we have
 $$ \Lb \chi(z)=z\chi(z)~~\mbox{and}~~
  \Lb^{\top}
  \chi(z^{-1})=z\chi(z^{-1})-ze_1,~~\mbox{with}~e_1=(1,0,0,...)^{\top}.
  $$
  Therefore, $p(z)=S\chi(z)$ and $q(z)=S^{\top
  -1}\chi(z^{-1})$ are eigenvectors, in the sense
  \bean
  Lp&=&S\Lb S^{-1}S\chi(z)=zS\chi(z)=zp  \\
  L^{\top}q&=& S^{\top -1}\Lb^{\top} S^{\top}
   S^{\top -1} \chi(z^{-1})=zS^{\top -1}\chi(z^{-1})
   -zS^{\top -1}e_1=zq-zS^{\top -1}e_1 .
   \eean
Then, using $L=L^{\top}$, one is lead to
$$((L-zI)p)_n=0,~~\mbox{for}~ n\geq 0~~~\mbox{and}~~~
 ((L-zI)q)_n=0,~~\mbox{for}~~ n\geq 1.$$

\noindent {\bf Toda lattice and Lie algebra
splitting:} The Lie algebra splitting of semi-infinite
matrices and the corresponding projections (used in
(3.1.8)), denoted by $(~)_{\frak s}$ and $(~)_{\frak
b}$ are defined as follows: $$ gl(\iy)={\frak s}
\oplus {\frak b}\left\{\begin{array}{l} {\frak s}
=\{\mbox{skew-symmetric matrices}\}\\ {\frak
b}=\{\mbox{lower-triangular matrices}\}
\end{array}
\right. .$$
 Conjugating the shift matrix $\Lambda$ by
$S(t)$ yields a matrix \bean L(t)&=&S(t)\Lb
S(t)^{-1}\\ &=&S \Lb S^{-1}S^{\top -1}S^{\top}\\
&=&S\Lb m_{\iy}S^{\top},\mbox{\,\,using (3.1.3)},\\
&=&Sm_{\iy}\Lb^{\top}S^{\top},\mbox{\,\,using $\Lb
m_{\iy}=m_{\iy}\Lb^{\top}$},\\ &=&S(S^{-1}S^{\top
-1})\Lb^{\top} S^{\top},\mbox{\,\,using again
(3.1.3)},\\ &=&(S\Lb S^{-1})^{\top}=L(t)^{\top}, \eean
which is symmetric and thus tridiagonal. Moreover,
from (3.1.3) one computes
  \bean
0&=&S\left(\Lb^km_{\iy}-\frac{\pl m_{\iy}}{\pl
t_k}\right)S^{\top}\\ &=&S\Lb^kS^{-1}-S\frac{\pl}{\pl
t_k}(S^{-1}S^{\top -1})S^{\top}\\ &=&L^k+\frac{\pl
S}{\pl t_k} S^{-1}+S^{\top -1}\frac{\pl S^{\top}}{\pl
t_k}. \eean
 Upon taking the $()_-$ and $()_0$ parts of this equation
  ($A_-$ means the lower-triangular part of the
matrix $A$, including the diagonal and $A_0$ the
diagonal part) leads to
$$
  (L^k)_-+\frac{\pl S}{\pl
t_k} S^{-1}+\left(S^{\top -1}\frac{\pl S^{\top}}{\pl
t_k}\right)_0=0~\mbox{and}~\left(\frac{\pl S}{\pl t_k}
S^{-1}\right)_0=-\frac{1}{2}(L^k)_0.
 $$
  Upon observing that for any
symmetric matrix
  $$
 \left(
\begin{array}{ll}
 a&c\\
  c&b
\end{array}
\right)_{\frak b}=
 \left(
\begin{array}{ll}
 a&0\\
  2c&b
\end{array}
\right)
 =2
 \left(
\begin{array}{ll}
 a&c\\
  c&b
\end{array}
\right)_- - \left(
\begin{array}{ll}
 a&c\\
  c&b
\end{array}
\right)_0,$$
 it follows that
 the matrices $L(t)$, $S(t)$ 
and the vector $p(t;z)=(p_n(t;z))_{n\geq
0}=S(t)\chi(z)$ satisfy the (commuting) differential
equations and the eigenvalue problem
\be
\frac{\pl S}{\pl t_k}=-\frac{1}{2}(L^k)_{\frak
b}S,~~~~L(t)p(t;z)=zp(t;z),  \ee and thus
 $$
 \frac{\pl L}{\pl
t_k}=-\left[\frac{1}{2}(L^k)_{\frak b},L\right],\quad
\frac{\pl p}{\pl t_k} =-\frac{1}{2}(L^k)_{\frak b} p.
$$ \hfill (Standard Toda lattice)

\noindent {\bf The bilinear identity:} The functions
$\tau_n(t)$ satisfy the following identity, for $n\geq
m+1, ~ t,t' \in \BC$, where one integrates along a
small circle about $\iy$,
\be
\oint_{z=\iy}\tau_n(t-[z^{-1}])\tau_{m+1}(t'+[z^{-1}])e^{\sum
(t_i-t'_i)z^i}z^{n-m-1}dz=0. \ee

 An elementary proof can be given by expressing the left
 hand side of (3.1.14), in terms of $p_n(t;z)$ and $p_m(t,z)$,
 using (3.1.9) and (3.1.10). One uses below the following
  identity (see \cite{AvM1})
 \be
\int_{\BR} f(z)g(z)dz =\left\la
f,\int_{\BR}\frac{g(u)}{z-u}du\right\ra_{\iy}, \ee
involving the residue pairing\footnote{ The residue
pairing about $z=\iy$ between $f=\sum_{i\geq 0} a_i
z^i
\in
\HR^+$ and $g =\sum_{j \in \BZ}b_j z^{-j-1}\in \HR$ is
defined as: $$ \la
f,g\ra_{\iy}=\oint_{z=\iy}f(z)g(z)\frac{dz}{2 \pi i}=
\sum_{i \geq 0}a_i b_i . $$}. So, modulo terms
depending on $t$ and $t'$
 only, the left hand side of (3.1.14) equals
 $$
\oint_{z=\iy}dz~z^{-n}p_n(t;z)e^{\sum_1^{\iy}(t_i-t'_i)z^i}z^{n-m-1}z^{m+1}
\int_{\BR}\frac{p_m(t';u)}{z-u}e^{\sum_1^{\iy}t'_iu^i}\rho(u)du
$$
 \vspace{-.8cm}\bea &=&\int_{\BR}p_n(t;z)e^{\sum(t_i-t'_i)z^i}
p_m(t';z)e^{\sum t'_iz^i}\rho(z)dz, \mbox{\,using
(3.1.15),}\nonumber\\
&=&\int_{\BR}p_n(t;z)p_m(t';z)e^{\sum t_iz^i}\rho(z)dz
=0,\mbox{\,\,when $m\leq n-1$.}
 \eea

\noindent {\bf The KP-hierarchy:} Setting $n=m+1$,
shifting $t\mapsto t-y,t'\mapsto t+y$, evaluating the
residue and Taylor expanding in $y_k$ and using the
Schur polynomials $s_n$, leads to (see footnote 14 for
the definition of $p(\pl_t)f\circ g$.)
\begin{eqnarray*}
0&=&\frac{1}{2\pi i} \oint dz~ e^{-\sum_1^{\iy} 2y_i
z^i}\tau_n(t-y-[z^{-1}])
 \tau_n(t+y+[z^{-1}])\\
&=&\frac{1}{2\pi i}\oint dz ~\left(\sum_0^{\iy}
z^{i}s_i(-2y)\right) \left(\sum_0^{\iy}
z^{-j}s_j(\tilde \pl)\right)
 e^{\sum_1^{\iy} y_k\frac{\pl}{\pl t_k}}\tau_n\circ\tau_n\\
&=& e^{\sum_1^{\iy} y_k\frac{\pl}{\pl t_k}}
 \sum_0^{\iy} s_i(-2y)s_{i+1}(\tilde \pl)\tau_n\circ\tau_n\\
&=& \left(1+\sum_1^{\iy}y_j \frac{\pl}{\pl t_j}
 +O(y^2)\right)\left(\frac{\pl}{\pl t_1}+
 \sum_1^{\iy}s_{i+1}(\tilde \pl)(-2y_i +O(y^2))
 \right) \tau_n\circ\tau_n\\
 &=&\left( \frac{\pl}{\pl t_1}+\sum_{1}^{\iy} y_{k}\left(\frac{\pl}{\pl t_{k}}
 \frac{\pl}{\pl t_1} -2s_{k+1}(\tilde \pl)
 \right)\right)\tau_n
 \circ\tau_n+O(y^2),
\end{eqnarray*}
thus yielding (3.1.6), taking into acount the fact
that $(\pl / \pl t_1)\tau \circ \tau=0$ and the
coefficient of $y_k$ is trivial for $k=1,2$.

\noindent {\bf The Riemann-Hilbert problem:}
 Observe that, as a function of $z$, the
integral (3.1.10) has a jump accross the real axis $$
\frac{1}{2\pi i}\lim_{{z'\rg z}\atop{{\cal I}
z'<0}}\int_{\BR}\frac{p_n(t;u)}{z'-u}
\rho_t(u)du=p_n(t,z)\rho_t(z)+\frac{1}{2\pi
i}\lim_{{z'\rg z}\atop{{\cal I} z'>0}}
\int_{\BR}\frac{p_n(t;u)}{z'-u} \rho_t(u)du, $$ and
thus we have: (see \cite{FIK,Bleher,AvM22})

\begin{corollary}
 The matrix
 $$ Y_n(z)=\left(\begin{array}{lll}
\frac{\tau_n(t-[z^{-1}])}{\tau_n(t)}z^n&
&\frac{\tau_{n+1}(t+[z^{-1}])} {\tau_n(t)}z^{-n-1}\\
 & & \\
\frac{\tau_{n-1}(t-[z^{-1}])}{\tau_n(t)}z^{n-1}&
&\frac{\tau_{n}(t+[z^{-1}])}{\tau_n(t)}z^{-n}
\end{array}
\right) $$ satisfies the Riemann--Hilbert
problem\footnote{$\Bbb C_+$ and $\Bbb C_{-}$ denote
the Siegel upper and lower half plane.}:

\medbreak

1. $Y_n(z)$ holomorphic on the $\Bbb C_+$ and $\Bbb
C_{-}$ \,

\medbreak

2.
$Y_{n-}(z)=Y_{n+}(z)\displaystyle{\left(\begin{array}{lll}
1& &2\pi i\rho_t(z)\\
 & & \\
 0& &1
\end{array}\right)}$ (Jump condition)

\medbreak

3. $Y'_n(z)\displaystyle{ \left(\begin{array}{lll}
z^{-n}& &0\\
 & & \\
 0& &z^{n}
 \end{array}\right)}=1+O(z^{-1})$, when $z\rg\iy$.
 \end{corollary}


\subsection {Pfaff Lattice and symmetric/symplectic
 matrix integrals}

\subsubsection{Pfaff lattice, factorization of
 skew-symmetric matrices and skew-orthogonal polynomials}

Consider an inner-product, with a skew-symmetric
weight $\tilde\rho(y,z)$,
  \be \la
f,g\ra_t=\int\!\int_{\BR^2}f(y)g(z)e^{\sum
t_i(y^i+z^i)}\tilde\rho(y,z)dy\,dz,\mbox{\,\,with
$\tilde\rho(z,y)=-\tilde\rho(y,z)$}. \ee
  Since $\la
f,g\ra_t=-\la g,f\ra_t, $  the moment matrix,
depending on $t=(t_1,t_2,\ldots)$, $$
m_{n}(t)=(\mu_{ij}(t))_{0\leq i,j\leq n-1}=(\la
y^i,z^j\ra_t)_{0\leq i,j\leq n-1} $$ is
skew-symmetric. It is clear from formula (3.2.1) that
the semi-infinite matrix $m_{\iy}$ evolves in $t$
according to the {\em commuting vector fields}:
\be
\frac{\pl\mu_{ij}}{\pl
t_k}=\mu_{i+k,j}+\mu_{i,j+k},\mbox{\,\,i.e.,}~~
\frac{\pl m_{\iy}}{\pl
t_k}=\Lb^km_{\iy}+m_{\iy}\Lb^{\top k}. \ee

Since $m_{\iy}$ is skew-symmetric, $m_{\iy}$ does not
admit a Borel factorization in the standard sense, but
$m_{\iy}$ admits a unique factorization, with an
inserted semi-infinite, skew-symmetric matrix $J$,
 with
$J^2=-I$, of the form (1.1.12): (see \cite{AvM1})
 $$ m_{\iy}(t)=Q^{-1}(t)J\,Q^{\top
-1}(t), $$ where
 \be
 Q(t)=\left(
\begin{array}{c@{}c@{}cc}
\ddots &&&0 \\
&&0& \\
 & \boxed{\begin{array}{cc}
 Q_{2n,2n} & 0 \\ 0 & Q_{2n,2n} \end{array}} &&\\
 &*& \boxed{\begin{array}{cc}
  Q_{2n+2,2n+2} & 0 \\ 0 & Q_{2n+2,2n+2} \end{array}} & \\
 &&& \ddots
 \end{array}
 \right) ~\in K.
 \ee
$K$ is the group of lower-triangular invertible
matrices of the form above, with Lie algebra
${\mathfrak k}$. Consider the Lie algebra splitting,
given by
 $$ gl(\iy)={\frak
k}\oplus{\frak n}\left\{
\begin{array}{l}
{\frak k}=\{\mbox{lower-triangular matrices of the
form (3.2.3)}\}\\ {\frak n}=sp(\iy)=\{\mbox{$a$ such
that $Ja^{\top}J=a$}\},
\end{array}
\right. $$ \vspace{-1.3cm}\be \ee
 with unique decomposition\footnote{$a_{\pm}$ refers
to projection onto strictly upper (strictly lower)
triangular matrices, with all $2\times 2$ diagonal
blocks equal zero. $a_{0}$ refers to projection onto
the ``diagonal", consisting of $2\times 2$ blocks. }
 \bea a&=&(a)_{{\frak k}}
+(a)_{{\frak n}} \nonumber\\
&=&\left((a_--J(a_+)^{\top}J)+\frac{1}{2}
(a_0-J(a_0)^{\top}J)\right)  \nonumber\\ &&~~~~~~~~~+
\left((a_++J(a_+)^{\top}J)+\frac{1}{2}(a_0+J(a_0)^{\top}
J)\right). \eea


 Consider as a special skew-symmetric
weight (3.2.1): (see \cite{ASV3})
\be
 \tilde\rho(y,z)= 2 D^{\al}\delta (y-z)\tilde \rho (y)
  \tilde \rho (z)~, \mbox{with}~\al=\mp 1
,~~~\tilde \rho(y)=e^{-\tilde V(y)},  \ee
 together with the associated
 inner-product\footnote{$\vr(x)=
 \mbox{sign} ~x, $
having the property $\vr^{\prime}=2\delta(x)$.} of
type (3.2.1):
 \bea
\la f,g\ra_t&=&\int\!\int_{\BR^2}f(y)g(z)e^{\sum
t_i(y^i+z^i)}2 D^{\alpha}\delta (y-z)\tilde \rho (y)
  \tilde \rho (z)dy\,dz  \\ &&\nonumber\\ &=& \left\{\begin{array}{ll}
\displaystyle{\int\!\!\!\int_{\BR^2}}
f(y)g(z)e^{\sum_1^{\iy}t_i(y^i+z^i)}\vr(y-z)
 \tilde\rho(y)\tilde\rho(z)dy\,dz
, &\mbox{\,\,for\,\,}\alpha =-1\\
 & \\
\displaystyle{\int_{\BR}}\{f,g\}(y)e^{\sum_1^{\iy}2t_i
y^i}
 \tilde\rho(y)^2dy,&
\mbox{\,\,for\,\,}\alpha =+1
\end{array}
\right.\nonumber
 \eea
 in terms of the Wronskian $\{ f,g \}:=\frac{\pl f}{\pl
y} g-f \frac{\pl g}{\pl y}.$ The moments with regard
to these inner-products (with that precise definition
of time $t$!) satisfy the differential equations
$\pl\mu_{ij}/\pl t_k=\mu_{i+k,j}+\mu_{i,j+k}$, as in
(3.2.2).

 It is well known
that the determinant of an odd skew-symmetric matrix
equals 0, whereas the determinant of an even
skew-symmetric matrix is the square of a polynomial in
the entries, the {\em Pfaffian} \footnote{with a sign
specified below. So $\det(m_{2n-1}(t))=0$ and
 \bean
 (\det
m_{2n}(t))^{1/2}&=&pf(m_{2n}(t))\\
&=&\frac{1}{n!}(dx_{0}\wedge dx_1\wedge \ldots\wedge
dx_{2n-1})^{-1} \left(\sum_{0\leq i<j\leq
2n-1}\mu_{ij}(t)dx_i\wedge dx_j\right)^n. \eean}.
 Define now the {\em ``Pfaffian $\tau$-functions"},
 defined with regard to the inner-products (3.2.7)
 :
 \be
 \tau_{2n}(t):=\left\{
 \begin{array}{l}
 \displaystyle{
  pf\left( \int\!\!\!\int_{\BR^2} y^kz^{\ell}
 \varepsilon(y-z)e^{\sum_1^{\iy}t_i (y^i+z^i)}\tilde\rho(y)
\tilde \rho(z) dydz\right)_{0\leq k,\ell\leq
2n-1}},\al=-1
\\
 \displaystyle{  pf \left(
\int_{\BR}\{y^k,y^{\ell}\}
 e^{\sum_1^{\iy}2t_i y^i}\tilde\rho^2(y)dy
  \right)_{0\leq k,\ell\leq 2n-1},~~~~~~~~~~~~~~~~\al=+1
  }
 \end{array}
\right.
 \ee
  Setting
  $$ \left\{
\begin{array}{ll} \tilde \rho(z)=\rho(z)I_E(z) &
\mbox{for}~\alpha=-1\\
 \tilde
\rho(z)=\rho^{1/2}(z)I_E(z),~t\rightarrow t/2 &
\mbox{for}~\alpha=+1,
 \end{array}
 \right.
  $$
in the identities (3.2.8) leads to the identities
(3.2.9) between integrals and Pfaffians, spelled out
in Theorem 3.3 below. Remember $c=(c_1,...,c_{2r})$
stands for the boundary points of the disjoint union
$E\subset \BR$. ${\cal S}_{n}(E)$ and ${\cal
T}_{n}(E)$ denotes the set of symmetric/symplectic
ensemble with spectrum in $E$.

\begin{theorem} {\em (Adler-Horozov-van Moerbeke
 \cite{AHV}, Adler-van Moerbeke \cite{AvM4})}
The following integrals $I_n(t,c)$ are Pfaffians:
 \bea
 I_n& =&\int_{E^n}|\Dt_n(z)|^{\beta}\prod_{k=1}^n
\left(e^{\sum_1^{\iy}t_i z_k^i}\rho(z_k)dz_k\right)
\nonumber\\ &=&\left\{
\begin{array}{l}
\underline{\beta=1}\\
 ~~~~~=
 \displaystyle{\int_{{\cal S}_{n}(E)}} e^{\Tr (-
V(X)+\sum_1^{\iy} t_iX^i)}dX,~\underline{\mbox{($n$
even)}}
 \\ ~~~~~~\displaystyle{=n! pf\left( \int\!\!\!\int_{E^2} y^kz^{\ell}
 \varepsilon(y-z)e^{\sum_1^{\iy}t_i (y^i+z^i)}\rho(y)
 \rho(z) dydz\right)_{0\leq k,\ell\leq n-1}} \\
 ~~~~~=n!\tau_n(t,c)\\ \\
\underline{\beta=4}\\
 ~~~~\displaystyle{=\int_{{\cal
T}_{2n}(E)} }e^{\Tr (-  V(X)+\sum_1^{\iy} t_iX^i)}dX,
~\underline{\mbox{($n$ arbitrary)}}\\
 ~~~~=\displaystyle{ n! pf~\left(
\int_E\{y^k,y^{\ell}\}
 e^{\sum_1^{\iy}t_i y^i}\rho(y)dy
  \right)_{0\leq k,\ell\leq 2n-1}
  }\\~~~~=n!\tau_{2n}(t/2,c).
\end{array}
\right. \eea
 The $I_n$ and $\tau_n$'s satisfy
\begin{description}

  \item[(i)] The \underline{Virasoro constraints}\footnote{
  here the $a_i$'s and $b_i$'s are defined in the
  usual way, in terms of $\rho(z)$; namely,
  $-\frac{\rho'}{\rho}=\frac{\sum b_i z^i}{\sum a_i z^i}$
  } (2.1.7)
for $\beta =1,4$,
 \be
\left(-\sum_1^{2r} c_i^{k+1}f(c_i)\frac{\pl}{\pl c_i}+
\sum_{i\geq 0}\left( a_i~
{}^{\beta}\BJ_{k+i,n}^{(2)}-b_i ~
 {}^{\beta}\BJ_{k+i+1,n}^{(1)}\right) \right)I_n(t,c) =0.
\ee

 \item[(ii)] The \underline{Pfaff-KP hierarchy}: (see
 footnote 14 for notation)
\be
\left({\bf
s}_{k+4}(\tilde\pl)-\frac{1}{2}\frac{\pl^2}{\pl t_1\pl
t_{k+3}}\right)\tau_{n}\circ\tau_{n}={\bf s}_k(\tilde
\pl)~\tau_{n+2}\circ\tau_{n-2} \ee \hfill $n
~\mbox{even},~ k=0,1,2,...~.$

 \noindent of which the first equation reads ($n$
 even)
 $$ \left(\left(\frac{\pl}{\pl t_1}
\right)^4+3\left(\frac{\pl}{\pl
t_2}\right)^2-4\frac{\pl^2}{\pl t_1 \pl
t_3}\right)\log\tau_n+6\left(\frac{\pl^2}{\pl
t^2_1}\log\tau_n
\right)^2=12\frac{\tau_{n-2}\tau_{n+2}}{\tau_{n}^2}.
$$

 \item[(iii)] The \underline{Pfaff Lattice}: The
time-dependent matrix $L(t)$, zero above the first
superdiagonal, obtained by dressing up $\Lambda$,
 \be
L(t)=Q(t)\Lb Q(t)^{-1}=\left(\begin{array}{llllll}
 *&~1&  \\
  &~*&(h_2/h_0)^{1/2}& & ~~~~{ O}
  \\
  & & ~~~* &\!\!\! 1 &  \\
  & &   & \!\!\!* &(h_4/h_2)^{1/2} \\
  & &   &         & ~~~*& 
  \\
  &*  &   &         &     & \!\!\!\!\!\!\!\!\!\!\!\!\!\!
  \ddots\\
  & &   &         &
  \\
 \end{array} \right)
 \ee
satisfies the
 Hamiltonian commuting equations
\be \frac{\pl L}{\pl t_{i}}=[-(L^i)_{\frak
k},L].~~~~~~~\mbox{\bf (Pfaff lattice) }\ee

\item[(iv)] \underline{Skew-orthogonal polynomials}:
 The vector of time-dependent polynomials
$q(t;z):=(q_n(t;z))_{n\geq 0}=Q(t)\chi(z)$ in $z$
 satisfy the eigenvalue problem
\be
L(t)q(t,z)=zq(t,z) \ee and
 enjoy the following representations
 (with $h_{2n}=\frac{\tau_{2n+2}(t)}{\tau_{2n}(t)}$)
  \bea
 q_{2n}(t;z) &=&\frac{h_{2n}^{-1/2}}{\tau_{2n}(t)}pf\left(
\begin{array}{llll|l}
 & &
 & &1
 \\
 &m_{2n+1}(t)& & &z\\
 & & & &\vdots\\
 & & & &z^{2n}\\
 \hline
 -1&-z&\ldots&-z^{2n}&0
\end{array}
\right)\nonumber\\ & &
\nonumber\\&=&z^{2n}h_{2n}^{-1/2}
 \frac{\tau_{2n}(t-[z^{-1}])}{\tau_{2n}(t)}=z^{2n} h_{2n}^{-1/2}+... ~~,~~~~
  \nonumber\\
 &&  \nonumber\\
 q_{2n+1}(t;z)&=&\frac{h_{2n}^{-1/2}}{\tau_{2n}(t)}pf\left(
\begin{array}{lll|ll}
 & & &1&\mu_{0,2n+1}\\
 &m_{2n}(t)& &z&\mu_{1,2n+1}\\
 & & &\vdots&\vdots\\
 & & &z^{2n-1}&\mu_{2n-1,2n+1}\\
 \hline
 -1&\ldots&-z^{2n-1}&0&-z^{2n+1}\\
 \mu_{2n+1,0}&\ldots&\mu_{2n+1,2n-1}&z^{2n+1}&0
\end{array}
\right).\nonumber\\ & & \nonumber
\\&=&z^{2n}h_{2n}^{-1/2}\frac{1}{\tau_{2n}(t)}\left(z+\frac{\pl}{\pl
t_1}\right) \tau_{2n}(t-[z^{-1}])=z^{2n+1}
h_{2n}^{-1/2}+... . \nonumber\\\eea
 They are skew-orthogonal polynomials in $z$; i.e.,
  $$
  \la q_i(t;z),q_j(t;z)\ra_t=J_{ij} .
  $$

\end{description}

\end{theorem}

The hierarchy (3.2.11) already appears in the work of
Kac-van de Leur \cite{KvdL} in the context of, what
they call the DKP-hierarchy, and interesting further
work has been done by van de Leur \cite{vdL}.


\subsubsection
 { Sketch of Proof}

\noindent{\bf Skew-orthogonal polynomials and the
Pfaff Lattice}: The equalities (3.2.9) between the
Pfaffians and the matrix integrals are based on two
identities \cite{ Mehta}, the first one due to de
Bruyn,
 \bean &&
\hspace{-1.5cm}\frac{1}{n!}\int_{\BR^n}\prod^n_1dy_i~~
\det\Bigl(F_i(y_1)~~G_i(y_1)~~\dots~~F_i(y_n)~~
G_i(y_n)\Bigr)_{0\leq i\leq 2n-1}\\
\hspace{1.5cm} &=& {\det}^{1/2} \left(\int_{\BR}
 (G_i(y)F_j(y)-F_i(y) G_j(y))dy\right)_{0\leq
i,j\leq 2n-1}\eean
 and (Mehta \cite{Mehta2})
  $$
\Delta_n^4(x)=\det\left(x_1^i~~~(x_1^i)'~~~x_2^i~~~(x_2^i)'~\dots~
~x_n^i~~~(x_n^i)'\right)_{0\leq i\leq 2n-1}.$$

On the one hand, (see Mehta \cite{Mehta}), setting in
the calculation below $\rho_{t}(z)=  \rho(z)e^{\sum
t_i z^i}I_E(z)$ and $$
F_i(x):=\int^x_{-\iy}y^i\rho_t(y)dy \quad\mbox{and}
\quad G_i(x):=F_i'(x)=x^i\rho_t(x), $$ one computes:
 ($\rho_t(z):=\rho(z) e^{\sum t_i z^i}$)
 \bean
\lefteqn{\frac{1}{(2n)!}\int_{\BR^{2n}}|\Delta_{2n}(z)|
\prod_{i=1}^{2n} \rho_t(z_i)dz_i}\\
&=&\int_{-\iy<z_1<z_2<\cdots<z_{2n}<\iy}\det\left(z_{j+1}^i
\rho_t(z_{j+1})\right)_{0\leq i,j\leq
2n-1}\prod_{i=1}^{2n} dz_i,\\
&=&\int_{-\iy<z_2<z_4<\cdots<z_{2n}<\iy}
\prod^n_{k=1}\rho_t(z_{2k})dz_{2k} \\ & &
\det\left(\int_{-\iy}^{z_2}z_1^i\rho_t(z_1)dz_1~,~z_2^i~,\dots,~
\int_{z_{2n-2}}^{z_{2n}}z_{2n-1}^i\rho_t(z_{2n-1})dz_{2n-1}~,
~z_{2n}^i \right)_{0\leq i\leq 2n-1}\\
&=&\int_{-\iy<z_2<z_4<\cdots<z_{2n}<\iy}
\prod^n_{k=1}\rho_t(z_{2k})dz_{2k}\\ & &
\det\left(F_i(z_2)~,~z_2^i~,~F_i(z_4)-F_i(z_2)~,~
z_4^i~,~\dots~,~F_i(z_{2n})-F_i(z_{2n-2})~,~z^i_{2n}\right)_{0\leq
i\leq 2n-1}\\
 &=&\int_{-\iy<z_2<z_4<\cdots<z_{2n}<\iy}\prod^n_1dz_i~~
\det\Bigl(F_i(z_2)~,~G_i(z_2)~,~\dots,~F_i(z_{2n})~,~
G_i(z_{2n})\Bigr)_{0\leq i\leq 2n-1},\\
 &=&\frac{1}{n!}\int_{\BR^n}\prod^n_1dy_i~~
\det\Bigl(F_i(y_1)~,~G_i(y_1)~,~\dots,~F_i(y_n)~,~
G_i(y_n)\Bigr)_{0\leq i\leq 2n-1},\\
&=&{\det}^{1/2}\left(\int_{\BR}
 (G_i(y)F_j(y)-F_i(y) G_j(y))dy\right)_{0\leq
i,j\leq 2n-1}\\
 &=& pf\left( \int\!\!\!\int_{E^2} y^kz^{\ell}
 \varepsilon(y-z)e^{\sum_1^{\iy}t_i (y^i+z^i)}\rho(y)
 \rho(z) dydz\right)_{0\leq k,\ell\leq 2n-1}=
 \tau_{2n}(t),
\eean
 establishing the first equation of (3.2.9), taking into account the results in section 2.1.

 On the other hand, upon setting,
  $$
 F_j(x)=x^j \rho(x)e^{\sum t_i x^i} ~~\mbox{and}~~
 G_j(x):=F'_j(x)=\left(x^j \rho(x)e^{\sum t_i x^i}\right)' ,
$$ one computes
 \bean \lefteqn{\frac{1}{n!}\int_{\BR^n}\prod_{1\leq
i,j\leq n}(x_i-x_j)^4
  \prod_{k=1}^{n}\left( \rho^2(x_k)e^{2\sum_{i=1}^{\iy}  t_i x_k^i}
  dx_k   \right)}\\
 &=&\frac{1}{n!}\int_{\BR^n}
\prod^n_{k=1}\left(\rho^2(x_k)e^{2\sum  t_i x_k^i}dx_k
 \right)\\ & &
\hspace{3cm}\det\left(x_1^i~~~(x_1^i)'~~~x_2^i~~~(x_2^i)'~\dots~
~x_n^i~~~(x_n^i)'\right)_{0\leq i\leq 2n-1}\\
&=&\frac{1}{n!}\int_{\BR^n}\prod^n_1dy_i~~
\det\Bigl(F_i(y_1)~~G_i(y_1)~~\dots~~F_i(y_n)~~
G_i(y_n)\Bigr)_{0\leq i\leq 2n-1},\\
&=& {\det}^{1/2} \left(\int_{\BR}
 (G_i(y)F_j(y)-F_i(y) G_j(y))dy\right)_{0\leq
i,j\leq 2n-1}\\ &=&\displaystyle{  pf~\left(
\int_E\{y^k,y^{\ell}\}
 e^{\sum_1^{\iy}2t_i y^i}\rho^2(y)dy
  \right)_{0\leq k,\ell\leq 2n-1}= \tau_{2n}(t)
  },
 \eean establishing the second
equation (3.2.9).

 The skew-orthogonality of the polynomials $q_k(t;z)$
  follows immediately from the skew-Borel decomposition of
  $m_{\iy}$: $$ \la
q_k(t,y),q_{\ell}(t,z)\ra_{k,\ell\geq 0}=Q(\la
y^i,z^j\ra)_{i,j\geq 0}Q^{\top}=Q\,m_{\iy}Q^{\top}=J.
$$
 \noindent with the $q_n$'s admitting the
representation (3.2.15) in terms of the moments.

 Using $L=Q\Lb Q^{-1}$,
$m_{\iy}=Q^{-1}JQ^{\top -1}$ and $J^2=-I$, one
computes from the differential equations (3.2.2)
 \bean
0&=&Q\left(\Lb^km_{\iy}+m_{\iy}\Lb^{\top k}-\frac{\pl
m_{\iy}}{\pl t_k} \right)Q^{\top}\\
&=&(Q\Lb^kQ^{-1})J-(JQ^{\top -1}\Lb^{\top k}Q^{\top
}J)J+\frac{\pl Q}{\pl t_k} Q^{-1}J-(JQ^{-1 \top}
 \frac{\pl Q^{\top }}{\pl t_k}J)J\\ &=&\left(L^k+\frac{\pl Q}{\pl
t_k}Q^{-1}\right)-J\left(L^k+\frac{\pl Q}{\pl t_k}
Q^{-1}\right)^{\top}J. \eean Then computing the +, -
and the diagonal part (in the sense of (3.2.4) and
(3.2.5)) of the expression leads to commuting
Hamiltonian differential equations for $Q$, and thus
for $L$ and $q(t;z)$, confirming (3.2.13):
\be
\frac{\pl Q}{\pl t_{i}}=-(L^i)_{\frak
k}Q,\quad\frac{\pl L}{\pl t_{i}}=[(L^i)_{\frak
n},L],\quad\frac{\pl q}{\pl t_{i}}=-(L^i)_{\frak k}q
.~\mbox{(Pfaff lattice)}\ee

\noindent {\bf The bilinear identities:} For all
$n,m\geq 0$, the $\tau_{2n}$'s satisfy the following
bilinear identity
 \bea &
&\oint_{z=\iy}\tau_{2n}(t-[z^{-1}])\tau_{2m+2}(t'+[z^{-1}])
e^{\sum(t_i-t'_i)z^i} z^{2n-2m-2}\frac{dz}{2\pi
i}\nonumber\\
&&\quad+\oint_{z=0}\tau_{2n+2}(t+[z])\tau_{2m}(t'-[z])
e^{\sum(t'_i-t_i)z^{-i}}z^{2n-2m}\frac{dz}{2\pi
i}=0.\nonumber\\ \eea
 The differential equation (3.2.2) on the moment matrix $m_{\iy}$ admits the following solution,
which upon using the Borel decomposition
$m_{\iy}=Q^{-1}J\,Q^{\top -1}$, leads to:
\be
m_{\iy}(0)=e^{-\sum_1^{\iy}t_k\Lb^k}m_{\iy}(t)e^{-\sum_1^{\iy}
t_k\Lb^{\top k}}
=\left(Q(t)e^{\sum_1^{\iy}t_k\Lb^k}\right)^{-1}J\left(Q(t)e^{
\sum_1^{\iy}t_k\Lb^k}\right)^{\top -1}, \ee and so the
right hand side of (3.2.11) is independent of $t$;
say, equal to the same expression with $t$ replaced by
$t'$. Upon rearrangement, one finds
$$ \left(Q(t)e^{\sum
t_k\Lb^k}\right)\left(JQ(t')e^{\sum
t'_k\Lb^k}\right)^{-1}= \left(JQ(t)e^{\sum
t_k\Lb^k}\right)^{\top -1}\left(Q(t')e^{\sum
t'_k\Lb^k} \right)^{\top} $$ and
therefore\footnote{using $\Lb\chi(z)=z\chi(z)$,
$\Lb^{\top}\chi(z)=z^{-1}\chi(z)$ and the following
matrix identities (see \cite{DJKM}) $$
U_1V_1=\oint_{z=\iy}U_1\chi(z)\otimes
V_1^{\top}\chi(z^{-1})\frac{dz}{2\pi iz},\quad
U_2V_2=\oint_{z=0}U_2\chi(z)\otimes
V_2^{\top}\chi(z^{-1})\frac{dz}{2\pi iz}. $$ } \bea &
&\hspace{-.5cm}\oint_{z=\iy}\left(Q(t)\chi(z)\otimes(JQ(t'))^{\top
-1}\chi(z^{-1})\right)e^{\sum_1^{\iy}(t_k-t'_k)z^k}\frac{dz}{2\pi
iz}\nonumber\\ &=&\oint_{z=0}\left((JQ(t))^{\top
-1}\chi(z)\otimes Q(t')\chi(z^{-1})\right)
e^{\sum_1^{\iy}(t'_k-t_k)z^{-k}}\frac{dz}{2\pi iz}.
\eea Setting $t-t'=[z_1^{-1}]+[z_2^{-1}]$ into the
exponential leads to
 \bean
e^{\sum_1^{\iy}(t_k-t'_k)z^k}&=&\left(1-\frac{z}{z_1}\right)^{-1}\left(1-\frac{z}{
 z_2}\right)^{-1}\\
 e^{\sum_1^{\iy}(t'_k-t_k)z^{-k}}
  &=&
\left(1-\frac{1}{zz_1}\right)\left(1-\frac{1}{
zz_2}\right)
 \eean
  and somewhat enlarging the
integration circle about $z=\iy$ to include the points
$z_1$ and $z_2$, the integrand on the left hand side
has poles at $z=z_1$ and $z_2$, whereas the integrand
on the right hand side is holomorphic. Combining the
identity obtained  and the one, with $z_2\nearrow\iy$,
one finds a functional relation involving a function
$\varphi(t;z)=1+O(z^{-1})$: $$
\frac{\vp(t-[z_2^{-1}];z_1)}{\vp(t;z_1)}=
\frac{\vp(t-[z_1^{-1}];z_2)}{\vp(t;z_2)}, \quad
t\in\BC^{\iy},z\in\BC. $$ Such an identity leads, by a
standard argument (see e.g. the appendix in
\cite{AvM5}) to the existence of a function $\tau(t)$
such that $$
\vp(t;z)=\frac{\tau(t-[z^{-1}])}{\tau(t)}. $$ This
fact combined with the bilinear identity (3.2.19)
leads to the bilinear identity (3.2.17).

\noindent {\bf The Pfaff-KP-hierarchy:} Shifting
$t\mapsto t-y,t'\mapsto t+y$ in (3.2.17), evaluating
the residue and Taylor expanding in $y_k$ leads to:

\begin{multline*}
\frac{1}{2\pi i}\oint_{z=\infty}
e^{-\sum_1^\infty2y_iz^i}
  \tau_{2n}(t-y-[z^{-1}])\tau_{2m+2}(t+y+[z^{-1}])
z^{2n-2m-2}dz
\\
{}+ \frac{1}{2\pi i}\oint_{z=0}
e^{\sum_1^\infty2y_iz^{-i}}
\tau_{2n+2}(t-y+[z])\tau_{2m}(t+y-[z]) z^{2n-2m}dz
\\
\begin{split}
&=\frac{1}{2\pi i} \oint_{z=\infty} \sum_{j=0}^\infty
z^j {\bf s}_j(-2y) e^{\sum-y_i\frac{\pl}{\pl
t_i}}\sum_{k=0}^\infty z^{-k} {\bf s}_k(-\tilde \pl)
\tau_{2n}\circ\tau_{2m+2} z^{2n-2m-2}dz
\\
&\phantom{={}}+\frac{1}{2\pi i}\oint_{z=0}
\sum_{j=0}^\infty z^{-j} {\bf s}_j(2y)
e^{\sum-y_i\frac{\pl}{\pl t_i}}\sum_{k=0}^\infty z^k
{\bf s}_k(\tilde \pl) \tau_{2n+2}\circ\tau_{2m}
z^{2n-2m}dz
\\
\end{split}
\end{multline*}
\begin{multline*}
\begin{split}
&= \sum_{j-k=-2n+2m+1} {\bf s}_j(-2y) e^{\sum-y_i
\frac{\pl}{\pl t_i} }{\bf s}_k(-\tilde \pl)
\tau_{2n}\circ\tau_{2m+2}
\\
&\hphantom{=2\pi i\biggl(}+ \sum_{k-j=-2n+2m-1} {\bf
s}_j(2y) e^{\sum-y_i\frac{\pl}{\pl t_i}}{\bf
s}_k(\tilde \pl) \tau_{2n+2}\circ\tau_{2m}\\
 &=   ...+y_{k}
 \left( \left(\frac{1}{2}\frac{\pl}{\pl t_1}
  \frac{\pl}{\pl t_{k}}-
 {\bf s}_{k+1}(\tilde \pl
 )
\right)\tau_{2n}\circ\tau_{2n}+{\bf s}_{k-3}(\tilde
\pl)\tau_{2n+2}\circ\tau_{2n-2} \right)+..., \\
\end{split}
\end{multline*}
establishing the Pfaff-KP hierarchy (3.2.11),
different from the usual KP hierarchy, because of the
presence of a right hand side.

\remark $L$ admits the following representation in
terms of $\tau$, much in the style of (3.1.7),

$$ L  =h^{-1/2} \left(\begin{tabular}{lllll}
  $ \hat L_{00}$ &  $\hat L_{01}$ & $0$ & $0$ &   \\
  $\hat L_{10}$ & $\hat L_{11}$   &$\hat L_{12}$&$0$&  \\
 $*$ & $\hat L_{21}$ & $\hat L_{22}$ & $\hat L_{23}$ &   \\
 $*$&$*$ & $\hat L_{32}$ & $\hat L_{33}$ & $\dots $  \\
  & & &$\vdots$&     \end{tabular}
\right) h^{1/2}, $$
 with the $2\times 2$ entries $\hat L_{ij}$ and
$h$, being a zero matrix, except for $2\times 2$
matrices along the diagonal:
 $$
h=\mbox{diag}(h_0I_2,h_2I_2,h_4I_2,\dots),
~h_{2n}=\tau_{2n+2}/\tau_{2n}. $$ For example ($
{}^.=\frac{\pl}{\pl t_1}$)

$$ \hat L_{nn}:=\left(\begin{tabular}{lll}
 $-(\log \tau_{2n})^.$ & & ~~~~$1$ \\
 $ -\frac{s_2(\tilde
\pl)\tau_{2n}}{\tau_{2n}} -\frac{s_2(-\tilde
\pl)\tau_{2n+2}}{\tau_{2n+2}}$ & &
 $(\log \tau_{2n+2})^.$  \end{tabular}
\right)
 $$
 $$~~~~~~\hat L_{n,n+1}:=\left(\begin{tabular}{ll}
$0$& $0$ \\
 $1$ & $0$   \end{tabular} \right)
 ~~~~~~\hat L_{n+1,n}:=\left(\begin{tabular}{lll}
$*$&$(\log \tau_{2n+2})^{..}$&
\\
$ *$ &~~~~~$*$\end{tabular} \right). $$

\vspace{0.4cm}


\subsection{2d-Toda lattice and coupled Hermitian
matrix integrals}

\subsubsection{2d-Toda lattice, factorization of
 moment matrices and bi-orthogonal polynomials}

Consider the inner-product,
  \be \la
f,g\ra_{t,s}=\int\!\!\int_{E \subset
\BR^2}f(y)g(z)e^{\sum_1^{\iy} ( t_i y^i-s_i z^i)+cyz }
dy\,dz,
 \ee
on a subset $E=E_1\times
E_2:=\cup^r_{i=1}[c_{2i-1},c_{2i}]\times
\cup^s_{i=1}[\tilde c_{2i-1},\tilde c_{2i}]\subset F_1
\times  F_2\subset\BR^2$. Define the customary moment
matrix, depending on
$t=(t_1,t_2,\ldots),~s=(s_1,s_2,\ldots)$, $$
m_{n}(t,s)=(\mu_{ij}(t,s))_{0\leq i,j\leq n-1}=(\la
y^i,z^j\ra_{t,s})_{0\leq i,j\leq n-1} $$ and its
factorization in lower- times upper-triangular
matrices \be m_{\iy}(t,s)= S_1^{-1}(t,s) S_2(t,s).\ee
Then $m_{\iy}$ evolves in $t,s$ according to the
equations
\be
\frac{\pl\mu_{ij}}{\pl
t_k}=\mu_{i+k,j},~~~\frac{\pl\mu_{ij}}{\pl
s_k}=-\mu_{i,j+k},\mbox{\,\,i.e.,}~~ \frac{\pl
m_{\iy}}{\pl t_k}=\Lb^km_{\iy},~~
 \frac{\pl
m_{\iy}}{\pl s_k}=-m_{\iy}\Lb^{\top k}. \ee
$dM$ in the integral (3.3.4) denotes properly
normalized Haar measure on ${\cal H}_n$.

\begin{theorem}
{\em (Adler-van Moerbeke \cite{AvM21,AvM2})}
 The integrals $I_n(t,s;c,\tilde c)$, with $I_0=1$,
  \bea
  \tau_n =\det
m_n \hspace{-.2cm}&=&\frac{1}{n!}I_n=\frac{1}{n!}\int
\hspace{-.3cm}\int_{E^n} \Delta_n(x)\Delta_n(y)
\prod^n_{k=1}
e^{\sum^{\iy}_{1}(t_ix_k^i-s_iy^i_k)+cx_ky_k}
 dx_kdy_k \nonumber\\
 \hspace{-.2cm} &=& \int \hspace{-.3cm}\int_{{\cal H}^2_{n}(E)}
e^{c\Tr(M_1M_2)}
e^{\Tr\sum_1^{\iy}(t_iM^i_1-s_iM_2^i)}dM_1dM_2, \eea
satisfy:
\begin{description}

  \item[(i)] \underline{Virasoro constraints}\footnote{For the Hirota symbol, see footnote
14.  The $J_k^{(i)}$'s are as in remark 1 at the end
of Theorem 2.1, for $\beta=1$ and $\tilde
J_k^{(i)}=J_k^{(i)}|_{t\rightarrow -s}$, with
 $$ 
 J_k^{(1)}=\frac{\pl}{\pl
t_k}+(-k)t_{-k}
 ,~~
  J^{(2)}_{k}=\sum_{i+j=k}\frac{\pl^2}{\pl
 t_{i}\pl t_{j}}+2\sum_{-i+j=k}it_{i}\frac{\pl}{\pl
 t_{j}}+\sum_{-i-j=k}it_{i}jt_{j}.
$$} (2.2.8) for $k \geq -1$,
 \bea
\left(-\sum_{i=1}^rc_i^{k+1}\frac{\pl}{\pl
c_i}+J_{k,n}^{(2)}\right)\tau_n^E+c\,{\bf s}
_{k+n}(\tilde\pl_t){\bf s}_n
(-\tilde\pl_s)\tau^E_1\circ\tau^E_{n-1}&=&0
\nonumber\\
  \left(-\sum_{i=1}^s\tilde
c_i^{k+1}\frac{\pl}{\pl \tilde c_i}+\tilde
J_{k,n}^{(2)}\right)\tau_n^E+c\,{\bf
s}_n(\tilde\pl_t){\bf s}_{k+n}
(-\tilde\pl_s)\tau^E_1\circ\tau^E_{n-1}&=&0,\nonumber\\
\eea with
 \bean
J_{k,n}^{(2)}&=&\frac{1}{2}(J_k^{(2)}
+(2n+k+1)J_k^{(1)}+n(n+1)J_k^{(0)}), \nonumber\\
\tilde J_{k,n}^{(2)}&=&\frac{1}{2}(\tilde J_k^{(2)}
+(2n+k+1)\tilde J_k^{(1)}+n(n+1)J_k^{(0)}). \eean

\item[(ii)] A \underline{Wronskian identity}
\footnote{in terms of the Wronskian $\{ f,g
\}_t=\frac{\pl f}{\pl t} g-f \frac{\pl g}{\pl t}.$}:
\be \left\{\frac{\pl^2 \log \tau_{n}}{\pl t_1 \pl
s_2},
 \frac{\pl^2 \log \tau_{n}}{\pl t_1 \pl s_1}   \right\}_{t_1}
  +\left\{\frac{\pl^2 \log \tau_{n}}{\pl s_1 \pl t_2},
  \frac{\pl^2 \log \tau_{n}}{\pl t_1 \pl s_1}   \right\}_{s_1}
  =0.
\ee

\item[(iii)] The \underline{2d-Toda lattice}: Given
the factorization (3.3.2), the matrices
$L_1:=S_1\Lambda
S_1^{-1}~~~\mbox{and}~~~L_2:=S_2\Lambda^{\top}
S_2^{-1}$, with $h_n=\frac{\tau_{n+1}}{\tau_n}$, have
the form (i.e., to be read as follows: the
$(k-\ell)$th subdiagonal is given by the diagonal
matrix in front of $\Lb^{k-\ell}$)
 \bea
 L_1^k&=&\sum_{\ell=0}^{\iy}\mbox{diag}~
 \left(\frac{{\bf s}_{\ell}(\tilde\pl_t)
\tau_{n+k-\ell+1}\circ\tau_n} {\tau_{n+k-\ell+1}
\tau_n}\right)_{n \in \BZ}\Lb^{k-\ell} \nonumber\\
hL_2^{\top k}h^{-1}&=& \sum_{\ell=0}^{\iy}\mbox{diag}~
 \left(\frac{{\bf s}_{\ell}(-\tilde\pl_s)
\tau_{n+k-\ell+1}\circ\tau_n} {\tau_{n+k-\ell+1}
\tau_n}\right)_{n\in\BZ}\Lb^{k-\ell},~~~~ 
 \eea
and satisfy the {\bf 2d-Toda Lattice}\footnote{$P_+$
and $P_-$ denote the upper (including diagonal) and
strictly lower triangular parts of the matrix $P$,
respectively. }
 \begin{equation}
\frac{\pl L_{i}}{\pl
t_n}=[(L^n_1)_+~,L_{i}]~~~\mbox{\,\,and\,\,} ~~~
\frac{\pl L_{i}}{\pl s_n}=[(L^n_2)_{-}~,L_{i}] ~,\quad
i=1,2,\end{equation}

\item[(iv)]  \underline{Bi-orthogonal polynomials}:
The expressions \bea
p_n^{(1)}(t,s;y)&:=&\left(S_1(t,s)\chi)y)\right)_n=
y^n
\frac{\tau_n(t-[y^{-1}],s)}{\tau_n(t,s)}\nonumber\\
 p_n^{(2)}(t,s;z)&:=&\left(hS_2^{\top
 -1}(t,s)\chi(z)\right)_n=
  z^n \frac{\tau_n(t,s+[z^{-1}])}{\tau_n(t,s)}
~~~~~~ \eea form a system of monic bi-orthogonal
polynomials in $z$: \be\la
p_n^{(1)}(t,s;y),p_m^{(2)}(t,s;z)\ra_{t,s}
=\dt_{n,m}h_n~~\mbox{with}~~h_n=\frac{\tau_{n+1}}{\tau_n},
 \ee
which also are eigenvectors of $L_1$ and $L_2$: \be
\hspace{-1cm}
zp_n^{(1)}(t,s;z)=L_1(t,s)p_n^{(1)}(t,s;z)
~~~\mbox{and}~~~
zp_n^{(2)}(t,s;z)=L_2^{\top}(t,s)p_n^{(2)}(t,s;z). \ee
\end{description}

\end{theorem}

\remark Notice that every statement can be dualized,
upon using the duality  $t\longleftrightarrow -s$,
$L_1\longleftrightarrow
 hL_2^{\top}h^{-1}$.

\subsubsection{Sketch of proof}

Identity (3.3.4) follows from the fact that the
product of the two Vandermonde appearing in the
integral (3.3.4) can be expressed as sum of
determinants:
 \be \Delta_n(u)\Delta_n(v)=\sum_{\sigma\in
S_n}\det\Bigl(u^{\ell-1}_{\sg(k)}
v^{k-1}_{\sg(k)}\Bigr)_{1\leq\ell,k\leq n}, \ee
together with the Harish-Chandra, Itzykson and Zuber
formula \cite{Harish,IZ}
 \begin{equation}
\int_{{ U}(n)}dU\,e^{c\Tr xUy\bar
U^{\top}}=\frac{(2\pi)^{\frac{n(n-1)}{2}}}{n!}
\frac{\det(e^{cx_iy_j})_{1\leq i,j\leq
n}}{\Delta_n(x)\Delta_n(y)}.
\end{equation}
 Moreover the $\tau_n$'s satisfy the following bilinear
identities, for all integer $m,n\geq 0$ and $t,s \in
\BC^{\iy}$: $$
\oint_{z=\iy}\tau_n(t-[z^{-1}],s)\tau_{m+1}(t'+[z^{-1}],s')
e^{\sum_1^{\iy}(t_i-t'_i)z^i} z^{n-m-1}dz $$ \be
=\oint_{z=0}\tau_{n+1}(t,s-[z])\tau_m(t',s'+[z])
e^{\sum_1^{\iy}(s_i-s'_i)z^{-i}}z^{n-m-1}dz. \ee

Again, the bi-orthogonal nature (3.3.10) of the
polynomials (3.3.9) is tantamount to the Borel
decomposition, written in the form
$S_1m_{\iy}(hS_2^{\top -1})^{\top}=h$. These
polynomials satisfy the eigenvalue problem (3.3.11)
and evolve in $t,s$ according to the differential
equations
 \bea
  \frac{\pl p^{(1)}}{\pl t_n}=-(L_1^n)_-p^{(1)} ~~~~~~~~~~~~~    &&
  \frac{\pl p^{(1)}}{\pl s_n}=-(L_2^n)_-p^{(1)}   \nonumber  \\
\frac{\pl p^{(2)}}{ \pl t_n}=-\left((h^{-1}L_1
h)^{\top n}\right)_-p^{(2)} &&
  \frac{\pl p^{(2)}}{\pl s_n}=
  \left((h^{-1}L_2 h)^{\top n}\right)_-p^{(2)}.    \eea

 From the representation (3.3.7)
and from the bilinear identity (3.3.14), it follows
that
\begin{equation}
\frac{p_{k-1}(\tilde\pl_t)\tau_{n+2}\circ \tau_n}
{\tau_{n+1}^2}=-\frac{\pl^2}{\pl s_1 \pl t_k}\log
\tau_{n+1} ,
 \end{equation}
and so, for $k=1$,
\begin{equation}
\frac{\tau_n \tau_{n+2} }
{\tau_{n+1}^2}=-\frac{\pl^2}{\pl s_1 \pl t_1}\log
\tau_{n+1} .
 \end{equation}
Thus, using (3.3.7), (3.3.16) and (3.3.17), we have
 \begin{eqnarray}
\left(L_1^k\right)_{n,n+1}&=&\frac{p_{k-1}
 (\tilde\pl_t )\tau_{n+2}
 \circ \tau_n}{\tau_{n+2} \tau_n}
 ~=~~\frac{\frac{\pl^2 \log \tau_{n+1}}{\pl s_1 \pl t_{k}}}
 {\frac{\pl^2 \log \tau_{n+1}}{\pl s_1 \pl t_{1}}}
\nonumber\\ \left(hL_2^{\top k}h^{-1}\right)_{n,n+1} &
=&\frac{p_{k-1}(-\tilde\pl_s )\tau_{n+2}
 \circ \tau_n}{\tau_{n+2} \tau_n}
 =~\frac{\frac{\pl^2 \log \tau_{n+1}}{\pl t_1 \pl s_{k}}}
 {\frac{\pl^2 \log \tau_{n+1}}{\pl s_1 \pl t_{1}}}.
\end{eqnarray}
Combining (3.3.18) with (3.3.17) for $k=2$ yields \bea
\left(L_1^2\right)_{n,n+1}
= \frac{\frac{\pl^2} {\pl s_1\pl
t_2}\log\tau_{n+1}}{\frac{\pl^2}{\pl s_1\pl t_1}
 \log\tau_{n+1}}&=&
 \frac{\pl}{\pl t_1}\log
  \left(-\frac{\tau_{n+2}}{\tau_n}  \right)\\ &=&
 \frac{\pl}{\pl t_1}\log\left(\left(\frac{\tau_{n+1}}{\tau_n }
  \right)^2
\frac{\pl^2}{\pl s_1\pl t_1}\log\tau_{n+1}
 \right). \nonumber\eea
 Then, subtracting $\pl/\pl s_1$
 of (3.3.19) from $\pl/\pl t_1$ of the dual of the
 same equation (see remark at the end of Theorem 3.4) leads to (3.3.6).
\qed

\subsection{The Toeplitz Lattice and Unitary matrix
integrals}

\subsubsection{Toeplitz lattice, factorization of
 moment matrices and bi-orthogonal polynomials}

Consider the inner-product
\be
\la f(z),g(z)\ra_{t,s}:=\oint_{S^1} \frac{dz}{2\pi i
z}f(z)g(z^{-1})
 e^{\sum_1^{\iy}(t_iz^i-s_iz^{-i})},~~~~t,s\in
 \BC^{\iy},
 \ee
where the integral is taken over the unit circle
$S^1\subset \BC$ around the origin. It has the
property
\be
\la z^k f,g\ra _{t,s}=     \la f, z^{-k}g \ra_{t,s}.
\ee The $t,s$ dependent semi-infinite moment matrix
$m_{\iy}(t,s)$, where\footnote{ A matrix is Toeplitz,
when its $(i,j)$th entry depends on $i-j$.}
 \bea m_n(t,s):=\left(\la
z^k,z^{\ell}\ra_{t,s}
 \right)_{0\leq k,\ell \leq n-1}\hspace{-.3cm}&=& \left(\oint_{S^1} \frac{\rho(z)dz}{2\pi i z}
z^{k-\ell}
e^{\sum_1^{\iy}(t_iz^i-s_iz^{-i})}\right)_{0\leq
k,\ell \leq n-1}\nonumber\\ \hspace{-.3cm}&=&
\mbox{Toeplitz matrix}
 \eea
satisfies the same differential equations, as in
(3.3.3):
 \be
 \frac{\pl m_{\iy}}{\pl t_n}=\Lb^n
m_{\iy}~~\mbox{and}~~ \frac{\pl m_{\iy}}{\pl s_n}=-
m_{\iy} \Lb^{\top n}.~~~~~\mbox{\sl(2-Toda
Lattice)}
 \ee
 As before, define
  $$
 \tau_n(t,s):= \det m_n(t,s).
 $$
 Also, consider the factorization $m_{\iy}(t,s)= S_1^{-1}(t,s)
 S_2(t,s)$, as in (3.3.2), from which one defines $L_1:=S_1\Lambda
S_1^{-1}~~~\mbox{and}~~~L_2:=S_2\Lambda^{\top}
S_2^{-1}$ and the bi-orthogonal polynomials
$p_i^{(k)}(t,s;z)$ for $k=1,2$. Since $m_{\iy}$
satisfies the same equations (3.3.3), the matrices
$L_1$ and $L_2$ satisfy the 2-Toda lattice equations;
the Toeplitz nature of $m_{\iy}$ implies a peculiar
``rank 2"-structure, with
$\frac{h_i}{h_{i-1}}=1-x_iy_i$ and $x_0=y_0=1$ :
$$h^{-1}L_1 h= \left(\begin{tabular}{lllll} $-x_1y_0$
& $1-x_1y_1$ & ~~ $0$      & ~~ $0$ &   \\ $-x_2y_0$ &
$-x_2y_1$  & $1-x_2y_2$& ~~ $0$   & \\ $-x_3y_0$ &
$-x_3y_1$ & $ -x_3y_2$&  $1-x_3y_3$ & \\ $ -x_4y_0$ &
$ -x_4y_1$ & $-x_4y_2$  & $ -x_4y_3$   &
\\
 & &  &    &  $\ddots$\\
\end{tabular}
\right) $$ and \be L_2= \left(\begin{tabular}{lllll}
$-x_0y_1$  &  $-x_0y_2$ & $-x_0y_3$     & $-x_0y_4$ &
\\ $1 -x_1y_1$ &  $-x_1y_2$  & $-x_1y_3$& $-x_1y_4$
& \\ ~~$0$       &  $1 -x_2y_2$ & $ -x_2y_3$&
$-x_2y_4$ & \\ ~~$0$       &  ~~$0$      & $ 1
-x_3y_3$  & $ -x_3y_4$   &  \\
 & &  &    &  $\ddots$\\
\end{tabular}
\right). \ee

\noindent Some of the ideas in the next theorem are
inspired by the work of Hisakado \cite{H}.

\begin{theorem} {\em (Adler-van Moerbeke \cite{AvM4})}
 The integrals $I_n(t,s)$, with $I_0=1$,
  \bea
 \tau_n(t,s)= \det m_n 
 =\frac{1}{n!}I_n:&=&
\frac{1}{n!} \int_{(S^1)^{n}}|\Dt_n(z)|^{2}
 \prod_{k=1}^n
\left(e^{\sum_1^{\iy}(t_i z_k^i-s_iz_k^{-i})}
 \frac{dz_k}{2\pi i z_k}\right)\nonumber\\
&=&\int_{U(n)}e^{\sum_1^{\iy}Tr (t_iM^i-s_i\bar
M^i)}dM \nonumber  \\ &=& \sum_{\{\mbox{\footnotesize
Young diagrams }\lambda~ |~ \hat\lambda_1\leq n\}}{\bf
s} _{\lambda}(t){\bf s}_{\lambda}(-s),
 \eea

\vspace{-.1cm}

 \noindent satisfy:
\begin{description}
  \item[(i)]  a SL(2,$\BZ$)-algebra of three \underline{Virasoro
  constraints} (2.3.2):
 \be\hspace{-.5cm}
\BJ_{k,n}^{(2)}(t,n)-
  \BJ_{-k,n}^{(2)}(-s,n) -k\left( \theta
\BJ_{k,n}^{(1)}(t,n)+(1-\theta)\BJ_{-k,n}^{(1)}(-s,n)
   \right)
   I_n(t,s)=0,
 \ee
\vspace{-.9cm}

\hspace{6cm}for $\left\{\begin{array}{l} k =-1
,~\theta= 0\\ k=0 , ~~\theta~~ {arbitrary}\\
  k=1 ,~
\theta= 1\end{array}\right.$

  \item[(ii)] \underline{2d-Toda identities}: The matrices
  $L_1$ and $L_2$, defined above, satisfy the 2-Toda
  lattice equations (3.3.8); in particular:
   $$ \frac{\pl^2}{\pl
s_1\pl t_1}\log\tau_{n}=
-\frac{\tau_{n-1}\tau_{n+1}}{\tau_n^2} $$
   and
\be
\frac{\pl^2}{\pl s_2\pl
t_1}\log\tau_{n}=-2\frac{\pl}{\pl
s_1}\log\frac{\tau_{n}}{\tau_{n-1}}~.\frac{\pl^{2}}{\pl
s_{1}\pl t_{1}}\log\tau_{n}- \frac{\pl^{3}}{\pl
s_{1}^2\pl t_{1}}\log\tau_{n}, \ee
 the first being equivalent to the \underline{\sl discrete sinh-Gordon
 equation}
$$
   \frac{\pl^2 q_{n}}{\pl
t_1 \pl s_1}  =e^{q_{n}-q_{n -1}}- e^{q_{n +1}-q_{n
}},~~~\mbox{where}~~q_n=\log \frac{\tau_{n+1}}{\tau_n}
.$$

  \item[(iii)] \underline{The Toeplitz
  lattice}
: The 2-Toda lattice solution is a very special one,
namely the matrices $L_1$ and $L_2$ have a "rank 2"
structure, given by (3.4.5), whose $x_n$'s and $y_n$'s
equal\footnote{Remember ${\bf s}(t_1,t_2,...)$ are
elementary Schur polynomials. }:
 \bea
 \lefteqn{x_n (t,s)}\nonumber\\&=&\frac{1}{\tau_n}\int _{U(n)}
 {\bf s}_n(-\Tr M,-\frac{1}{2}\Tr M^2,-\frac{1}{3}\Tr M^3,...)
 e^{\sum_1^{\iy} \Tr (t_iM^i-s_i\bar
M^i)} dM \nonumber\\ &=&\frac{{\bf
s}_n(-\frac{\pl}{\pl t_1},-\frac{1}{2}\frac{\pl}{\pl
t_2},-\frac{1}{3}\frac{\pl}{\pl
t_3},...)\tau_n(t,s)}{\tau_n(t,s)}=p_n^{(1)}(t,s;0)
\nonumber\\
 \lefteqn{y_n(t,s)}\nonumber\\
  &=&\frac{1}{\tau_n}\int _{U(n)}
 {\bf s}_n(-\Tr \bar M,-\frac{1}{2}\Tr \bar M^2,
  -\frac{1}{3}\Tr \bar M^3,...)
 e^{\sum_1^{\iy} \Tr (t_iM^i-s_i\bar
M^i)} dM \nonumber
\\ &=&\frac{{\bf s}_n(\frac{\pl}{\pl
t_1},\frac{1}{2}\frac{\pl}{\pl
t_2},\frac{1}{3}\frac{\pl}{\pl
t_3},...)\tau_n(t,s)}{\tau_n(t,s)}
 =p_n^{(2)}(t,s;0),
  \eea
and satisfy the following integrable Hamiltonian
system \bea \frac{\pl x_n}{\pl
t_i}=(1-x_ny_n)\frac{\pl H^{(1)}_i}{\pl y_n}  &~~~~~&
\frac{\pl y_n}{\pl t_i}=-(1-x_ny_n)\frac{\pl
H^{(1)}_i}{\pl x_n} \nonumber \\ \frac{\pl x_n}{\pl
s_i}=(1-x_ny_n)\frac{\pl H^{(2)}_i}{\pl y_n}  &~~~~~&
\frac{\pl y_n}{\pl s_i}=-(1-x_ny_n)\frac{\pl
H^{(2)}_i}{\pl x_n}
  , \\ &&\hspace{2.5cm}(\mbox{\bf Toeplitz lattice})
 \nonumber \eea
with \underline{initial condition}
$x_n(0,0)=y_n(0,0)=0$ for $n\geq 1$ and
\underline{boundary} \underline{condition}
$x_0(t,s)=y_0(t,s)=1$. The traces
 $$ H^{(k)}_i=-
\frac{1}{i}Tr~L _k^i,~~i=1,2,3,...,~~k=1,2. $$ of the
matrices ${ L}_i$ in (3.4.5) are integrals in
involution with regard to the symplectic structure $
\omega := \sum_0^{\iy} (1-x_ky_k)^{-1}dx_k \wedge
dy_k.$ The Toeplitz nature of $m_{\iy}$ leads to
identities between $\tau$'s, the simplest being (due
to Hisakado \cite{H}) :
\be
\left(1+\frac{\pl^2}{\pl s_1 \pl t_1}\log \tau_{n+1}
 \right)\left(1+\frac{\pl^2}{\pl s_1 \pl t_1}\log \tau_{n}
 \right)
=-\frac{\pl}{\pl t_{1}}\log
\frac{\tau_{n+1}}{\tau_n}\frac{\pl}{\pl s_{1}}\log
\frac{\tau_{n+1}}{\tau_n}. \ee

\end{description}

\end{theorem}
 \remark The first equation in the hierarchy above
reads:
 \bean \frac{\pl x_n}{\pl
t_1}=x_{n+1}(1-x_ny_n) &~~~~~& \frac{\pl y_n}{\pl
t_1}=-y_{n-1}(1-x_ny_n)  \\ \frac{\pl x_n}{\pl
s_1}=x_{n-1}(1-x_ny_n)  &~~~~~& \frac{\pl y_n}{\pl
s_1}=-y_{n+1}(1-x_ny_n) . \\
 \eean

\subsubsection{Sketch of Proof}

The identity (3.4.6) between the determinant and the
moment matrix uses again the Vandermonde identity
(3.3.12),
  \bea
\lefteqn { \int_{U(n)}e^{\sum_1^{\iy}Tr
(t_iM^i-s_i\bar M^i)}  dM}\nonumber\\
 &=&
 \int_{(S^1)^{n}}|\Dt_n(z)|^{2}
 \prod_{k=1}^n
\left(e^{\sum_1^{\iy}(t_i z_k^i-s_iz_k^{-i})}
 \frac{ dz_k}{2\pi i z_k}\right)\nonumber\\
&=&
 \int_{(S^1)^{n}}\Dt_n(z)\Dt_n(\bar z)
 \prod_{k=1}^n
\left(e^{\sum_1^{\iy}(t_i z_k^i-s_iz_k^{-i})}
 \frac{dz_k}{2\pi i z_k}\right)\nonumber\\
&=&
 \int_{(S^1)^{n}}\sum_{\sigma\in S_n}\det\left(z_{\sigma(m)}^{\ell-1}
 \bar z_{\sigma(m)}^{m-1}  \right)_{1\leq \ell,m\leq n}
 \prod_{k=1}^n
\left(e^{\sum_1^{\iy}(t_i z_k^i-s_iz_k^{-i})}
 \frac{dz_k}{2\pi i z_k}\right)
 \nonumber\\
&=&
 \sum_{\sigma\in S_n} \det\left(\oint_{S^1}z_{k}^{\ell-1}
 \bar z_{k}^{m-1}
e^{\sum_1^{\iy}(t_i z_k^i-s_iz_k^{-i})}
 \frac{dz_k}{2\pi i z_k}\right)_{1\leq \ell,m\leq n}
  \nonumber\\
  &=&
 n!\det\left(\oint_{S^1}z^{\ell-m}
e^{\sum_1^{\iy}(t_i z^i-s_iz^{-i})}
 \frac{dz}{2\pi i z}\right)_{1\leq \ell,m\leq n}
 =n!\det m_n(t,s)=n! \tau_n(t)\nonumber
 \eea

Using $z^{k\top}=z^{-k}$ (see (3.4.2)), one shows that
the polynomials $p^{(1)}_{n+1}(z)-zp_n^{(1)}(z)$ and
$p^{(1)}_{n+1}(0)z^np_n^{(2)}(z^{-1})$ are
perpendicular to the monomials $z^0,z^1,...,z^n$ and
that they have the same $z^0$-term; one makes a
similar argument, by dualizing $1\leftrightarrow 2$.
Therefore, we have the Hisakado identities between the
following polynomials:
 \bea
p^{(1)}_{n+1}(z)-zp_n^{(1)}(z)&=&p^{(1)}_{n+1}(0)z^np_n^{(2)}(z^{-1})\nonumber\\
p^{(2)}_{n+1}(z)-zp_n^{(2)}(z)&=&p^{(2)}_{n+1}(0)z^np_n^{(1)}(z^{-1}).
\eea The rank $2$ structure (3.4.5) of $L_1$ and
$L_2$, with $x_n=p^{(1)}_n(t,s;0)$ and
$y_n=p^{(2)}_n(t,s;0)$, is obtained by taking the
inner-product of $p^{(1)}_{n+1}(z)-zp_n^{(1)}(z)$ with
itself, for different $n$ and $m$, and using the fact
that $ zp_n^{(1)}(z)=L_1 p_n^{(1)}(z)$.

To check the first equation in the hierarchy (see
remark at the end of theorem 3.5), consider, from
(3.4.9),
 \bean
\frac{\pl x_n}{\pl t_1}
 &=& \left. \frac{\pl p_n^{(1)}(t,s;z)}{\pl t_1}\right|_{z=0}
 \\
 &=& \left. -\left((L_1)_-p^{(1)}\right)_n \right|_{z=0}
  ~,~~~ \mbox{using (3.3.15),}\\
 &=& h_n p^{(1)}_{n+1}(t,s;0)\sum_{i=0}^{n-1} \frac
 { p^{(1)}_{i}(t,s;0) p^{(2)}_{i}(t,s;0)}{h_i}
  ~,~~~ \mbox{using (3.4.5),}\\
 &=& h_n x_{n+1} \sum_{i=0}^{n-1} \frac
 { x_{i} y_{i}}{h_i}\\
 &=&  h_n x_{n+1} \sum_{i=0}^{n-1}\left(
 \frac{1}{h_i}-\frac{1}{h_{i-1}}\right)~,~~~
 \mbox{using}~ \frac{h_i}{h_{i-1}}=1-x_iy_i,\\
 &=& x_{n+1}\frac{h_n}{h_{n-1}}\\
 &=& x_{n+1}  (1-x_ny_n),
 \eean
 and similarly for the other coordinates.
 From (3.3.7) and (3.4.5), upon making the products of the corresponding
 diagonal entries of $L_1$ and $hL_2^{\top}h^{-1}$, one
 finds (3.4.11):
\bean \frac{\pl}{\pl t_{1}}\log
\frac{\tau_{n+1}}{\tau_n}\frac{\pl}{\pl s_{1}}\log
\frac{\tau_{n+1}}{\tau_n}=-x_{n+1}y_n x_ny_{n+1}&=&-
x_n y_n  x_{n+1} y_{n+1}\\
  &=&-\left(1-\frac{h_n}{h_{n-1}}\right)
 \left (1-\frac{h_{n+1}}{h_{n}}\right).
\eean \qed


\section{Ensembles of finite random matrices}

\subsection{PDE's defined by the probabilities in
Hermitian, symmetric and symplectic random ensembles}

As used earlier, the disjoint union
$E=\cup_1^{2r}~[c_{2i-1},c_{2i}]\subset \BR$, and the
weight $\rho(z)=e^{-V(z)}$, with
$-\rho'/\rho=V^{\prime}=g/f$ define
 an algebra of differential operators, ($k \in \BZ$)
 $${\cal
  B}_k=\sum_1^{2r}c_i^{k}f(c_i)\frac{\pl}{\pl c_i}.$$
The aim of this section is to find PDE's for the
following probabilities in terms of the boundary
points $c_i$ of $E$ (see (1.1.9), (1.1.11) and
(1.1.18)), i.e.
 \bea
 P_n(E):&=& P_n(\mbox{ all
spectral points of }~M \in E) \nonumber\\
  &=&\frac{
  \int_{{\cal H}_n(E),~{\cal S}_n(E)~\mbox{\tiny{or}}
  ~{\cal T}_n(E)}e^{-tr ~V(M)}dM}
  {\int_{{\cal H}_n (\BR),~{\cal S}_n(\BR)~\mbox{\tiny{or}}
  ~{\cal T}_n(\BR)}e^{-tr ~ V(M)}dM}
  \\ &=&
  \frac{\int_{E^n}|\Dt_n(z)|^{\beta}\prod_{k=1}^n
e^{-V(z_k)}dz_k}{\int_{\BR^n}|\Dt_n(z)|^{\beta}\prod_{k=1}^n
e^{-V(z_k)}dz_k} ~~~~~ \beta=2,1,4
~~\mbox{respectively}, \nonumber \eea
 involving the classical weights below. In
 anticipation, the equations obtained in
 Theorems 4.1, 4.2 and 4.3 are closely
 related to three of the six Painlev\'e differential
 equations:
 $$
 \begin{array}{l|l|l}
 \mbox{weight}&\rho(z)& \mbox{Painlev\'e } \\
 \hline
 \mbox{Gauss} & e^{-bz^2} & ~~~~ IV\\
  \mbox{Laguerre} & z^a e^{-bz} & ~~~~V \\
  \mbox{Jacobi } &   (1-z)^a(1+z)^b &~~~~VI\\
  \end{array}
  $$
  For $\beta=2$, the
probabilities satisfy partial differential equations
in the boundary points of $E$, whereas in the case
$\beta=1,4$, the equations are inductive. Namely, for
$\beta=1$ (resp. $\beta=4$) , the probabilities
$P_{n+2}$ (resp. $P_{n+1}$) are given in tems of
$P_{n-2}$ (resp. $P_{n-1}$) and a differential
operator acting on $P_n$.
  The weights above involve the parameters $\beta, a$,
$b$ and $$ \delta^{\beta}_{1,4}:=2
\left(\left(\frac{\beta}{2}\right)^{1/2}
-\left(\frac{\beta}{2}\right)^{-1/2}  \right)^2=
\left\{ \begin{array}{l} 0 ~~ \mbox{for}~~\beta=2\\  1
~~ \mbox{for}~~\beta=1,4.\\
\end{array}\right.
 $$
As a consequence of the duality (2.1.12) between
$\beta$-Virasoro generators under the map $\beta
\mapsto 4/\beta$, and the equations (2.1.7), the PDE's
obtained have a remarkable property: the coefficients
$Q$ and $Q_i$ of the PDE's are functions in the
variables $n,\beta,a,b$, having the invariance
property under the map
 $$
  n \rightarrow -2
n,~a \rightarrow -\frac{  a}{2},~b \rightarrow -\frac{
b}{2} ;
 $$
 to be precise,
  \be\left.
Q_i(-2n,\beta ,-\frac{a}{2},-\frac{b}{2}
)\right|_{\beta =1}=\left.Q_i(n, \beta,a,b
)\right|_{\beta =4}.
 \ee

The results in this section are mainly due to
Adler-Shiota-van Moerbeke \cite{ASV1} for $\beta=2$
and Adler-van Moerbeke \cite{AvM3} for $\beta=1,4$.
For more detailed references, see the end of section
4.2.

\subsubsection{Gaussian Hermitian, symmetric and symplectic
ensembles}

Given the disjoint union $E$ and the weight $e^{-b
z^2}$, the differential operators ${\cal
  B}_k$ take on the form
   $${\cal
  B}_k=\sum_1^{2r}c_i^{k+1}\frac{\pl}{\pl c_i}.$$
Define the {\em invariant} polynomials
 $$
  Q= 12
b^2n\left(n+1-\frac{2}{\beta}\right)~~~~\mbox{ and
}~~~~ Q_2=4 (1+\delta^{\beta}_{1,4}) b
\left(2n+\delta^{\beta}_{1,4}
 (1-\frac{2}{\beta})\right).
 $$

\vspace{.3cm}

\begin{theorem}  The following probabilities ($\beta=2,1,4$)
 \be P_n(E)= \frac{\int_{E^n}|\Dt_n(z)|^{\beta}\prod_{k=1}^n
e^{-b
z_k^2}dz_k}{\int_{\BR^n}|\Dt_n(z)|^{\beta}\prod_{k=1}^n
e^{-b z_k^2}dz_k}
,  \ee
  satisfy the PDE's ($F:=F_n=\log P_n$):

   \bigbreak

  \noindent$\displaystyle{\delta^{\beta}_{1,4} Q
\left(\frac{P_{n-{2 \atop 1}}P_{n+{2 \atop 1}
}}{P_n^2}-1\right)}  \hspace{1.8cm}\mbox{with index}
\left\{
\begin{array}{l} 2 ~~ \mbox{when $n$ even and}~~\beta=1\\
  1 ~~\mbox{when $n$ arbitrary and}
~\beta=4\\
\end{array}\right.   $
 \be
=\left({\cal B} _{-1}^4+(Q_2+6{\cal B}_{-1}^2F){\cal
B}_{-1}^2+4(2-\delta^{\beta}_{1,4})\frac{b^2}{\beta}(3{\cal
B}^2_0-4
 {\cal B}_{-1}{\cal B}_1+6{\cal B}_0)\right)F.
\ee

\end{theorem}

\subsubsection{Laguerre Hermitian, symmetric and symplectic
ensembles} Given the disjoint union $E\subset \BR^+$
and the weight $z^{a}e^{-b z}$, the ${\cal B}_k$ take
on the form
  $${\cal
  B}_k=\sum_1^{2r}c_i^{k+2}\frac{\pl}{\pl c_i}.$$
Define the polynomials, also respecting the duality
(4.1.2),
 \bean
  Q&=&\left\{\begin{array}{ll}
\displaystyle{\frac{3}{4} n(n-1)( n+2a)(
n+2a+1)},~~~\mbox{for}~~\beta=1\\  \\ \displaystyle{
\frac{3}{2}n(2n+1)(2n+a)(2n+a-1)}
,~~~\mbox{for}~~\beta=4
\end{array}\right.\\
Q_2&=&\left(
  3\beta n^{2}-\frac{a^2}{\beta}+6a n+4(1-\frac{\beta}{2})
  a
  +3\right)\delta^{\beta}_{1,4}+(1-a^2)(1-\delta^{\beta}_{1,4})\\
  Q_1&=&\left(
  \beta n^{2}+2a n+(1-\frac{\beta}{2}) a\right),~~~~
  Q_0= b (2-\delta^{\beta}_{1,4})(n+
 \frac{a}{\beta}).
  \eean

\vspace{.1cm}

\begin{theorem} 
   The following probabilities
 \be P_n(E)=\frac{\int_{E^n}|\Dt_n(z)|^{\beta}\prod_{k=1}^n
z_k^{a} e^{-b z_k}dz_k
}{\int_{\BR_+^n}|\Dt_n(z)|^{\beta}\prod_{k=1}^n
z_k^{a} e^{-b z_k}dz_k}
 \ee
  satisfy the PDE
  \footnote{with the same convention on the
  indices
  $n\pm 2$ and $n\pm 1$, as in (4.1.4)}
  : ($F:=F_n=\log P_n$)

\noindent $\displaystyle{
 \delta^{\beta}_{1,4} Q\left(\frac{P_{n-{2\atop 1}}P_{n+{2\atop 1}}}{P_n^2}-1 \right)}$

$$ =\Bigl( {\cal B}_{-1}^4 -2(\delta^{\beta}_{1,4}
+1){\cal B}_{-1}^3 + (Q_2+ 6 {\cal
B}_{-1}^2F-4(\delta^{\beta}_{1,4}+1) {\cal B}_{-1}F)
{\cal B}_{-1}^2 -3 \delta^{\beta}_{1,4} (Q_1-{\cal
B}_{-1}F)
  {\cal B}_{-1}$$
  \vspace{-.8cm}\be +\frac{b^2}{\beta}(2-\delta^{\beta}_{1,4})(3
 {\cal B}^2_{0}-4{\cal B}_{1}{\cal B}_{-1}
-2{\cal B}_{1})
  +Q_0(2{\cal B}_{0}{\cal B}_{-1}
    -{\cal B}_{0}) \Bigr)
   F.
 \ee
\end{theorem}

\vspace{.1cm}

\subsubsection{Jacobi Hermitian, symmetric and symplectic
ensembles}

  In terms of $E\subset [-1,+1]$ and the Jacobi weight
$(1-z)^{a}(1+z)^{b}$, the differential
 operators ${\cal B}_k$ take on the form
   $${\cal
  B}_k=\sum_1^{2r}c_i^{k+1}(1-c_i^2)\frac{\pl}{\pl c_i}.$$
 Introduce the following variables, which themselves have
 the invariance property (4.1.2)
($b_0=a-b,~b_1=a+b$):
 $$
  r=\frac{4}{\beta}(b_0^2+(b_1+2-\beta)^2)
%
~~~~~~s=\frac{4}{\beta}b_0(b_1+2-\beta)
$$ $$ q_n= \frac{4}{\beta}(\beta n+b_1 +2-\beta)(\beta
n
 +b_1),
$$
 and the following {\em invariant} polynomials in $q,r,s$:
 \bean
 Q&=&\frac{3}{16}\left(
(s^{2}-qr+q^{2})^{2}-4(rs^{2}-4q
 s^{2}-4s^{2}+q^{2}r)\right)
\\
  Q_1&=& 3s
 ^{2}-3qr-6r+2q^{2}+23q+24\\
 Q_2&=&3qs^{2}
 +9s^{2}-4q^{2}\,r+2qr+4q^{3}+10q^{2}
 \\
 Q_3&=&3qs^{2}+6s^{2}-3q^{2}
 r+q^{3}+4q^{2}
 \\
 Q_4&=& 9s
 ^{2}-3qr-6\,r+q^{2}+22q+24=Q_1+(6s^2-q^2-q).\\
 \eean
\vspace{-2cm} \be \ee

\begin{theorem} 
   The following probabilities
 \be P_n(E)=\frac{\int_{E^n}|\Dt_n(z)|^{\beta}\prod_{k=1}^n
(1-z_k)^{a}(1+z_k)^{b}dz_k
}{\int_{[-1,1]^n}|\Dt_n(z)|^{\beta}\prod_{k=1}^n
(1-z_k)^{a}(1+z_k)^{b}dz_k}
 \ee
  satisfy the PDE ($F=F_n=\log P_n$):

  \bigbreak

\noindent \underline {for $\beta=2$:}

 \bea &&\Bigl(2 {\cal
B}_{-1}^4+(q-r+4){\cal B}^2_{-1}-(4{\cal B}_{-1}F-s)
{\cal B}_{-1}+3q {\cal B} _{0}^2 - 2q {\cal
B}_{0}+8{\cal B}_0 {\cal B}_{-1}^2 \nonumber \\
&&-4(q-1){\cal B}_{1}{\cal B}_{-1} +(4{\cal
B}_{-1}F-s) {\cal B}_1 +2(4{\cal B}_{-1}F-s) {\cal
B}_{0} {\cal B}_{-1}+
 2 q {\cal B}_2 \Bigr)F
\nonumber\\ && \hspace{2cm}
  +4{\cal B}^2_{-1}F\left( 2{\cal B}_{0}F
   +3{\cal B}^2_{-1}F\right)
=0\eea

\noindent \underline {for $\beta=1,4$:}

\bigbreak

{\footnotesize

 $\displaystyle{Q\left(\frac{P_{n+{2\atop
1}}P_{n-{2\atop 1}}} {P_n^2}-1\right)}$ \bea
&&=(q+1)\Bigl(4q {\cal
 B}^4_{-1}
 +12(4{\cal B}_{-1}F-s)
 {\cal B}^3_{-1}+2\left(q+12\right)(4{\cal B}_{-1}F-s)
 {\cal B}_{0}{\cal B}_{-1} \nonumber\\
 &&+3q^{2}
 {\cal B}_{0}^2 -4\left(q-4\right)q {\cal B}_{1}{\cal
B}_{-1}+q(4{\cal B}_{-1}F-s)  {\cal B}_{1}+20q{\cal
B}_{0}{\cal B}^2_{-1}+2q^{2}{\cal B}_{2}\Bigr) F
 \nonumber\\
 &&+\Bigl(Q_2{\cal B}^2_{-1}-sQ_1{\cal
 B}_{-1}+Q_3{{\cal B}_0}\Bigr)F +48({\cal B}_{-1}F)^{4}-48s({\cal
B}_{-1}F)^{3}+2Q_4({\cal B}_{-1}F)^{2} \nonumber\\
 &&+12\,q^{2}({{\cal
B}_0}F)^{2}+16\,q\,\left(2\,q-1\right)
 ({\cal B}^2_{-1}F)({{\cal B}_0}F)+24\left(q-1\right)q({\cal B}^2_{-1}F)^{2}
\nonumber \\
 &&+24  \Bigl(2{\cal B}_{-1}F
 -s\Bigr)\Bigl(  (q+2)
  {{\cal B}_0}F
 + (q+3){\cal B}^2_{-1}F\Bigr){\cal B}_{-1}F.
\eea }

\end{theorem}

\vspace{-.7cm}

The 
{\it Proof\/} of these three theorems will be sketched
in subsections 4.3, 4.4 and 4.5.


\subsection{ODE's, when $E$ has one boundary point}
Assume the set E consists of one boundary point $c=x$,
besides the boundary of the full range; thus, setting
respectively $E=[-\iy,x], E=[0,x], E=[-1,x]$ in the
PDE's (4.1.4), (4.1.6) and (4.1.9), (4.1.10), leads to
the equations in $x$ below. Notice that, for $\beta
=2$, the equations obtained are ODE's and, for
$\beta=1,4$, these equations express $P_{n+2}$ in
terms of $P_{n-2}$ and a differential operator acting
on $P_n$:

\medbreak

 \noindent (1) {\em Gauss ensemble} ~
 ($\beta=2,1,4$) :
 $f_n(x)=\frac{d}{dx}\log P_n(\max_{i} ~
 \lb_i \leq x)$
  satisfies
 \bea \lefteqn{
 \delta^{\beta}_{1,4}Q\left(\frac{P_{n-{2\atop 1}}P_{n+{2\atop 1}}}{P_n^2}-1 \right)
 }
 \\
  &=& f_n^{\prime\prime\prime}+6 f_n^{\prime 2 }+
\left(4\frac{b^2x^2}{\beta}(\delta^{\beta}_{1,4}-2)
 +Q_2\right)f_n^{\prime
}-4\frac{b^2x}{\beta}(\delta^{\beta}_{1,4}-2)f_n
.\nonumber\eea

\noindent (2) {\em Laguerre ensemble} ($\beta=2,1,4$):
 $f_n(x)=x\frac{d}{dx}\log P_n(\max_{i} ~
 \lb_i \leq x
 )$ (with all eigenvalues $\lb_i\geq 0$)
satisfies:
 \bea
 \lefteqn{ \delta^{\beta}_{1,4}
Q\left(\frac{P_{n-{2\atop 1}}P_{n+{2\atop
1}}}{P_n^2}-1 \right)
 -\left(3\delta^{\beta}_{1,4} f_n
-\frac{b^2 x^2 }{\beta}(\delta^{\beta}_{1,4}-2)
-Q_0x-3 \delta^{\beta}_{1,4} Q_1\right)f_n
}\nonumber\\
 &=&
 x^3f_n^{\prime\prime\prime}-(2\delta^{\beta}_{1,4}
-1) x^2f_n^{\prime\prime}+6x^2f_n^{\prime 2}
\nonumber\\ &&~~-x\left(4(\delta^{\beta}_{1,4}+1)f_n-
\frac{b^2x^2}{\beta}(\delta^{\beta}_{1,4}-2)
-2Q_0x-Q_2+2\delta^{\beta}_{1,4}+1\right) f_n^{\prime}
.\nonumber\\ \eea

\noindent (3) {\em Jacobi ensemble}:
$f:=f_n(x)=(1-x^2)\frac{d}{dx}\log P_n(\max_{i} \lb_i
\leq x)$ (with all eigenvalues $-1 \leq \lb_i\leq 1$)
satisfies:
\begin{itemize}
  \item

 for \underline{$\beta=2$}: \bea
 &&\hspace{-1.3cm} 2(x^2-1)^2f^{\prime\prime\prime}
 +4(x^2-1)\left(xf^{\prime\prime}
 -3f^{\prime 2}\right)
 +\left(16
xf-q(x^{2}-1)-2sx-r \right)f^{\prime}
 \nonumber\\&&~~~~~-f\left(4f-qx-s \right)=0,
 \eea

\item
for \underline{$\beta=1,4$}:

$\displaystyle{Q\left(\frac{P_{n+{2\atop
1}}P_{n-{2\atop 1}}} {P_n^2}-1\right)}$
 {\footnotesize
\bea &=& 4(q+1)(x^2-1)^{2}\Bigl(-q(x^2-1)
f^{\prime\prime\prime}+
 (12f -qx-3s)
  f^{\prime\prime}+6 q(q-1)
 f^{\prime 2}\Bigr)\nonumber\\
  &&-(x^2-1)f^{\prime}\Bigl(24f(q+3)(2f-s)+8fq
  (5q-1)x-
   q(q+1)(qx^{2}+2sx+8)+
 Q_2
 \Bigr)
 \nonumber\\
 && +f\Bigl(48f^{3}+48f^{2}(qx+2x-s)+2f\left(8q^{2}x^{2}
 +2qx^{2}-12qsx-24s
 x+Q_4\right)\nonumber\\
  &&~~~~~-q(q+1)x(3qx^{2}+sx-2
 qx-3q)+Q_3x- Q_1s\Bigr). \eea
 }

\end{itemize}

For $\beta=2$, the term containing the ratio
$\frac{P_{n+{2\atop 1}}P_{n-{2\atop 1}}} {P_n^2}-1$ on
the left hand side of (4.2.1), (4.2.2) and (4.2.4)
vanishes and one thus obtains the ODE's:
\begin{itemize}
  \item {\bf Gauss}: $f_n(x):=\frac{d}{d x}
  \log P_n(\max_i \lb_i \leq x)$
  satisfies:
$$
 f^{\prime\prime\prime} +  6
~f^{\prime 2}+4b (2n-bx^2) f^{\prime} + 4b^2 x~f =0$$

  \item
  {\bf Laguerre}: $f_n(x):=x\frac{d}{d x}
   \log P_n(\max_i \lb_i\leq x)$
satisfies
\begin{eqnarray*}
\hspace{-1cm} x^2f^{\prime\prime\prime} +
xf^{\prime\prime}
 + 6
xf^{ \prime 2} - 4 ff^{\prime}-((a-bx)^2 - 4nbx )
f^{\prime}- b(2n + a - bx) f = 0 .
\end{eqnarray*}
  \item {\bf Jacobi}: $f_n(x)=(1-x^2)\frac{d}{dx}\log P_n(\max_{i} \lb_i \leq
x)$ satisfies:

  \bean &&\hspace{-2cm}{ 2(x^2-1)^2f^{\prime\prime\prime}
 +4(x^2-1)\left(xf^{\prime\prime}
 -3f^{\prime 2}\right)
 +\left(16
xf-q(x^{2}-1)-2sx-r \right)f^{\prime}  }
 \nonumber\\&&\hspace{5cm}-f\left(4f-qx-s \right)=0.
 \eean
\end{itemize}
 Each of these three equations is of the Chazy form
 (see section 9)
 \be
f^{\prime \prime\prime}+\frac{P'}{P}f^{
\prime\prime}+\frac{6}{P}f^{\prime
2}-\frac{4P'}{P^2}ff' +\frac{P^{\prime\prime}}{P^2}
f^2
+\frac{4Q}{P^2}f'-\frac{2Q'}{P^2}f+\frac{2R}{P^2}=0,
\ee
 with $c=0$ and $P,Q,R$ having the form:
  $$
\begin{array}{llll}
Gauss & P(x)=1&~4Q(x)=-4b^2x^2+8bn&~R=0 \\ Laguerre &
P(x)=x&~4Q(x)=-(bx-a)^2+4bnx
&~R=0
\\ Jacobi &
P(x)=1-x^2&~4Q(x)=-\frac{1}{2}(q(x^2-1)+2sx+r)&~
 R=0 \end{array}
 $$

 Cosgrove shows such a third order equation (4.2.5) in $f(x)$,
 with $P(x)$, $Q(x)$, $R(x)$ of respective degrees $3,2,1$,
has a first integral (9.0.2), which is second order in
$f$ and quadratic in $f^{\prime\prime}$, with an
integration constant $c$. Equation (9.0.2) is a master
Painlev\'e equation, containing the 6 Painlev\'e
equations. If $f(x)$ satisfies the equations above,
then the new (renormalized) function $g(z)$, defined
below,

$$
\begin{array}{ll}
Gauss & g(z)=b^{-1/2}f(zb^{-1/2})+\frac{2}{3}nz\\
Laguerre &
g(z)=f(z)+\frac{b}{4}(2n+a)z+\frac{a^2}{4}\\ Jacobi &
g(z):=-\frac{1}{2}f(x)|_{x=2z-1}-\frac{q}{8}z+\frac{q+s}{16}
 \end{array}
 $$
satisfies the canonical equations, which then can be
transformed into the standard Painlev\'e equations;
these canonical equations are respectively:

\begin{itemize}
  \item $  g^{\prime\prime 2}=-4g^{\prime 3
}+4(zg^{\prime}-g)^2+A_1g^{\prime}+A_2$
\hspace{3cm}({\bf Painlev\'e IV})

  \item $ (zg^{\prime\prime })^2=(zg^{\prime}-g)\Bigl(-4g^{\prime 2
}+A_1(zg^{\prime}-g)+A_2\Bigr)+A_3g^{\prime}+A_4
 $

 \hspace{10cm}({\bf
Painlev\'e V})
  \item $ (z(z-1)g^{\prime\prime
})^2=(zg^{\prime}-g)\Bigl(4g^{\prime
2}-4g^{\prime}(zg^{\prime}-g)
 +A_2\Bigr) +A_1g^{\prime
2}+A_3g^{\prime}+A_4 $

\hspace{10cm}({\bf Painlev\'e
VI})
\end{itemize}

\vspace{-1cm}

\noindent with respective coefficients
 \begin{itemize}
  \item
  $A_1=3\left( \frac{4n}{3}\right)^2,~A_2=-\left(
  \frac{4n}{3}\right)^3 $

\item
$
 A_1=b^2,~A_2=b^2((n+\frac{a}{2})^2+\frac{a^2}{2})
 ,~A_3=-a^2b(n+\frac{a}{2}),~
 A_4=\frac{(ab)^2}{2}\\
 ~~~~~~~~~~~~~~~~~~~~~~~~~~~~~~~~~~~~~~~~~~~~~
 ~~~~~~~~~~~~~~~~~
 ~~~~ .((n+\frac{a}{2})^2+\frac{a^2}{8})
 \\
 $
  \item

$A_1=\frac{2q+r}{8},~A_2=\frac{qs}{16},~A_3=\frac{(q-s)^2
+2qr}{64},~ A_4=\frac{q}{512}(2s^2+qr).
$

\end{itemize}

 For $\beta=1~\mbox{and}~4$,
 the inductive partial differential equations (4.1.4),
  (4.1.6),
 (4.1.10), and the derived differential equations
 (4.2.1), (4.2.2) and (4.2.4) are due to Adler-van Moerbeke
 \cite{AvM3}. For $\beta =2$ and for general $E$, they
 were first computed by
 Adler-Shiota-van Moerbeke \cite{ASV1}, using the
method of the present paper. For $\beta=2$ and for $E$
having one boundary point, the equations obtained here
coincide with the ones first obtained by Tracy-Widom
in \cite{TW1}, who recognized them to be Painlev\'e IV
and V for the Gaussian and Laguerre distribution
respectively. In his Louvain doctoral dissertation,
J.P. Semengue, together with L. Haine
 \cite{Haine}, were lead to
 Painlev\'e VI for the Jacobi ensemble,
 for $\beta=2$ and $E$
 having one boundary point,
 upon subtracting the Tracy-Widom
 differential equation (\cite{TW1}) from the one
 computed with the
 Adler-Shiota-van Moerbeke method (\cite{ASV1}).
 Cosgrove's (\cite{Cosgrove}) and
Cosgrove-Scoufis's classification (\cite{CS}, (A.3),)
leads directly to these results.

\subsection{Proof of Theorems 4.1, 4.2 and 4.3}

\subsubsection{Gaussian and Laguerre ensembles}

The three first Virasoro equations, as in (2.1.29) and
(2.1.32), are differential equations, involving
partials in $t\in\BC^{\iy}$ and partials ${\cal
D}_1,{\cal D}_2,{\cal D}_3$ in
$c=(c_1,...,c_{2r})\in\BR^{2r}$, for $F:=F_n(t,c)=\log
I_n$; they have the general form:
\be
{\cal D}_kF=\frac{\pl F}{\pl t_k}+\sum_{-1\leq
j<k}\ga_{kj}V_j(F)+\ga_k+\delta_kt_1,\quad k=1,2,3,
\ee with first $V_j(F)$'s given by:
\be
V_j(F)=\sum_{i,i+j\geq 1}it_i\frac{\pl F}{\pl
t_{i+j}}+\frac{\beta}{2}\delta_{2,j}\left(\frac{\pl^2
F}{\pl t_1^2}+\left(\frac{\pl F}{\pl
t_1}\right)^2\right),\quad -1\leq j\leq 2. \ee In
(4.3.1) and (4.3.2), $\beta
>0,\gamma_{kj},\gamma_k, \delta_k$ are arbitrary
parameters; also $\delta_{2j}=0$ for $j\neq 2$ and
$=1$ for $j=2$. The claim is that the equations
(4.3.1) enable one to express all partial derivatives,
\be
 \left.\frac{\pl^{i_1+...+i_k} F(t,c)}{\pl t_1^{i_1}
 ...\pl t_k^{i_k}}\right|_{\cal L},~\mbox{along}~{\cal
 L}:=\{\mbox{all $t_i=0, ~c=(c_1,...,c_{2r})$
 arbitrary}\},
 \ee
uniquely in terms of polynomials in
 $${\cal D}_{j_1}...{\cal D}_{j_r}F(0,c).$$
 Indeed, the method consists of expressing
$\pl F /\pl t_k \bigl|_{t=0}$ in terms of ${\cal
D}_kf\bigl|_{t=0}$, using (4.3.1). Second derivatives
are obtained by acting on ${\cal D}_kF$ with ${\cal
D}_{\ell}$, by noting that ${\cal D}_{\ell}$ commutes
with all $t$-derivatives, by using the equation for
${\cal D}_{\ell}F$, and by setting in the end $t=0$:
{\footnotesize \bea {\cal D}_{\ell}{\cal D}_kF&=&{\cal
D}_{\ell}\frac{\pl F}{\pl t_k} +\sum_{-1\leq
j<k}\ga_{kj}{\cal D}_{\ell}(V_j(F))\nonumber\\
&=&\left(\frac{\pl }{\pl t_k} +\sum_{-1\leq
j<k}\ga_{kj}V_j\right){\cal
D}_{\ell}(F),\quad\mbox{provided $V_j(F)$ does
not}\nonumber
\\ &
&  \hspace{5cm}\mbox{contain non-linear
terms}\nonumber\\ &=&\left(\frac{\pl }{\pl t_k}
+\sum_{-1\leq j<k}\ga_{kj}V_j\right)\left(\frac{\pl
F}{\pl t_{\ell}} +\sum_{-1\leq j<\ell}\ga_{\ell
j}V_j(F)+\dt_{\ell}t_1\right)\nonumber\\
&=&\frac{\pl^2F}{\pl t_k\pl
t_{\ell}}+\mbox{\,lower-weight terms.}\nonumber \eea }
 When the non-linear term is present, it is taken
care as follows:
\begin{eqnarray*} {\cal
D}_{\ell}\left( \frac{\pl F}{\pl t_1} \right)^2
=2\frac{\pl F}{\pl t_1} {\cal D}_{\ell} \frac{\pl
F}{\pl t_1}&=&2\frac{\pl F}{\pl t_1} \frac{\pl }{\pl
t_1} \left( \frac{\pl F}{\pl t_{\ell}}+\sum_{-1\leq
j<{\ell}}\ga_{{\ell}j}V_j(F)+\ga_{\ell}+\delta_{\ell}t_1
\right).
\end{eqnarray*}
Higher derivatives are obtained in the same way. We
only record here, for future use, the few partials
appearing in the KP equation (3.1.6): {\footnotesize
\begin{eqnarray*}
 \frac{\pl^2 F}{\pl
t^2_1}\Biggl|_{\cal L}&=&\left({\cal
D}_1^2-\ga_{10}{\cal
D}_1\right)F+\ga_{10}\ga_1-\dt_1\nonumber\\
 \frac{\pl^4 F}{\pl
t^4_1}\Biggl|_{\cal L}&=&\left({\cal
D}_1^4-6\ga_{10}{\cal D}_1^3+11\ga_{10}^2{\cal
D}_1^2-6\ga_{10}^3 {\cal
D}_1\right)F-6\ga_{10}^2(\dt_1-\ga_1\ga_{10})\nonumber
\\
  \frac{\pl^2 F}{\pl t^2_2}\Biggl|_{\cal
L}&=&\biggl({\cal D}_2^2-2\ga_{20}{\cal
D}_2+\beta\ga_{21} \ga_{32}{\cal D}^2_1-((2\ga_1
+\ga_{10})\ga_{21}\ga_{32}\beta+2\ga_{2,-1}) {\cal
D}_1\\ &&-2\ga_{21}{\cal D}_3\biggr)F
+\beta\ga_{21}\ga_{32}({\cal D}_1F)^2\\ &&
+\beta\ga_{21}\ga_{32}(\ga^2_1+\ga_{10}
\ga_1-\dt_1)+2(\ga_{21}\ga_3+\ga_{20}\ga_2+
 \ga_1\ga_{2,-1})\nonumber\end{eqnarray*}
\begin{eqnarray*}
\frac{\pl^2 F}{\pl t_1\pl t_3}\Biggl|_{\cal
L}&=&\biggl({\cal D}_1{\cal
D}_3-\frac{\beta}{2}\ga_{32}{\cal D}_1^3
+\beta\ga_{32}( \ga_1+2\ga_{10}){\cal
D}_1^2-\frac{3\beta}{2}
\ga_{10}\ga_{32}(2\ga_1+\ga_{10}){\cal D}_1\nonumber\\
& &-3\ga_{1,-1}{\cal D}_2-3\ga_{10}{\cal
D}_3\biggr)F+\frac{3\beta}{2}\ga_{10}\ga_{32}({\cal
D}_1F)^2 -\beta\ga_{32}({\cal D}_1F)({\cal
D}_1^2F)\nonumber\\ & &+
\frac{3}{2}(2\ga_{10}\ga_3+\beta\ga_{32}\ga_{10}(\ga^2_1+\ga_{10}\ga_1-\dt_1)+2\ga_{1,-1}
\ga_2).\nonumber
\end{eqnarray*}  }

\vspace{-1cm}

\subsubsection{Jacobi ensemble}

Here, from the Virasoro constraints (2.1.35), one
proceeds in the same way as before, by forming
$\left.{\cal B}_{i}F\right|_{t=0},~\left.{\cal B}_{i}
{\cal B}_{j}F\right|_{t=0}$, etc..., in terms of $t_i$
partials. For example, from the expressions $
\left.{\cal B}_{-1}F\right|_{t=0}, \left. {\cal
B}^2_{-1}F\right|_{t=0}, \left.{\cal
B}_{0}F\right|_{t=0} ,$
  one extracts $$
 \left.\frac{\pl F}{\pl t_1}\right|_{t=0},
 \left.\frac{\pl^2 F}{\pl t_1^2}\right|_{t=0},
 \left.\frac{\pl F}{\pl t_2}\right|_{t=0}.
$$ 
From the expressions $ \left.{\cal
B}^3_{-1}F\right|_{t=0}, \left.{\cal B}_{0} {\cal
B}_{-1}F\right|_{t=0}, \left.{\cal
B}_{1}F\right|_{t=0} $,
 and using the previous
information, one extracts $$
 \left.\frac{\pl F}{\pl t_3}\right|_{t=0},
 \left.\frac{\pl^2 F}{\pl t_1^3}\right|_{t=0},
 \left.\frac{\pl^2 F}{\pl t_1\pl t_2}\right|_{t=0}.
$$
 Finally, from the expressions $ \left.{\cal
B}_{2}F\right|_{t=0}, \left.{\cal B}_{1} {\cal
B}_{-1}F\right|_{t=0}, \left.{\cal
B}^2_{0}F\right|_{t=0},
 \left.{\cal
B}_{0} {\cal B}^2_{-1}F\right|_{t=0},$ $ \left. {\cal
B}^4_{-1}F\right|_{t=0} $, one deduces
 \be
 \left.\frac{\pl ^4 F}{\pl t_1^4}\right|_{t=0},
 \left.\frac{\pl F}{\pl t_4}\right|_{t=0},
 \left.\frac{\pl^3 F}{\pl t^2_1\pl t_2}\right|_{t=0},
 \left.\frac{\pl^2 F}{\pl t_1\pl t_3}\right|_{t=0},
 \left.\frac{\pl^2 F}{\pl t_2^2}\right|_{t=0}.
\ee This provides all the partials, appearing in the
KP equation (3.1.6).


\subsubsection{Inserting partials into the integrable equation}

From Theorem 3.3, the integrals $I_n(t,c)$ , depending
on  $\beta=2,1,4$, on $t=(t_1,t_2,...)$ and on
 the boundary points
$c=(c_1,...,c_{2r})$ of $E$, relate to
$\tau$-functions, as follows:
 \bea
 I_n(t,c)& =&\int_{E^n}|\Dt_n(z)|^{\beta}\prod_{k=1}^n
\left(e^{\sum_1^{\iy}t_i z_k^i}\rho(z_k)dz_k\right)
\nonumber\\ &=&\left\{
\begin{array}{lll}
 n!\tau_n(t,c),& n~\mbox{arbitrary}, &\beta=2
 \\
n!\tau_n(t,c),& n~\mbox{even}, &\beta=1
 \\ n!
\tau_{2n}(t/2,c), & n~\mbox{arbitrary},&\beta=4
\end{array}
\right. ,\eea
 where $\tau_n(t,c)$ satisfies the KP-like equation
 \be
 12\frac{\tau_{n-2}(t,c)\tau_{n+2}(t,c)}{\tau_n(t,c)^2}
\delta^{\beta}_{1,4}=(KP)_t\log \tau_n(t,c),~~
\left\{\begin{array}{l} n~\mbox{arbitrary for}
~\beta=2\\
 n~\mbox{even for}~ \beta=1,4
   \end{array}\right.
    \ee
 with
  $$
(KP)_t  F:= \left(\left(\frac{\pl}{\pl t_1}
\right)^4+3\left(\frac{\pl}{\pl
t_2}\right)^2-4\frac{\pl^2}{\pl t_1 \pl
t_3}\right)F+6\left(\frac{\pl^2}{\pl t^2_1}F
\right)^2.
 $$

 \noindent \underline{\em Evaluating the left hand
 side of (4.3.6):}
Here $I_n(t)$ will refer to the integral (4.3.5) over
the full range. For \underline{$\beta =2$}, the left
hand side is zero. For \underline{$\beta =1$}, the
left hand side can be evaluated in terms of the
probability $P_n(E)$, as follows: taking into account
$P_n:=P_n(E)=I_n(0,c)/I_n(0)$,
 \bean
\left.12\frac{\tau_{n-2}(t,c)\tau_{n+2}(t,c)}{\tau_n(t,c)^2}
\right|_{t=0}&=&12\frac{(n!)^2}{(n-2)!(n+2)!}
\left.\frac{I_{n-2}(t,c)I_{n+2}(t,c)}{I_n(t,c)^2}
\right|_{t=0}\\&=&12\frac{n(n-1)}{(n+1)(n+2)}
\frac{I_{n-2}(0)I_{n+2}(0)}{I_{n}(0)^2}
 \frac{P_{n-2}P_{n+2}}{P_n^2}
  \\& =&12b_n^{(1)}\frac{P_{n-2}(E)P_{n+2}(E)}{P_n^2(E)}
,\eean
 with $b^{(1)}_n$ given by\footnote{this calculation is based on
Selberg's integrals (see Mehta \cite{Mehta}, p 340).
For instance, in the Jacobi case, one uses \bean
I_n^{(\beta)}&=&\int_{[-1,1]^n} \Delta_{n}
(x)^{\beta}\prod^n_{j=1}(1-x_j)^{a}(1+x_j)^{b}dx_j
\\&=&2^{n(2a+2b+\beta(n-1)+2)/2}
 \prod^{n-1}_{j=0} \frac{\Gamma(a+j\beta /2
+1)\Gamma(b+j\beta /2+1) \Gamma((j+1)\beta /2 +1)}
{\Gamma(\beta /2 +1)\Gamma(a+b +(n+j-1)\beta /2+2)}.
\eean }
 \bea
    b_n^{(1)}
 =\frac{n(n-1)}{(n+2)(n+1)}
\frac{I_{n-2}(0)I_{n+2}(0)}{I_n(0)^2} =\left\{
\begin{array}{l}
 \frac{n(n-1)}{16b^2}  ~~~~\mbox{ (Gaussian)}\\  \\
 \frac{n(n-1)(n+2a)(n+2a+1)}{16b^4}\\ \hspace{2cm}
\mbox{ (Laguerre)}\\
 \frac{Q}{Q_6^{\pm}} ~~~~~~~~~~~  \mbox{
(Jacobi)}\\
\end{array}\right. \nonumber\\
 \eea

\noindent For \underline{\em $\beta=4$}, we have:
  \bean
\left.12\frac{\tau_{2n-2}(t/2,c)
 \tau_{2n+2}(t/2,c)}{\tau_{2n}(t/2,c)^2}
\right|_{t=0}&=&12\frac{(n!)^2}{(n-1)!(n+1)!}
\left.\frac{I_{n-1}(t,c)I_{n+1}(t,c)}{I_n(t,c)^2}
\right|_{t=0}\\&=& 12 \frac{n }{ (n+1)}
\frac{I_{n-1}(0)I_{n+1}(0)}{I_n(0)^2}\frac{P_{n-1}P_{n+1}}{P_n^2}\\
 &=& 12b_n^{(4)}\frac{P_{n-1}(E)P_{n+1}(E)}{P_n^2(E)}
,\eean with
 \bea
 b_n^{(4)}&:=&\frac{(n!)^2}{(n-1)!(n+1)!}
\frac{I_{n-1}(0)I_{n+1}(0)}{I_n^2(0)}
  =\left\{
\begin{array}{l}
 \frac{2n(2n+1)}{4b^2}  ~~~~\mbox{ (Gauss)}\\  \\
 \frac{2n(2n+1)(2n+a)(2n+a-1)}{b^4}\\ \hspace{2cm}
\mbox{ (Laguerre)}\\
 \frac{Q}{Q_6^{\pm}} ~~~~~~~~~~~  \mbox{
(Jacobi)}\\
\end{array}\right. \nonumber\\
 \eea
where $Q$ is precisely the expression appearing
 on the left hand side of (4.1.10), and where $Q_6^{\pm}$ is given by
 {\footnotesize \bea Q^{\pm}_6
 &=& 3
 q\left(q+1\right)\left(q-3\right)\left(q+4\pm4\sqrt{q+1}
 \right)~~~~~~~~  \left\{\begin{array}{lll} +& \mbox{for}& \beta=1\\
                             - & \mbox{for} & \beta=4
                             \end{array}\right.
 \eea}
The exact formulae $ b_n^{(4)}$ and
 $b_{n}^{(1)}$ show they satisfy the duality property (4.1.2):
$$
b_n^{(4)}(a,b,n)=b_{n}^{(1)}(-\frac{a}{2},-\frac{b}{2},-2n).
$$

\noindent \underline{\em Evaluating the right hand
side of (4.3.6)}:
 From section 2.4, it also follows
that $F_n(t;c)=\log I_n(t;c)$ satisfies Virasoro
constraints, corresponding precisely to the situation
(4.3.1), with

\noindent {\em Gaussian ensemble}\footnote{Remember
from section 2.1, $\sigma_1=\beta (n-1)+2$.}: $$
\left\{\begin{array}{l}
\ga_{1,-1}=-\frac{1}{2},\ga_{1,0}=\ga_1=0,\dt_1=-\frac{n}{2}\\
\ga_{2,-1}=0,\ga_{2,0}=-1/2,\ga_{2,1}=0,
\ga_2=-\frac{n}{4} \sigma_1 
 ,\dt_2=0\\
\ga_{3,-1}=-\frac {1}{4} \sigma_1
,\ga_{3,0}=0,\ga_{3,1}=-\frac{1}{2},\ga_{3,2}=\ga_3=0,
\dt_3=-\frac{n}{4}\sigma_1
 .
\end{array}\right.
$$

\noindent {\em Laguerre ensemble}\footnote{Remember
from section 2.1, $\sigma_1=\beta (n-1)+a+2$ and
$\sigma_2=\beta (n-\frac{3}{2})+a+3$.}:
  $\dt_1=\dt_2=\dt_3=0$,
and $$ \left\{\begin{array}{l}
\ga_{1,-1}=0,\ga_{1,0}=-1,\ga_1=-\frac{n}{2}(\sigma_1+a),\\
\ga_{2,-1}=0,\ga_{2,0}=-\sigma_1,\ga_{2,1}=-1,
 \ga_2=-\frac{n}{2}\sigma_1(\sigma_1+a)
,\\ \ga_{3,-1}=0,\ga_{3,0}=-\sigma_1 \sigma_2
,\ga_{3,1}=-\sigma_2, \ga_{3,2}=-1, \ga_3=
 -\frac{n}{2}\sigma_1\sigma_2(\sigma_1+a) .
\end{array}\right.
$$

\noindent {\em Jacobi ensemble}: see (4.3.4).

They lead to expressions for $$
 \left.\frac{\pl ^4 F}{\pl t_1^4}\right|_{t=0},
 \left.\frac{\pl^2 F}{\pl t_2^2}\right|_{t=0},
 \left.\frac{\pl^2 F}{\pl t_1\pl t_3}\right|_{t=0},
  \left.\frac{\pl^2 F}{\pl t_1^2}\right|_{t=0},
$$ in terms of ${\cal D}_k$ and ${\cal B}_k$, which
substituted in the right hand side of (4.3.6) - i.e.
in the KP-expressions - leads to the right hand side
of (4.1.4),(4.1.6), (4.1.9) and (4.1.10). In the
Jacobi case, the right hand side of (4.3.6) contains
the same coefficient $1/Q_6^{\pm}$ as in (4.3.9),
which therefore cancels with the one appearing on the
left hand side; see the expression $b_n^{1,4}$ in
(4.3.7) and (4.3.8).


\section{Ensembles of infinite random matrices: Fredholm
determinants, as $\tau$-functions of the KdV equation}

Infinite Hermitian matrix ensembles typically relate
to the Korteweg-de Vries hierarchy, itself a reduction
of the KP hierarchy; a brief sketch will be necessary.
The {\em KP-hierarchy} is given by $t_n$-deformations
of a pseudo-differential operator\footnote{In this
section, given $P$ a pseudo-differential operator,
$P_+$ and $P_-$ denote the differential and the
(strictly) smoothing part of $P$ respectively. } $L$:
(commuting vector fields)
  \be
   \frac{\pl
L}{\pl t_n}=[(L^n)_+,L],~~L=D+a_{-1}D^{-1}+...,\quad
\mbox{with}\quad D=\frac{\pl}{\pl x}.
 \ee
  Wave and
adjoint wave functions are eigenfunctions
$\Psi^+(x,t;z)$ and $\Psi^{-}(x,t;z)$, depending on $x
\in \BR,~t \in \BC^{\iy},~ z\in \BC$, behaving
asymptotically like (5.0.3) below and satisfying:
\be
\begin{array}{ll}
z\Psi^+=L\Psi^+,&\displaystyle{\frac{\pl\Psi^+}{\pl
t_n}}=(L^n)_+\Psi^+,  \\
z\Psi^-=L^{\top}\Psi^-,&\displaystyle{\frac{\pl\Psi^-}{\pl
t_n}}=-(L^{\top n})_+\Psi^-.
\end{array}
\ee
 According to Sato's theory, $\Psi^+$ and
$\Psi^-$ have the following representation in terms of
a $\tau$-function (see \cite{DJKM}): \bea
\Psi^{\pm}(x,t;z)&=&e^{\pm (xz+\sum_1^{\iy}t_iz^i)}
\frac{\tau(t\mp
[z^{-1}])}{\tau(t)}\nonumber\\&=&e^{\pm
(xz+\sum_1^{\iy}t_iz^i)}(1+O(z^{-1})), ~\mbox{for}~z
\nearrow \iy, \eea where $\tau$ satisfies and is
characterized by the following bilinear relation
 \bea
  \oint
e^{\sum_1^{\iy}
(t_i-t'_i)z^i}\tau(t-[z^{-1}])\tau(t'+[z^{-1}])dz=
0,~~\mbox{for all $t,t^{\prime} \in \BC^{\iy}$}; \eea
the integral is taken over a small circle around
$z=\iy$. From the bilinear relation, one derives the
KP-hierarchy, already mentioned in Theorem 3.1, of
which the first equation reads as in (3.1.6).

 We consider the $p$-reduced KP hierarchy, i.e.,
the reduction to pseudo- differential $L$'s such that
$L^p=D^p+\ldots$ is a differential operator for some
fixed $p\geq 2$. Then $(L^{kp})_+=L^{kp}$ for all
$k\geq 1$ and thus $\pl L/\pl t_{kp}=0$, in view of
the deformation equations (5.0.1) on $L$. Therefore
the variables $t_p,t_{2p},t_{3p},...$ are not active
and can thus be set $=0$. The case $p=2$ is
particularly interesting and leads to the KdV
equation, upon setting all even $t_i=0$.

For the time being, take the integer $p\geq 2$
arbitrary. The arbitrary linear
combinations\footnote{$\zeta_p:=\{\om \mbox{ such that
} \om^p=1\}$}
\be
\Phi^{\pm}(x,t;z):=\sum_{\om\in\zeta_p}a^{\pm}_{\om}
\Psi^{\pm}(x,t;\om z)\ee are the most general solution
of the spectral problems $L^p\Phi^+=z^p\Phi^+$ and
$L^{\top p}\Phi^-=z^p\Phi^-$ respectively, leading to
the definition of the kernels:
  \be
k_{x,t}(y,z):=\int^xdx~\Phi^-(x,t;y)\Phi^+(x,t;z),\mbox{
and } k^E_{x,t}(y,z):=k_{x,t} (y,z)I_E(z), \ee
where the integral is taken from a fixed, but
arbitrary origin in $\BR$. In the same way that
$\Psi^{\pm}(x,t,z)$ has a $\tau$-function
representation, so also does $k^E_{x,t}(y,z)$ have a
similar representation, involving the vertex operator:
\be
Y(x,t;y,z):=\sum_{\om,\om'\in\zeta_p}a^-_{\om}a^+_{\om'}X(x,t;\om
y,\om'z), \ee where (see \cite{DJKM,ASV1,ASV2})\be
X(x,t;y,z):=\frac{1}{z-y}e^{(z-y)x+\sum_1^{\iy}(z^i-y^i)t_i}
e^{\sum_1^{\iy}(y^{-i}-z^{-i})\frac{1}{i}\frac{\pl}{\pl
t_i}}. \ee
 A condition
$
\sum_{\om\in\zeta_p}\frac{a^+_{\om}a^-_{\om}}{\om}=0
  $
is needed to guarantee that the right hand side of
(5.0.7) is free of singularities in the positive
quadrant $\{y_i\geq 0$ and $z_j\geq 0$ with
$i,j=1,...,n\}$ and $\lim_{y\rg z}Y(x,t;y,z)$ exists.
Indeed, using Fay identities and higher degree Fay
identities, one shows stepwise the following three
statements, the last one being a statement about a
Fredholm determinant\footnote{The Fredholm determinant
of a kernel $A(y,z)$ is defined by $$ \det(I- \lb
A)=1+\sum^{\infty}_{m=1} (-\lb)^m
\mathop{\int\!\ldots\!\int}\limits_{z_1\leq \cdots
\leq z_m} \det \left(A(z_i,z_j) \right)_{1\leq i,j\leq
m}dz_1\ldots dz_m . $$}: \bea
k_{x,t}(y,z)&=&\frac{1}{\tau(t)}Y(x,t;y,z)\tau(t)
\nonumber\\ \det\Bigl(k_{x,t}(y_i,z_j)\Bigr)_{1\leq
i,j\leq n}
&=&\frac{1}{\tau}\prod^k_{i=1}Y(x,t;y_i,z_i)\tau
\nonumber\\ \det(I-\lb
k_{x,t}^E)&=&\frac{1}{\tau}e^{-\lb\int_Edz\,Y(x,t;z,z)}
\tau=:\frac{\tau(t,E)}{\tau(t)}.
 \eea

 The kernel (5.0.12) at $t=0$ will define the ststistics of
 a random Hermitian ensemble, when the size
  $n \nearrow \iy$. The next
 theorem is precisely a statement about Fredholm
 determinants of
  kernels of the form (5.0.12); it will be identified at $t=0$ with
 the probability that no eigenvalue belongs to a
 subset E; see section 1.2. The
 initial condition that Virasoro annihilates $\tau_0$, as
 in sections 2.1.2 (Proof of Theorem 2.1), is now replaced by
 the {\em initial condition} (5.0.11) below.

\begin{theorem}  {\em (Adler-Shiota-van Moerbeke \cite{ASV1,
 ASV2})}
 Consider Virasoro generators $J^{(2)}_{\ell}$
 satisfying
\be
 \frac{\pl}{\pl z} z^{\ell
+1}  Y(x,t;z,z) =
\left[\frac{1}{2}J^{(2)}_{\ell}(t),Y(x,t;z,z)\right],
\ee where $Y(x,t;z,z)$ is defined in (5.0.7), and a
$\tau$-function satisfying the Virasoro constraint,
with an arbitrary constant $c_{kp}$:
  \be
\left(J_{kp}^{(2)}-c_{kp}\right)\tau =0\quad\mbox{for
a fixed $k\geq -1$. }
 \ee
 Then, given the disjoint union $E
\subset \BR^+$,
the Fredholm determinant of \be
K_{x,t}^E(\lb,\lb^{\prime}):=\frac{1}{p}
\frac{k_{x,t}(z,z^{\prime})}
{z^{\frac{p-1}{2}}z^{\prime\frac{p-1}{2}}}I_E(\lb^{\prime}),\quad\lb
=z^p,\lb^{\prime}=z^{\prime p}, \ee satisfies the
following constraint for that same $k\geq -1$:
  \be
\left(-\sum^{2r}_{i=1}c_i^{k+1}\frac{\pl}{\pl
c_i}+\frac{1}{2p}(J_{kp}^{(2)}-c_{kp})\right)\tau\det(I-\mu
 K_{x,t}^E)=0.\ee

\end{theorem}

\noindent The generators $J_n^{(2)}$ take on the
following precise form: {\footnotesize \bea
 J_n^{(1)}&:=&
\frac{\pl}{\pl t_n} + (-n)t_{-n}
 \nonumber\\
 J_n^{(2)}&:=&\sum_{i+j=n}
 \mathopen:J_i^{(1)}J_j^{(1)}\mathclose:-(n+1)J_n^{(1)}\nonumber\\
& =& \sum_{i+j=n}\frac{\pl^2}{\pl
 t_{i}\pl t_{j}}+2\sum_{-i+j=n}it_{i}\frac{\pl}{\pl
 t_{j}}+\sum_{-i-j=n}it_{i}jt_{j}-(n+1)J_n^{(1)}.\eea
  }

\remark For KdV (i.e., $p=2$),
$(L^2)^{\top}=L^2=D^2-q(x)$ holds, and thus the
adjoint wave function has the simple expression:
$\Psi^-(x,t;z)=\Psi^+(x,t;-z)$. In the next two
examples, which deal with KdV, set $$ \Psi(x,t;z):=
\Psi^+(x,t;z).$$

\subsubsection*{ Example 1: Eigenvalues of large random
 Hermitian matrices near the ``soft edge" and the Airy kernel}

Remember from section 1.2, the spectrum of the
Gaussian Hermitian matrix ensemble has, for large size
$n$, its edge at $\pm\sqrt{2n}$, near which the
scaling is given by $\sqrt{2}n^{1/6}$. Therefore, the
eigenvalues in Theorem 5.2 must be expressed in that
new scaling. Define the disjoint union $E=\cup_1^r
[c_{2i-1},c_{2i}]$, with $c_{2r}$ possibly $\iy$.

\begin{theorem} Given the spectrum $z_1\geq z_2 \geq
...$ of the large random Hermitian matrix $M$, define
the ``eigenvalues" in the new scale:
 \be
u_i=2n^{\frac{2}{3}}\left(\frac{z_i}{\sqrt{2n}}-1
\right)~~\mbox{for $n \nearrow \iy$}.
 \ee
The probability of the ``eigenvalues"
 \be  P(E^c):=
 P(~\mbox{all eigenvalues}~u_i \in E^c)\ee
 satisfies the partial differential equation
  $(\mbox{setting}~{\cal B}_k:=\sum_{i=1}^{2r}
c_i^{k+1} \frac{\pl}{\pl c_i})$\footnote{When
$c_{2r}=\iy$, that term in ${\cal B}_k$ is absent.}
\be
 \left({\cal B}_{-1}^3 - 4 ({\cal B}_0 -
\frac{1}{2})\right) {\cal B}_{-1}\log P(E^c) + 6({\cal
B}_{-1}^2 \log P(E^c))^2 = 0.
 \ee
In particular, the statistics of the largest
``eigenvalue" $u_1$ (in the new scale) is given by
   \be P(u_1\leq x) 
   =\exp\left(-\int^{\iy}_{x}(\al-x)g^2(\al)d\al
\right),\ee
 with
  \be \left\{\begin{array}{l}
g''=xg+2g^3 ~~
 (\mbox{\bf Painlev\'e II})
 \\ g(x)\cong
  -\frac{
  e^{-\frac{2}{3}  x^{\frac{3}{2}}}}{2\sqrt \pi x^{1/4}}
\mbox{\,\,for\,\,}x\nearrow \infty.
\end{array}\right..
 \ee
 \end{theorem}

  The partial differential equation (5.0.17) is due to
  Adler-Shiota-van Moerbeke \cite{ASV1,ASV2}.
   The equation (5.0.19) for the largest eigenvalue
   is a special case of (5.0.17), but was first derived by
   Tracy-Widom \cite{TW1}, by methods of
functional analysis.

\proof Remember from section 1.2, the statistics of
the eigenvalues is governed by the Fredholm
determinant of the kernel (1.2.4), for the Hermite
polynomials. In the limit, $$
\lim_{n\nearrow\infty}\frac{1}{\sqrt{2}n^{1/6}}K_n\left(\sqrt{2n}+\frac{u}{\sqrt{2}n^{1/
6}}, \sqrt{2n}+\frac{v}{\sqrt{2}n^{1/6}}
\right)=K(u,v), $$ where \be
K(u,v)=\int^{\iy}_{0}A(x+u)A(x+v)dx,~~~A(u)=\int^{\iy}_{
-\iy}e^{iux-x^3/3}dx. \ee
 Then
 \be P(E^c):=P( \mbox{ all
eigenvalue}~ u_i \in E^c)=\det (I-K(u,v)I_E(v)). \ee
In order to compute the PDE's of this expression, with
regard to the endpoints $c_i$ of the disjoint union
$E$, one proceeds as follows:

Consider the KdV wave function $\Psi(x,t;z)$, as in
(5.0.2),
 with initial condition:
\be
\Psi(x,t_0;z)=z^{\frac{1}{2}}A(x+z^2)=e^{xz+{2\over
3}z^{3}}(1+O(z^{-1})) ~ ,~z \rightarrow
\iy,~t_0=(0,0,\frac{2}{3},0,...),\ee in terms of the
Airy function\footnote{The $i$ in the definition of
the Airy function is omitted here. }, which, by
stationary phase, has the asymptotics:
 $$
A(u):=\frac{1}{\sqrt{\pi}}\int^{\infty}_{-\infty}
e^{-\frac{y^{3}}{3}+yu}dy
=u^{-\frac{1}{4}}e^{\frac{2}{3}u^{\frac{3}{2}}}(1+
O(u^{-\frac{3}{2}})). $$
 The definition of $A(u)$ is slightly changed,
 compared to (5.0.20). $A(u)$ satisfies the
 differential equation $A(y)^{\prime\prime}=yA(y)$,
 and thus
 the wave function $\Psi(x,t_0;z)$ satisfies
$(D^2-x)\Psi(x,t_0;z)=z^2\Psi(x,t_0;z)$. Therefore
$
L^2|_{t=t_0}=SD^2S^{-1}|_{t=t_0}=D^2-x,
$
so that $L^2$ is a differential operator, and $\Psi$
is a KdV wave function, with $\tau(t)$
satisfying\footnote{Although not used here, the
$\tau$-function is
 Kontsevich's integral:((see \cite{Kont,AMSV})
 $$
 \tau(t)=\frac{\int_\HR
dY\,e^{-\mbox{Tr}(\frac{1}{3}Y^3+Y^2Z)}}{\int_\HR
dY\,e^{-\mbox{Tr}Y^2Z}}\mbox{ with }
t_n=-\frac{1}{n}\Tr(Z^{-n}) + {2\over 3}\dt_{n,3},~
 Z=\mbox{diagonal matrix}. $$
  }
 the Virasoro constraints (5.0.11) with $c_{2k}=-\frac{1}{4}\delta_{k0}$.
The argument to prove these constraints is based on
the fact that the linear span (a point in an
infinite-dimensional Grassmannian) $${\cal
W}=\mbox{span} {}_{\BC}
\left\{ \psi_n(z):= 
e^{-\frac{2}{3} z^3} \sqrt{z} \left.\frac{\pl^n}{\pl
u^n}A(u)\right|_{u=z^2},~n=0,1,2,... \right\} $$ is
invariant under multiplication by $z^2$ and the
operator $\frac{1}{2z} (\frac{\pl}{\pl z}
+2z^2)-\frac{1}{4}z^{-2}$.

 Define for
$\lb=z^2, \lb'=z^{\prime 2}$, the kernel
$K_t(\lb,\lb')$,
\be
K_t(\lb,\lb')
:=\frac{1}{2z^{1\over2}z^{\prime{1\over2}}}\int^{\iy}_0
 \Psi(x,t;
z)\Psi(x,t;z')dx ,\ee
  which flows off the Airy kernel, by (5.0.22),
$$ K_{t_0}(\lb,\lb')=\frac{1}{2} \int^{\iy}_0
A(x+\lb)A(x+\lb')dx. $$
 Thus
 $\tau \det
(I-K_{x,t}^E)$ satisfies (5.0.13), with that same
constant $c_{kp}$, for $k=-1,0,1,\dots$\,:
\be
\left(-\sum_{i=1}^{2r}c_i^{k+1}\frac{\pl}{\pl
c_i}+\frac{1}{4}J_{pk}^{(2)} +{1\over 16
}\delta_{k,0}\right) \tau \det (I-K_{t}^E)=0. \ee

 Upon shifting $t_3 \mapsto t_3+2/3$, in view of
(5.0.22), the two first Virasoro constraints for
$k=-1$ and $k=0$ read: (${\cal B}_k:=\sum_{i=1}^{2r}
c_i^{k+1} \frac{\pl}{\pl c_i}$)
 {\footnotesize  \bea {\cal B}_{-1} \log
\tau(t,E)&=& \left(\frac{\pl}{\pl
t_1}+\frac{1}{2}\sum_{i\geq 3} it_i \frac{\pl}{\pl
t_{i-2}} \right)\log \tau(t,E) + \frac{t_1^2}{4}
\nonumber\\ {\cal B}_0 \log \tau(t,E)&=&
\left(\frac{\pl}{\pl t_3}+\frac{1}{2}\sum_{i\geq 1}
it_i \frac{\pl}{\pl t_{i}} \right) \log \tau(t,E) +
\frac{1}{16}. \eea  }
 The same method as in section 4 enables
one to express all the $t$-partials, appearing in the
KdV equation,
 $$
\left(\frac{\pl^4}{\pl t_1^4} -4\frac{\pl^2}{\pl t_1
\pl t_3}\right)\log\tau(t,E)+6\left(\frac{\pl^2}{\pl
t^2_1}\log\tau(t,E)\right)^2=0,
 $$
 in terms of $c$-partials, which upon substitution
 leads to the partial differential equation
$
 \left({\cal B}_{-1}^3 - 4 ({\cal B}_0 -
\frac{1}{2})\right) f + 6({\cal B}_{-1} f)^2 = 0
 $ (announced in (5.0.17))
 for
  $$ f: = {\cal B}_{-1} \log P(E^c) = \sum^{2r}_1
\frac{\partial }{\partial c_i} \log
P(E^c),~~\mbox{where}~ P(E^c)=\det(I-K^E)=\frac{\tau
(t,E)}{\tau (t)}.$$ When $E = (- \infty, x)$, this PDE
 reduces to an ODE:
\be
 f^{\prime\prime\prime} - 4 xf^{\prime} + 2f +
6{f^{\prime}}^2 = 0,~~~\mbox{with}~~~f =
\frac{d}{dx}\log P(\max_{i} ~
 \lb_i \leq x).
 \ee
 According to Appendix on
Chazy classes (section 9) , this equation can be
reduced to
 \be  f^{\prime\prime
2}+4f^{\prime}(f^{\prime 2 }-xf^{\prime}+f )=0,
\hspace{.8cm}(\mbox{\bf Painlev\'e II})\ee
 which can be solved by setting
 $$ f^{\prime}=-g^2~~\mbox{and}~~ f=g^{\prime
 2}-xg^2-g^4.
 $$
An easy computation shows $g$ satisfies the equation
 $
g^{\prime\prime} =2g^3+xg$ (Painlev\'e II), thus
leading to (5.0.19).


\subsubsection*{ Example 2: Eigenvalues of large random
 Laguerre Hermitian matrices near the ``hard edge" and
 the Bessel kernel}

 Consider the ensemble of $n\times n$ random matrices
for the Laguerre probability distribution, thus
corresponding to
 (1.1.9) with $\rho(dz)=z^{\nu/2}e^{-z/2}dz$. Remember
 from section 1.2,
 the density of eigenvalues near the ``hard edge"
 $z=0$ is given by $4n$ for very large $n$. At this
 edge, the kernel (1.2.4) with Laguerre polynomials $p_n$
 tends to the Bessel kernel
  \cite{Nagao, Forrester}:
 \be
  \lim_{n\nearrow\infty}\frac{1}{4n}K^{(\nu)}_n
   \left(\frac{u}{4n}, \frac{v}{4n}
\right)=K^{(\nu)}(u,v):= {1\over2}\int_0^1
xJ_\nu(xu)J_\nu(xv)dx. \ee

Therefore, the eigenvalues in the theorem below will
be expressed in that new scaling. Define, as before,
the disjoint union $E=\cup_1^r [c_{2i-1},c_{2i}]$.

\begin{theorem} Given the spectrum $0\leq z_1\leq z_2 \leq
...$ of the large random Laguerre-distributed
Hermitian matrix $M$, define the ``eigenvalues" in the
new scale:
 \be
u_i=4nz_i~~\mbox{for $n \nearrow \iy$}.
 \ee
The statistics of the ``eigenvalues"
 \be  P(E^c):=
 P(~\mbox{all ``eigenvalues"}~u_i \in E^c)\ee
 leads to the following PDE for $F=\log P(E^c)$: $(\mbox{setting}~{\cal B}_k:=\sum_{i=1}^{2r}
c_i^{k+1} \frac{\pl}{\pl c_i})$
 \be
  \left({\cal B}_0^4 -
2 {\cal B}_0^3 + ( 1 - \nu^2) {\cal B}_0^2 + {\cal
B}_1 \left({\cal B}_0 - \frac{1}{2}\right)\right) F -
4({\cal B}_0 F) ({\cal B}_0^2 F) + 6({\cal B}_0^2 F)^2
= 0.
 \ee
 In particular, for very large $n$, the statistics of
 the smallest eigenvalue is governed by
 $$ P(u_1 \geq x)=
 \exp \left(-\int_0^x \frac{f(u)}{u} du\right),  ~~~~
 u_1 \sim 4nz_1,
 $$
with $f$ satisfying
 \be
  (xf^{\prime\prime})^2-4\bigl(xf^{\prime}-f \bigr
)f^{\prime 2}+\bigl((x-\nu^2)f^{\prime}-f \bigr
)f^{\prime}=0 .\, \mbox{ \bf (Painlev\'{e} V)}
 \ee

 \end{theorem}

  The equation (5.0.32) for the smallest eigenvalue, first derived by
   Tracy-Widom \cite{TW1}, by methods of
functional analysis,
   is a special case of the partial differential
   equation (5.0.31), due to \cite{ASV1,ASV2}.

\remark This same theorem would hold for the Jacobi
 ensemble, near the ``hard edges" $z=\pm 1$.

\proof Define a wave function $\Psi(x,t;z)$, flowing
off
 $$ \Psi(x,0;z)=e^{xz}B(-xz)=e^{xz}(1+O(z^{-1}),
$$ where $B(z)$ is the Bessel function\footnote{
$\vr=i\sqrt{\pi/2}\,e^{i\pi\nu/2}$, $-1/2<\nu<1/2$.}
  $$
B(z)=\vr\sqrt{z}e^zH_\nu(iz)
={e^z2^{\nu+1/2}\over\Gamma(-\nu+1/2)}\int_1^\infty
{z^{-\nu+1/2}e^{-uz}\over(u^2-1)^{\nu+1/2}}du
=1+O(z^{-1}). $$
   As the operator $$
L^2|_{t=0}=D^2-{\nu^2-1/4\over x^2} $$ is a
differential operator, we are in the KdV situation;
again one may assume $t_2=t_4=\dots=0$ and we have
 $$ \Psi^-(x,t;-z)=
\Psi^+(x,t;z)=e^{xz+\sum
t_iz^i}{\tau(t-[z^{-1}])\over\tau(t)}, $$ in terms of
a $\tau$-function\footnote{ $\tau(t)$ is given by the
Adler-Morozov-Shiota-van Moerbeke double Laplace
matrix transform, with $t_n$ given in a similar way as
in footnote 32  (see \cite{AMSV}): $$
\tau(t)=c(t)\int_{\HR_N^+}dX\det X^{\nu-1/2}
e^{-\Tr(Z^2X)} \int_{\HR_N^+}dYS_0(Y)e^{-\Tr(XY^2)}~.
$$ }  satisfying the following Virasoro constraints
\be
J^{(2)}_{2k}\tau=((2\nu)^2-1)\delta_{k0}\tau.
\label{bessel} \ee
 Set $p=2$,
$a^-_1=a^+_{-1}={1\over4\pi}ie^{i\pi\nu/2}$ and
$a^-_{-1}=a^+_1={1\over4\pi}e^{-i\pi\nu/2}$ in
(5.0.5); this defines the kernel (5.0.6) and so
(5.0.12), which in terms of $\lb=z^2$ and
$\lb^{\prime}=z^{\prime 2}$, takes on the form:
\begin{eqnarray*}
K_{x,t}^{(\nu)}(\lb,\lb')&=&{1\over4\pi\sqrt{zz'}}\int^x
\left(ie^{i\pi\nu\over2}\Psi^*(x,t,z)+e^{-{i\pi\nu\over2}}\Psi^*(x,t,-z)
\right)\cdot\\
&&\hphantom{{1\over4\pi\sqrt{zz'}}\int_0^i}\cdot
\left(e^{-{i\pi\nu\over2}}\Psi(x,t,z') +
ie^{i\pi\nu\over2}\Psi(x,t,-z') \right) dx ,
\end{eqnarray*}
which flows off the {\bf Bessel kernel},
\begin{eqnarray*}
K_{x,0}^{(\nu)}(\lb,\lb')&=&
{1\over2}\int_0^xxJ_\nu(x\sqrt{\lb})J_\nu(x\sqrt{\lb'})dx.
\\
 &=& \frac{J_{\nu} (\sqrt{\lambda})
  \sqrt{\lambda^{\prime}}
  J'_{\nu} (\sqrt{\lambda^{\prime}}) -
J_{\nu} (\sqrt{\lambda^{\prime}}) \sqrt{\lambda}
J'_{\nu} (\sqrt{\lambda})} {2(\lambda -
\lambda^{\prime})}~~\mbox{for}~~ x=1.
\end{eqnarray*}
The Fredholm determinant satisfies for $E\subset
\BR_+$ and for $k=0,1,\dots$\,:
\be
\left(-\sum_1^{2r}c_i^{k+1}{\pl\over\pl c_i} +
{1\over4}J^{(2)}_{2k}
+\left({1\over4}-\nu^2\right)\dt_{k,0}\right)
\left(\tau\det(I-K_{x,t}^{(\nu)E}\right)=0. \ee

Upon shifting $t_1 \mapsto t_1+\sqrt{-1}$ and using
the same ${\cal B}_{i}$ as in (5.0.25), the equations
for $k=0$ and $k=1$ read
 \bea {\cal B}_{0} \log \tau(t,E) &=&
\frac{1}{2} \left(\sum_{i\geq 1} it_i \frac{\pl}{\pl
t_{i}} + \sqrt{-1}\frac{\pl}{\pl t_1} \right)\log
\tau(t,E) + \frac{1}{4}(\frac{1}{4}-\nu^2) \nonumber\\
{\cal B}_1\log \tau(t,E)&=& \frac{1}{2}
\left(\sum_{i\geq 1} it_i \frac{\pl}{\pl t_{i+2}}+
\frac{1}{2}\frac{\pl^2}{\pl
t_1^2}+\sqrt{-1}\frac{\pl}{\pl t_3} +
\frac{1}{2}\frac{\pl}{\pl t_1}  \right) \log \tau(t,E)
. \nonumber\\ \eea

 Expressing the $t$-partials (5.0.20), appearing in
the KdV-equation at $t=0$ (see formula below (5.0.25))
 in terms of the $c$-partials applied to $\log \tau(0,E)$,
 leads to the following PDE for $F=\log P(E^c)$:
 \be
  \left({\cal B}_0^4 -
2 {\cal B}_0^3 + ( 1 - \nu^2) {\cal B}_0^2 + {\cal
B}_1 \left({\cal B}_0 - \frac{1}{2}\right)\right) F -
4({\cal B}_0 F) ({\cal B}_0^2 F) + 6({\cal B}_0^2 F)^2
= 0.
 \ee
  Specializing this equation to the interval $E =
(0, x)$ leads to an ODE for $f: = -x~ \pl F / \pl x,$
namely
 \be f^{\prime\prime\prime} +\frac{1}{x}
f^{\prime\prime}-\frac{6}{x}f^{\prime 2}
 + \frac{4}{x^2}ff^{\prime} + \frac{(x - \nu^2)}{x^2}
f^{\prime} - \frac{1}{2x^2}f=0,
 \ee
  which is an
equation of the type (9.0.1); changing
$x\curvearrowright -x$ and $f\curvearrowright -f$
leads again to an equation of type (9.0.1), with
$P(x)=x, ~4Q(x)=-x-\nu^2$ and $R=0$. According to
Cosgrove-Scoufis \cite{CS} (see the Appendix on Chazy
classes), this equation can be reduced to the equation
(9.0.2), with the same $P,Q,R$ and with $c=0$. Since
$P(x)=x$, this equation is already in one of the
canonical forms (9.0.3), which upon changing back $x$
and $f$, leads to
$$(xf^{\prime\prime})^2+4\bigl(-xf^{\prime}+f \bigr
)f^{\prime 2}+\bigl((x-\nu^2)f^{\prime}-f \bigr
)f^{\prime}=0 .\, \mbox{ \bf (Painlev\'{e} V)}$$

\subsubsection*{ Example 3: Eigenvalues of large random
 Gaussian Hermitian matrices in the bulk and
 the sine kernel}

 Setting $\nu=\pm1/2$ yields kernels related to
the sine kernel:
  \bean
        K_{x,0}^{(+1/2)}(y^2,z^2)&=&{1\over\pi}\int_0^x
        {\sin xy~\sin xz
        \over y^{1/2}z^{ 1/2}}dx\nonumber\\
       & =&
       \frac{1}{2 \pi} \left(\frac{\sin x(y - z)}{y - z} -
         \frac{\sin x(y+ z)}{y + z}\right) \nonumber\\
         && \nonumber\\
        K_{x,0}^{(-1/2)}(y^2,z^2)&=&{1\over\pi}\int_0^x
        {\cos xy~\cos xz\over
        y^{1/2}z^{ 1/2}}dx\nonumber\\ &=&
        \frac{1}{2 \pi} \left(\frac{\sin x(y - z)}{y - z} + \frac{\sin x(y
+ z)}{y + z}\right).
 \eean

Therefore the sine-kernel obtained in the context of
the bulk-scaling limit (see 1.2.7) is the sum $
K_{x,0}^{(+1/2)} + K_{x,0}^{(-1/2)}$. Expressing the
Fredholm determinant of this sum in terms of the
Fredholm determinants of each of the parts, leads to
the Painlev\'e V equation (1.2.8).

\section{Coupled random Hermitian ensembles}

Consider a product ensemble $(M_1,M_2)\in {\cal H}_n^2
:={\cal H}_n\times{\cal H}_n$ of $n\times n$ Hermitean
matrices, equipped with a Gaussian probability
measure,
\begin{equation}
c_ndM_1dM_2\,e^{-\frac{1}{2}{\rm
Tr}(M^2_1+M_2^2-2cM_1M_2)},
\end{equation}
where $dM_1dM_2$ is Haar measure on the product ${\cal
H}_n^2$, with each $dM_i$,
\begin{equation}
dM_1=\Dt_n^2(x)\prod^n_1dx_idU\mbox{\,\,and\,\,}
dM_2=\Dt_n^2(y)\prod^n_1dy_idU \label{2}
\end{equation}
decomposed into radial and angular parts. In terms of
the coupling constant $c$, appearing in (6.0.1), and
the boundary of the set
\begin{equation}
 E=E_1\times E_2:=\cup^r_{i=1}[a_{2i-1},a_{2i}]\times
\cup^s_{i=1}[b_{2i-1},b_{2i}]\subset \BR^2,
 \end{equation}
define differential operators ${\cal A}_k$, $ {\cal
B}_k$ of ``weight" $k$,
$$
\begin{tabular}{ll}
${\cal
A}_1=\displaystyle{\frac{1}{c^2-1}\left(\sum^r_1\frac{\pl}{\pl
a_j}+c\sum^s_1\frac{\pl}{\pl b_j}\right)}$&${\cal
B}_1=\displaystyle{\frac{1}{1-c^2}\left(c\sum^r_1\frac{\pl}{\pl
a_j}+\sum^s_1\frac{\pl}{\pl b_j}\right)}$\\ ${\cal
A}_2=\displaystyle{\sum^r_{j=1}a_j^{}\frac{\pl}{\pl
a_j}-c\frac{\pl}{\pl c}}$&${\cal
B}_2=\displaystyle{\sum^s_{j=1} b_j^{}\frac{\pl}{\pl
b_j}-c\frac{\pl}{\pl c}},$
\end{tabular}
$$
forming a closed Lie algebra\footnote{
\begin{eqnarray*}
&& [{\cal A}_1,{\cal B}_1]=0\quad [{\cal A}_1,{\cal
A}_2]=\frac{1+c^2}{1-c^2}{\cal A}_1 \quad [{\cal
A}_2,{\cal B}_1]=\frac{2c}{1-c^2}{\cal A}_1 \\
&&[{\cal A}_2,{\cal B}_2]=0\quad [{\cal A}_1,{\cal
B}_2]=\frac{-2c}{1-c^2}{\cal B}_1 \quad\enspace [{\cal
B}_1,{\cal B}_2]=\frac{1+c^2}{1-c^2}{\cal B}_1    .
\nonumber
\end{eqnarray*}
}. The following theorem follows, via similar
 methods, from the Virasoro constraints (3.3.5) and
  the 2-Toda equation (3.3.6):


\begin{theorem} {\em ({\bf Gaussian probability}) (Adler-van Moerbeke \cite{AvM2})}
  The joint statistics
  \bean
  P_n(M \in   {\cal H}_n^2(E_1\times
  E_2))
 & =&  \frac{\displaystyle{\int\!\!\int_{{\cal H}_n^2(E_1\times
  E_2)}}dM_1dM_2\,e^{-\frac{1}{2}{\rm
Tr}(M^2_1+M_2^2-2cM_1M_2)}} {\displaystyle{\int \!\!
\int_{{\cal H}_n^2}} dM_1dM_2\,e^{-\frac{1}{2}{\rm
Tr}(M^2_1+M_2^2-2cM_1M_2)}}\\ &=&\frac{
\displaystyle{\int\!\!\int_{E^n}}
\Delta_n(x)\Delta_n(y) \prod^n_{k=1}e^{-\frac{1}{2}(
x_k^2+y_k^2-2cx_ky_k) }dx_k dy_k} {
\displaystyle{\int\!\!\int_{\BR^{2n}}}
\Delta_n(x)\Delta_n(y) \prod^n_{k=1}e^{-\frac{1}{2}(
x_k^2+y_k^2-2cx_ky_k) }dx_k dy_k}
 \eean
  satisfies the non-linear third-order partial
  differential equation\footnote{in terms of the Wronskian
  $\{f,g\}_X=Xf.g-f.Xg$, with regard to a first order
  differential operator $X$.} (independent of $n$) $(~F_n: =\frac{1}{n}\log P_n(E)~)$:
\be
 \left\{{\cal B}_2   {\cal A}_1 F_n  ~,~
  {\cal B}_1 {\cal A}_1 F_n +\frac{c}{c^2-1} \right\}
_{ {\cal A}_1} ~-~ \left\{{\cal A}_2   {\cal B}_1 F_n
~,~ {\cal A}_1 {\cal B}_1 F_n+\frac{c}{c^2-1}\right\}
_{ {\cal B}_1 }=0
 .
 \ee
\end{theorem}


\section{Random permutations }

The purpose of this section is to show that the
generating function of the probability
 $$
 P(L(\pi_n)\leq\ell)= \frac{1}{n!} \# \{\pi_n \in S_n
 ~\bigl|~ L(\pi_n)\leq \ell \}
 $$
is closely related to a special solution of the
Painlev\'e V equation, with peculiar initial
condition. Remember from section 1.4, $L(\pi_n)$ is
the length of the longest increasing sequence in the
permutation $\pi_n$.

\begin{theorem} For every $\ell\geq 0$,
  \bea
&&\sum^{\iy}_{n=0}\frac{x^n}{n!}P(L(\pi_n)\leq\ell)
 = \int_{U(\ell )}e^{\sqrt{x}~tr
(M+\bar M)}dM=\exp {\int_0^x \log \left(\frac{x}{u}
\right)g_{\ell}(u)du},
 \nonumber \\
\eea
  with $  g_{\ell}$
satisfying the initial value problem for {\bf
Painlev\'e V}:
\be
\left\{
\begin{array}{l}
\displaystyle{g^{\prime\prime}-
 \frac{g^{\prime 2}}{2}\left(\frac{1}{g-1}+
  \frac{1}{g}\right)+ \frac{g^{\prime}}{u}
   +\frac{2}{u}g(g-1)-
   \frac{\ell^2}{2u^{2}}\frac{g-1}{g}=0}
 \\  \\
\displaystyle{g_{\ell}(u)=1-\frac{u^{\ell}}{(\ell!)^2}
 +O(u^{\ell+1}),~~\mbox{near} ~u=0
  .}
\end{array}
\right. \ee

\end{theorem}

 The Painlev\'e V equation for this integral was first found
 by Tracy-Widom \cite{TW3}. This systematic
 derivation and the initial condition are due to
 Adler-van Moerbeke \cite{AvM4}.

\proof Upon inserting $(t_1,t_2,...)$ and
$(s_1,s_2,...)$ variables in the $U(n)$-integral
(7.0.1), the integral
 \bea
 I_n(t,s)&=&\int_{U(n)}e^{Tr \sum_1^{\iy}(t_iM^i-s_i\bar
 M^i)}dM \\ &=& n!
   \det\left(\int_{S^1}z^{k-\ell}
e^{\sum_1^{\iy}(t_i z^i-s_iz^{-i})}
 \frac{dz}{2\pi i z}\right)_{0\leq k,\ell\leq
 n-1}=n!\tau_n(t,s)  \nonumber
 \eea
puts us in the conditions of Theorem 3.5. It deals
with semi-infinite matrices $L_1$ and
$hL^{\top}_2h^{-1}$ of ``rank $2$", having diagonal
elements
 $$
 b_n:=\frac{\pl}{\pl t_1}\log
 \frac{\tau_n}{\tau_{n-1}}=(L_1)_{n-1,n-1}
 ~~~\mbox{and}~~~
 b^*_n:=-\frac{\pl}{\pl s_1}\log \frac{\tau_n}{\tau_{n-1}}
 =(hL_2^{\top}h^{-1})_{n-1,n-1}.
 $$ To summarize Theorem 3.5,
$I_n(t,s)$ satisfies the following three types of
identities:

\medbreak

\noindent (i) {\bf Virasoro}: ($F:=\log \tau_n$) (see
(3.4.7)) {\footnotesize\begin{eqnarray} 0=\frac{{\cal
V}_{-1}\tau_n}{\tau_n}&=&\left(\sum_{i\geq
1}(i+1)t_{i+1}\frac{\pl}{\pl t_{i}}-\sum_{i\geq
2}(i-1)s_{i-1}\frac{\pl}{\pl s_{i}}+n\frac{\pl}{\pl
s_{1}}\right)F+nt_1\nonumber\\0=\frac{ {\cal
V}_{0}\tau_n}{\tau_n}&=&\sum_{i\geq
1}\left(it_{i}\frac{\pl}{\pl
t_{i}}-is_{i}\frac{\pl}{\pl s_{i}}\right)F\nonumber\\
0=\frac{{\cal
V}_{1}\tau_n}{\tau_n}&=&\left(-\sum_{i\geq
1}(i+1)s_{i+1}\frac{\pl}{\pl s_{i}}+\sum_{i\geq
2}(i-1)t_{i-1}\frac{\pl}{\pl t_{i}}\frac{\pl}{\pl t_1}
\right)F+ns_1.\nonumber\\ 0= \frac{\pl}{\pl
t_{1}}\frac{{\cal V}_{-1}\tau_{n}}{\tau_{n}}
 &=&\left(\sum_{i\geq
1}(i+1)t_{i+1}\frac{\pl^{2}}{\pl t_{1}\pl t_{i}}-
\sum_{i\geq 2} (i-1)s_{i-1}\frac{\pl^2}{\pl t_{1}\pl
s_{i}}+n\frac{\pl^2}{\pl t_1\pl
s_1}\right)F+n\nonumber\\
\end{eqnarray}
}
  \noindent (ii) {\bf two-Toda}: (see (3.4.8))
\begin{eqnarray}
\frac{\pl^2\log\tau_{n}}{\pl s_2\pl
t_1}&=&-2\frac{\pl}{\pl
s_1}\log\frac{\tau_{n}}{\tau_{n-1}}\frac{\pl^{2}}{\pl
s_{1}\pl t_{1}}\log\tau_{n}- \frac{\pl^{3}}{\pl
s_{1}^2\pl t_{1}}\log\tau_{n} \nonumber\\
 &=& 2b_n^* ~\frac{\pl^{2}}{\pl s_{1}\pl
t_{1}}\log\tau_{n}- \frac{\pl^{3}}{\pl s_{1}^2\pl
t_{1}}\log\tau_{n},
\end{eqnarray}

\noindent (iii) {\bf Toeplitz}: (see (3.4.11))
 {\footnotesize
\begin{eqnarray} {\cal T}(\tau)_{n}&=&\frac{\pl}{\pl
t_{1}}\log \frac{\tau_n}{\tau_{n-1}} \frac{\pl}{\pl
s_{1}}\log\frac{\tau_n}{\tau_{n-1}}\nonumber\\
 && \hspace{0cm}+
 \left(1+\frac{\pl^{2}}{\pl s_{1}\pl t_{1}}
  \log\tau_n\right)\left(1+\frac{\pl^{2}}{\pl s_{1}\pl
t_{1}}\log\tau_n-\frac{\pl}{\pl s_{1}}
\left(\frac{\pl}{\pl t_{1}}\log\frac{\tau_{n}}
{\tau_{n-1}}\right)\right)
 \nonumber\\
 &=&-b_n b^*_n+\left( 1+\frac{\pl^{2}}{\pl s_{1}\pl t_{1}}\log\tau_n
   \right)  \left(1+ \frac{\pl^{2}}{\pl s_{1}\pl t_{1}}\log\tau_n
   -\frac{\pl}{\pl s_1}b_n  \right)=0. \nonumber\\
\end{eqnarray}  }
 Defining the locus ${\cal
L}=\{~\mbox{all}~t_{i}=s_{i}=0$, except
$t_{1},s_{1}\neq 0\}$, and using the second relation
(7.0.4), we have on ${\cal L}$,
 $$
   \frac{{\cal
V}_{0}\tau_{n}}{\tau_{n}}\Big|_{\LR}
=\left(t_{1}\frac{\pl}{\pl t_{1}}-s_{1}\frac{\pl}{\pl
s_{1}}\right)\log\tau_{n}\Big|_{\LR}=0, $$ implying
$\tau_{n}(t,s)\Big|_{\LR}$ is a function of
$x:=-t_{1}s_{1}$ only. Therefore we may write
$\tau_{n}\Big|_{\LR}=\tau_{n}(x)$, and so, along
$\LR$, we have  $ \frac{\pl}{\pl
t_{1}}=-s_{1}\frac{\pl}{\pl x},\quad\frac{\pl}{\pl
s_{1}}=-t_{1}\frac{\pl}{\pl x},\quad\frac{\pl^2}{\pl
t_{1}\pl s_{1}}=-\frac{\pl}{\pl x}x\frac{\pl}{\pl x}.
$ Setting \be
f_{n}(x)= \frac{\pl}{\pl x}x\frac{\pl}{\pl
x}\log\tau_{n}(x)=- \frac{\pl^2}{\pl t_{1}\pl
s_{1}}\log\tau_{n}(t,s)\Big|_{\LR}, \ee and using
$x=-t_1s_1$, the two-Toda relation (7.0.5) takes on
the form \bean s_1\frac{\pl^2\log\tau_n}{\pl s_2\pl
t_1}\Big|_{\LR}
 &=&s_1\left(2b_n^* ~\frac{\pl^2}{\pl s_1\pl t_1}
  \log\tau_n-
\frac{\pl}{\pl s_1}\left(\frac{\pl^2\log\tau_n} {\pl
s_1\pl t_1}\right)\right)\nonumber\\
&=&x(2\frac{b_n^*}{t_1}f_n+f^{\prime}_n). \eean
Setting relation (7.0.7) into the Virasoro
 relations (7.0.4) yields
\begin{eqnarray*}
0=\frac{{\cal V}_{0}\tau_{n}}{\tau_{n}} -\frac{{\cal
V}_{0}\tau_{n-1}}{\tau_{n-1}}\Big|_{\LR}
&=&\left(t_{1}\frac{\pl}{\pl
t_{1}}-s_{1}\frac{\pl}{\pl
s_{1}}\right)\log\frac{\tau_{n}}{\tau_{n-1}}\Big|_{\LR}=
 t_1b_n+s_1b_n^* \nonumber\\ &&\\
0=  \frac{\pl}{\pl t_{1}}\frac{{\cal
V}_{-1}\tau_{n}}{\tau_{n}}\Big|_{\LR} &=&
\left(-s_{1}\frac{\pl^2}{\pl s_{2}\pl t_{1}}+n
\frac{\pl^2}{\pl t_{1}\pl
s_{1}}\right)\log\tau_{n}\Big|_{\LR}+n \nonumber
\\ & &\nonumber\\
 &=& -x\left( 2\frac{b_n^*}{t_1} f_n(x)+f'_n(x)
   \right) +n (-f_n(x)+1).\nonumber
\end{eqnarray*}
This is a system of two linear relations in $b_n$ and
$b_n^*$, whose solution, together with its
derivatives, are given by: $$
\frac{b_n^*}{t_1}=-\frac{b_n}{s_1}=-
 \frac{n(f_n-1)+xf^{\prime}_n}{2xf_n},\quad
\frac{\pl b_n}{\pl s_1}=\frac{\pl }{\pl
x}x\frac{b_n}{s_1}=\frac{x(f_n
f^{\prime\prime}_n-f^{\prime
2}_n)+(f_n+n)f^{\prime}_n}{2f^2_n}. $$ Substituting
the result into the Toeplitz relation (7.0.6), namely
 $$ b_nb_n^* =(1-f_n)\left(1-f_n-\frac{\pl}{\pl
s_1}b_n\right), $$ leads to $f_n$ satisfying
Painlev\'e equation (7.0.2), with $g=f_n$, as in
(7.0.7) and $u=x$.
 Note, along the locus
$\LR$, we may set $t_{1}=\sqrt{x}$ and
$s_{1}=-\sqrt{x}$, since it respects $t,s_{1}=-x$.
Thus,  $\left.I_{n}(t,s)\right|_{\LR}$ equals (7.0.1).

The initial condition (7.0.2) follows from the fact
that as long as $0\leq n \leq \ell$, the inequality
$L(\pi_n)\leq \ell$ is always verified, and so
 \bean
\sum_0^{\iy} \frac{x^n}{(n!)^2}\#\{ \pi \in
S_n~\bigl|~L(\pi_n)\leq \ell\}
  &=& \sum_0^{\ell} \frac{x^n}{n!}+
  \frac{x^{\ell+1}}{(\ell+1)!^2}
   ((\ell+1)!-1)  +O(x^{\ell+2})\\
  &=&\exp\left(x-\frac{x^{\ell+1}}{(\ell+1)!^2}
  +O(x^{\ell+2})\right),
 \eean
thus proving Proposition 7.1. \qed

 \remark Setting $$ f_n(x)=\frac{g(x)}{g(x)-1} $$ leads to
standard Painlev\'e V, with $\alpha=\delta=0,~
\beta=-n^2/2, \gamma=-2$.

\section{Random involutions }

This section deals with a generating function for the
distribution of the length of the longest increasing
sequence of a fixed-point free random involution
$\pi_{2k}^0$, with the uniform distribution:
 $$ P\left( L(\pi_{2k}^{0})\leq
\ell+1,~ \pi_{2k}^{0}\in S^{0}_{2k} \right)
=\frac{2^kk!}{(2k)!}
 \#\{\pi_{2k}^0\in
S^0_{2k}~\bigl|~L(\pi_{2k}^0)\leq \ell+1\}. $$

\begin{proposition} {\em (Adler-van Moerbeke \cite{AvM4})}
The generating function
 \bea
\lefteqn{2\sum^{\iy}_{k=0}\frac{(x^2/2)^k}{k!}
 P(L(\pi^0_{2k})  \leq \ell +1  )
 }\nonumber\\& =& E_{O(\ell +1)_-} e^{x Tr M} + E_{O(\ell +1)_+}
 e^{x Tr M}  \nonumber\\
&=&\exp \left({\int_0^x\frac{f^-_{\ell
 }(u)}{u}du}\right)+\exp \left({\int_0^x\frac{f^+_{\ell
}(u)}{u}du}\right), 
 \eea
 where $f=f^{\pm}_{\ell}$, satisfies the initial value
problem:
 \be
 \left\{\begin{array}{l}
\displaystyle{f^{\prime\prime\prime}+
 \frac{1}{u}
f^{\prime\prime}
 +\frac{6}{u}{f^{\prime} }^2-\frac{4}{u^2} f f^{\prime}
 -\frac{16 u^2+\ell^2}{u^2}f^{\prime}
+\frac{16}{u}f  +\frac{2(\ell^2-1)}{u}=0} \\
   \\
\displaystyle{ f_{\ell}^{\pm}(u)=u^2 \pm
\frac{u^{\ell+1}}{\ell !}+O(u^{\ell+2}),~\mbox{near}~
u=0 .}
 ~~~~~~~~~~~~~~~~~~~~ (\mbox{\bf Painlev\'e V})
\end{array}
 \right. \ee

\end{proposition}

\proof The first equality in (8.0.1) follows
immediately from proposition 1.1. The results of
section 1.3 lead to
 \be
 \int_{O(2n+1)_{\pm}}e^{x Tr M }dM   =   e^{\pm x}
 \int_{[-1,1]^n}\Delta_n(z)^2
 \prod_{k=1}^n e^{2xz_k} (1-z_k)^{a}(1+z_k)^{b}
 dz_k,
   \ee
   with $ a=\pm 1/2,b=\mp 1/2$, (with corresponding
signs). Inserting $t_i$'s in the integral, the
perturbed integral, with $e^{\pm x}$ removed and with
$t_1=2x$, reads
 \be
I_n(t)= \int_{[-1,1]^n}\Dt_n(z)^{2}\prod_{k=1}^n
(1-z_k)^{a}(1+z_k)^{b}e^{\sum_{1}^{\iy}t_iz_k^i}dz_k=n!\tau_n(t);
 \ee
this is precisely integral (3.1.4) of section 3.1.1
and thus it satisfies the Virasoro constraints
(3.1.5),
 but without boundary
contribution ${\cal B}_iF$. Explicit Virasoro
expressions appear in (2.1.35), upon setting $\beta
=2$.  Also, $\tau_n(t)$, as in (3.1.4), (see Theorem
3.1) satisfies the KP equation (3.1.6).
Differentiating the Virasoro constraints in $t_1$ and
$t_2$, and restricting to the locus
 $$ \LR:=\{t_1=x, ~\mbox{all other}~t_i=0 \},$$ lead to a
linear system of five equations, with
$b_0=a-b,~b_1=a+b$,
  \bean
  && \frac{1}{I_n}\left.\left(\BJ_{k+2}^{(2)}-\BJ_{k}^{(2)}+
b_0\BJ_{k+1}^{(1)}
+b_1\BJ_{k+2}^{(1)}\right)I_n\right|_{\LR}=0,~~
k=-1,0\\
 &&\frac{\pl}{\pl t_1} \frac{1}{I_n}\left.\left(\BJ_{k+2}^{(2)}-\BJ_{k}^{(2)}+
b_0\BJ_{k+1}^{(1)}
+b_1\BJ_{k+2}^{(1)}\right)I_n\right|_{\LR}=0,~~ k=-1,0
\\
 &&\frac{\pl}{\pl t_2} \frac{1}{I_n}\left.\left(\BJ_{k+2}^{(2)}-\BJ_{k}^{(2)}+
b_0\BJ_{k+1}^{(1)}
+b_1\BJ_{k+2}^{(1)}\right)I_n\right|_{\LR}=0,~~ k=-1
 \eean
 in five unknowns ($F_n=\log \tau_n$) $$ \left.\frac{\pl F_n}{\pl
t_2}\right|_{\LR},
 \quad\left.\frac{\pl F_n}{\pl t_3}\right|_{\LR},\quad\left.
 \frac{\pl^2F_n}{\pl
t_1\pl
t_2}\right|_{\LR},\quad\left.\frac{\pl^2F_n}{\pl
t_1\pl t_3}\right|_{\LR},\quad
\left.\frac{\pl^2F_n}{\pl t_2^2}\right|_{\LR}.\quad $$
  Setting $t_1=x$
and $F'_n=\pl F_n/ \pl x$, the solution is given by
the following expressions,
 {\footnotesize \bean
\left.\frac{\pl F_n}{\pl
t_2}\right|_{\LR}&=&-\frac{1}{x}
\Bigl((2n+b_1)F^{\prime}_n+n(b_0-x)\Bigr)\\ & &
\\ \left.\frac{\pl F_n}{\pl
t_3}\right|_{\LR}&=&-\frac{1}{x^2}\Bigl(x\left(F^{\prime\prime}_n+
F_n^{\prime
2}+(b_0-x)F^{\prime}_n+n(n+b_1)\right)-(2n+b_1)
\left((2n+b_1)F^{\prime}_n+b_0n\right)\Bigr)\\ & & \\
\left.\frac{\pl^2 F_n}{\pl t_1\pl
t_2}\right|_{\LR}&=&-\frac{1}{x^2}
 \Bigl((2n+b_1)(xF^{\prime\prime}_n-F^{\prime}_n)-b_0n\Bigr)
\\
\left.\frac{\pl^2 F_n}{\pl t_1\pl t_3}\right|_{\LR}
 &=&-\frac{1}{x^3}\Bigl(
x^2(F^{\prime\prime\prime}_n+2
F_n^{\prime}F_n^{\prime\prime})-x\left((x^2-b_0x+1)
F_n^{\prime\prime}+F_n^{\prime
2}+b_0F_n^{\prime}+(2n+b_1)^2F_n^{\prime\prime}\right.\\
 &&\hspace{3cm} \left.
+n(n+b_1)\right)+2(2n+b_1)^2F_n^{\prime}+2b_0n(2n+b_1)\Bigr)\\
& & \\ \left. \frac{\pl^2 F_n}{\pl t_2^2}\right|_{\LR}
 &=&\frac{1}{x^3}\Bigl(
x\left(2F^{\prime 2}_n+2b_0
F^{\prime}_n+((2n+b_1)^2+2)F^{\prime\prime}_n+2n(n+b_1)\right)\\
 & &\hspace{4cm} -3(2n+b_1)^2F^{\prime}_n
  -3b_0n(2n+b_1)\Bigr).
\eean  }

\vspace{-.6cm}

\noindent Putting these expressions into KP and
setting $t_1=x$, one finds:
 \bean
  0&=& \left(\left(\frac{\pl}{\pl
t_1} \right)^4+3\left(\frac{\pl}{\pl
t_2}\right)^2-4\frac{\pl^2}{\pl t_1 \pl
t_3}\right)F_n+6\left(\frac{\pl^2}{\pl t^2_1}F_n
\right)^2\\
 &=&\frac{1}{x^3}\Bigl(
 ~x^3F^{\prime\prime\prime\prime}+4x^2F^{\prime\prime\prime}
+x\left(-4x^2+4b_0x+2-(2n+b_1)^2\right)F^{\prime\prime}
+8x^2F^{\prime}F^{\prime\prime}\nonumber\\&&+6x^3
 {F^{\prime\prime}}^2
+2x{F^{\prime}}^2+\left(2b_0x-(2n+b_1)^2\right)F^{\prime}+n(2x-b_0)(n+b_1)-b_0n^2
\Bigr). \nonumber\\
 \eean

Finally, the function
$H(x):=x\frac{d}{dx}F(x)=x\frac{d}{dx} \log \tau_n(x)$
satisfies \bea
&&x^2H^{\prime\prime\prime}+xH^{\prime\prime}
 +6x{H^{\prime} }^2-\left(4H+4x^2-4bx
  +(2n+a)^2\right)H^{\prime}+(4x-2b)H\nonumber\\
&&\hspace{7cm}+2n(n+a)x-bn(2n+a)=0.\nonumber\\ \eea
This 3rd order equation is Cosgrove's
\cite{CS,Cosgrove} equation, with
$P=x,~4Q=-4x^2+4bx-(2n+a)^2,~2R=2n(n+a)x-bn(2n+a)$.
So, this 3rd order equation can be transformed into
the Painlev\'e V equation (9.0.3) in the appendix. The
boundary condition $f(0)=0$ follows from the
definition of $H$ above, whereas, after an elementary,
but tedious computation,
$f^{\prime}(0)=f^{\prime\prime}(0)=0$ follows from the
differential equation (8.0.5) and the Aomoto extension
\cite{A} (see Mehta \cite{Mehta}, p. 340) of Selberg's
integral:\footnote{where $Re\,\gamma$, $Re\,\delta
>-1$, $Re\,\beta >
-2\min\displaystyle{\left(\frac{1}{n},
\frac{Re\,\gamma +1}{n-1}, \frac{Re\,\delta
+1}{n-1}\right)}$}

{\footnotesize $$\displaystyle{\frac{\int^1_{0}\ldots
\int^1_{0}x_1\ldots
x_m\left|\Delta(x)\right|^{\beta}\prod^n_{j=1}
 x^{\gamma}_{j}
(1-x_{j})^{\delta}dx_1...dx_n}{\int^1_{0}\ldots
\int^1_{0}\left|\Delta(x)\right|^{\beta}
\prod^n_{j=1}x^{\gamma}_{j}
(1-x_{j})^{\delta}dx_1...dx_n}}=\prod^m_{j=1}\frac{\gamma
+1+(n-j)\beta/2}{\gamma+\delta+2+(2n-j-1)\beta/2}. $$}
 However, the initial condition (8.0.2) is a much
stronger statement, again stemming from the fact that
as long as $0\leq n \leq \ell$, the inequality
$L(\pi_n)\leq \ell$ is trivially verified, thus
leading to
 $$E_{O_{\pm}(\ell  +1)} e^{x Tr M}  =
\exp \left(
 \frac{x^2}{2}\pm\frac{x^{\ell+1}}{(\ell+1)!}
 +O(x^{\ell+2})\right),$$
ending the proof of Proposition 8.1.\qed


\section{Appendix: Chazy classes}

Most of the differential equations encountered in this
survey belong to the general Chazy class $$
f^{\prime\prime\prime}=F(z,f,f^{\prime},f^{\prime\prime}),~\mbox{where
$F$ is rational in $f,f^{\prime},f^{\prime\prime}$ and
locally analytic in z,}$$ subjected to the requirement
that the general solution be free of movable branch
points; the latter is a branch point whose location
depends on the integration constants. In his
classification Chazy found thirteen cases, the first
of which is given by
\be
f^{\prime \prime\prime}+\frac{P'}{P}f^{
\prime\prime}+\frac{6}{P}f^{\prime
2}-\frac{4P'}{P^2}ff' +\frac{P^{\prime\prime}}{P^2}
f^2
+\frac{4Q}{P^2}f'-\frac{2Q'}{P^2}f+\frac{2R}{P^2}=0
\ee
 with arbitrary polynomials $P(z), Q(z), R(z)$ of
degree $3,2,1$ respectively. Cosgrove and Scoufis
\cite{CS,Cosgrove}, (A.3), show that this third order
equation
 has a first integral, which is second order in $f$
and quadratic in $f^{\prime\prime}$, \bea &f^{
\prime\prime 2}& +\frac{4}{P^2}
 \left( (Pf^{\prime 2}+Q f^{\prime}+R)f^{\prime}
 - (P' f^{\prime 2}+\frac{}{}Q' f^{\prime}+R')f^{}
  \right. \nonumber\\ && \hspace{2cm} \left.
   +\frac{1}{2}(P^{\prime\prime}f^{\prime
}+Q^{\prime\prime} )f^2
 -\frac{1}{6} P^{\prime\prime\prime}f^3 +c\right)=0;
\eea $c$ is the integration constant. Equations of the
general form $$ f^{ \prime\prime 2}=G(x,f,f^{
\prime})$$ are invariant under the map $$ x\mapsto
\frac{a_1z+a_2}{a_3z+a_4}~~ \mbox{and}~~ f\mapsto
\frac{a_5f+a_6z+a_7}{a_3z+a_4}.$$ Using this map, the
polynomial $P(z)$ can be normalized to $$
P(z)=z(z-1),~z,~\mbox{or} ~1.$$

 In this way, Cosgrove shows (9.0.2) is a master Painlev\'e
 equation, containing the 6 Painlev\'e equations.  In
each of the cases, the canonical equations are
respectively:

\begin{itemize}

 \item $  g^{\prime\prime 2}=-4g^{\prime 3
}-2  g^{\prime}(zg^{\prime}-g)+A_1 $ \hspace{4cm}({\bf
Painlev\'e II})

  \item $  g^{\prime\prime 2}=-4g^{\prime 3
}+4(zg^{\prime}-g)^2+A_1g^{\prime}+A_2$
\hspace{3cm}({\bf Painlev\'e IV})

  \item $ (zg^{\prime\prime })^2=(zg^{\prime}-g)\Bigl(-4g^{\prime 2
}+A_1(zg^{\prime}-g)+A_2\Bigr)+A_3g^{\prime}+A_4
 $

 \hspace{10cm}({\bf
Painlev\'e V})
  \item $ (z(z-1)g^{\prime\prime
})^2=(zg^{\prime}-g)\Bigl(4g^{\prime
2}-4g^{\prime}(zg^{\prime}-g)
 +A_2\Bigr) +A_1g^{\prime
2}+A_3g^{\prime}+A_4 $

\hspace{10cm}({\bf Painlev\'e VI})
\end{itemize}
\vspace{-1cm}\be \ee
 Painlev\'e II equation above can
be solved by setting
 \bean
   g(z)&=&\frac{1}{2} (u^{\prime})^2-\frac{1}{2}
(u^2+\frac{z}{2})^2-(\alpha+\frac{\vr_1}{2})u\\
g^{\prime}(z)&=&-\frac{\vr_1}{2}
u^{\prime}-\frac{1}{2} (u^2+\frac{z}{2})\\ A_1&=&
\frac{1}{4}(\al +(u^2+\frac{z}{2})^2
\vr_1)^2,~~(\vr=\pm 1).
 \eean
Then $u(z)$ satisfies yet another version of the
Painlev\'e II equation
 $$
 u^{\prime\prime}=2u^3+zu+\al. \hspace{4cm}\mbox{({\bf
Painlev\'e II})}
  $$
 Now, each of these Painlev\'e II,IV,V,VI equations can be transformed
into the standard Painlev\'e equations, which are all
differential equations of the form $$ f^{\prime\prime
}=F(z,f,f^{\prime}), \mbox{rational in $f,~f^{\prime
}$, analytic in $z$,} $$ whose general solution has no
movable critical points. Painlev\'e showed that this
requirement leads to 50 types of equations, six of
which cannot be reduced to known equations.

\end{document}